\magnification=1000
\hsize=11.7cm
\vsize=18.9cm
\lineskip2pt \lineskiplimit2pt
\nopagenumbers

\hoffset=-1truein
\voffset=-1truein

\advance\voffset by 4truecm
\advance\hoffset by 4.5truecm

\newif\ifentete

\headline{\ifentete\ifodd	\count0 
      \rlap{\head}\hfill\tenrm\llap{\the\count0}\relax
    \else
        \tenrm\rlap{\the\count0}\hfill\llap{\head} \relax
    \fi\else
\global\entetetrue\fi}

\def\entete#1{\entetefalse\gdef\head{#1}}
\entete{}

\input amssym.def
\input amssym.tex

\def\-{\hbox{-}}
\def\.{{\cdot}}

\def\P{{\cal P}}

\def\R{{\cal R}}

\def\B{{\cal B}}

\def\D{{\cal D}}

\def\Gr{\frak G\frak r}

\def\int{\frak i\frak n\frak t}

\def\qq{\quad{\rm and}\quad}

\def\too{\longrightarrow}

 3
 2
\font\large=cmr10  scaled \magstep 2
 2
\font\larti=cmti10  scaled \magstep 2
 2
\font\cds=cmr7

\font\it=cmti10

\centerline{\large Weight parameterization of simple modules}
\smallskip
\centerline{\large  for {\larti p}-solvable groups}

\vskip 0.5cm

\centerline{\bf Lluis Puig }
\medskip
\centerline{\cds CNRS, Institut de Math\'ematiques de Jussieu}
\smallskip
\centerline{\cds 6 Av Bizet, 94340 Joinville-le-Pont, France}
\smallskip
\centerline{\cds puig@math.jussieu.fr}

\vskip 0.5cm
\noindent
{\bf £1. Introduction}
\bigskip
£1.1. The {\it weights\/} for a finite group $G$ with respect to a prime number $p$ where introduced by Jon Alperin in [1] in order to formulate his celebrated conjecture. Explicitly, a 
{\it weight\/} of $G$  is a pair $(R,Y)$ formed by a $p\-$subgroup $R$ of $G$ and by an isomorphism class $Y$ of simple $kN_G (R)\-$modules with vertex~$R\,;$ then, Alperin's Conjecture affirms that the number of $G\-$conjugacy classes of weights of $G$ coincides with the number of isomorphism classes of simple $kG\-$modules, where $k$
is an algebraically closed field of characteristic~$p\,.$ More precisely, Alperin's Conjecture involves the blocks of $G$ as
we explain below.

\medskip
£1.2. In the case that $G$ is $p\-$solvable,  thirty years ago Tetsuro Okuyama [8] already proved that,
for any $p\-$subgroup $R$ of $G\,,$ the number of isomorphism classes $Y$ of simple $kN_G (R)\-$modules
with vertex $R$ coincides with the number of isomorphism classes of simple $kG\-$modules of vertex $R\,,$
which clearly shows Alperin's Conjecture restricted to $p\-$solvable groups. Once again, Okuyama's result actually 
involves the blocks of $G\,.$ Note that, setting $\bar N_G (R) = N_G (R)/R\,,$ a simple $kN_G(R)\-$module of vertex
$R$ is just the restriction of a {\it simple projective $k\bar N_G (R)\-$module\/}.

\medskip 
1.3. On the other hand, in [11,~6.4] we introduce a {\it multiplicity module\/} for any indecomposable $kG\-$module
$M\,,$ and in [11,~Lemma~9.9] we prove that  $M$ is determined by the triple formed by
a {\it vertex\/}~$R\,,$ an {\it $R\-$source\/} $E$ and a {\it multiplicity module\/} $V$ of $M$ --- an 
indecomposable projective $k_*\skew3\hat{\bar N}_G (R)_E\-$module where ${\bar N}_G (R)_E$ is the stabilizer
of the isomorphism class of $E$ in $\bar N_G (R)\,,$ $\skew3\hat{\bar N}_G (R)_E$ is the 
{\it central $k^*\-$extension\/} of $\bar N_G (R)_E$ determined by the action on ${\rm End}_k(E)\,,$ and 
$k_*\skew3\hat{\bar N}_G (R)_E$ is the corresponding {\it twisted\/} group algebra (cf.~£2.5 below) --- 
and that this correspondence actually defines a bijection
between the set of isomorphism classes of indecomposable $kG\-$modules and the set of $G\-$conjugacy classes
of triples $(R,E,V)$ formed by a $p\-$subgroup $R$ of $G\,,$ an indecomposable $kR\-$module $E$
of vertex $R$ and an indecomposable projective $k_*\skew3\hat{\bar N}_G (R)_E\-$module $V\,.$

 \medskip
 £1.4. Moreover, if $M$ is a simple $kG\-$module then it follows from [9,~Proposition~1.6]
 that $V$ is actually a {\it simple projective $k_*\skew3\hat{\bar N}_G (R)_E\-$module\/}. But, in the case that $G$ is 
 $p\-$solvable and $M$ is {\it primitive\/} --- namely, not induced from any proper subgroup --- it is well-known
 [17,~Lemma~30.4] that there is a $G\-$stable finite $p'\-$subgroup $K$ of ${\rm End}_k (M)^*$ generating 
 the $k\-$al-gebra~${\rm End}_k (M)\,.$ Consequently, in this case ${\rm End}_k(M)$ is actually a {\it Dade $R\-$algebra\/}
 [13,~1.3]; in particular, $\bar N_G (R)\-$stabilizes the isomorphism class of~$E$ [13,~1.8] and it follows from
 [15,~Theorem~9.21] that the central $k^*\-$extension $k_*\skew3\hat{\bar N}_G (R)$ above is {\it split\/}
 --- we are more explicit from £2.13 to £2.17  below.

\medskip
£1.5. That is to say, if $G$ is $p\-$solvable and $M$ a primitive simple $kG\-$mo-dule, then the pair formed by a 
vertex $R$ and by the isomorphism class of~the~restriction to $N_G (R)$ of a multiplicity $k\bar N_G (R)\-$module
$V$ --- after a choice~of~a splitting for the corresponding central $k^*\-$extension --- is actually a~{\it weight\/} of $G\,.$
More generally, since any simple $kG\-$module is certainly induced from a primitive simple $kH\-$module for 
some subgroup $H$ of $G\,,$ if~$G$ is $p\-$solvable then ${\rm End}_k (E)$ is always a {\it Dade $R\-$algebra\/} for any vertex $R$ and any $R\-$source $E$ of~$M\,;$ hence, in this case, the central $k^*\-$extension 
$k_*\skew3\hat{\bar N}_G (R)_E$ is always split and the corresponding multiplicity module $V$ becomes a
{\it simple projective $k\bar N_G (R)_E\-$module.\/}

\medskip
£1.6. In this paper, for a systematic  choice~of those splittings {\it via\/} a {\it polari-zation\/} [15,~9.5], on the one hand
we exhibit  a {\it natural bijection\/} --- namely compatible with the action of the group of {\it outer automorphisms\/} 
of $G$ --- between the sets of isomorphism classes of simple $kG\-$modules $M$ and of $G\-$conjugacy
classes of weights $(R,Y)$ of $G\,.$ On the other hand, we  determine the relationship between  a multiplicity 
$k_*\skew3\hat{\bar N}_G (R)_E\-$module $V$ and a simple $kN_G (R)\-$module $U$ with vertex~$R$ in the class
$Y$ of the corresponding {\it weight\/} of $G\,;$ explicitly,  there is a subgroup $N$ of $N_G (R)_E$ containing $R\,,$ a simple  $kN\-$module $W$ of vertex $R$ and, setting $\bar N = N/R\,,$ a group homomorphism $\theta\,\colon \bar N\to k^*$ 
in such a way that, denoting by $\bar W$ the corresponding $k\bar N\-$module and setting $\bar W_\theta
 = k_\theta\otimes_k \bar W\,,$   we have
$$U \cong {\rm Ind}_{N}^{N_G (R)} (W)\qq V \cong  {\rm Ind}_{\bar N}^{\bar N_G (R)_E} (\bar W_\theta) 
\eqno £1.6.1.$$
 The tools to carry out our purpose are mainly the {\it Fong reduction theorems\/} developed in~[16];  as in that paper, it is handy --- but {\it not\/} more general! --- to work systematically with
{\it $k^*\-$groups with finite $k^*\-$quotient $G$\/} [11,~\S5] --- namely, with central $k^*\-$extensions of $G\,.$

\medskip
£1.7. In 1994, when talking about this work at Beijing University, Zhang Jiping pointed out to us that Gabriel Navarro [7]  already had given a bijection between the above sets of isomorphism classes of simple $kG\-$modules and of $G\-$conjugacy classes of {\it weights\/} for finite groups of {\it odd\/} order, and therefore {\it solvable\/}. In our Appendix we show that Navarro's bijection corresponds indeed to the bijection obtained for a particular choice of the {\it splittings\/} above, a choice which is only possible for groups of odd order.

\bigskip
\noindent
{\bf £2. Notations and quoted results}

\bigskip
£2.1. We fix a prime number $p$ and an algebraically closed field $k$ of characteristic $p\,.$ We call
$k^*\-$group a group $X$ endowed with an injective group homomorphism $\theta\,\colon k^*\to Z(X)$
[11,~\S5], and call {\it $k^*\-$quotient\/} of $(X,\theta)$ the group~$X/\theta (k^*)\,;$ we denote by $X^\circ$
the $k^*\-$group formed by $X$ and by the composition of $\theta$ with the automorphism $k^*\cong k^*$
mapping $\lambda\in k^*$ on $\lambda^{-1}\,;$ we say that a $k^*\-$group is {\it finite\/} whenever
 its $k^*\-$quotient is finite. Usually, we denote by $\hat G$ a $k^*\-$group and by $G$ its $k^*\-$quotient,
 and we write $\lambda\.\hat x$ for the product of $\hat x\in \hat G$ and the image of $\lambda\in k^*$
 in $\hat G\,.$

 \medskip
 £2.2. If $\hat G'$ is a second $k^*\-$group, we denote by $\hat G\,\hat\times\, \hat G'$ the quotient
 of the direct product $\hat G\times \hat G'$ by the image in $\hat G\times \hat G'$ of the {\it inverse\/}
 diagonal of~$k^*\times k^*\,,$ which has an obvious structure of $k^*\-$group with $k^*\-$quotient
 $G\times G'\,;$ moreover, if $G = G'$ then we denote by $\hat G \star\hat G'$ the $k^*\-$group obtained from
 the inverse image of $\Delta (G)\i G\times G$ in $\hat G\,\hat\times\, \hat G'\,,$ which is nothing but
 the so-called {\it sum\/} of both central $k^*\-$extensions of $G\,;$  in particular, we have a {\it canonical\/} 
 $k^*\-$group isomorphism
 $$\hat G \star \hat G^\circ \cong k^*\times G
 \eqno £2.2.1.$$
 A $k^*\-$group homomorphism $\varphi\,\colon \hat G\to \hat G'$ is a group homomorphism which preserves
 the {\it $k^*\-$multiplication\/}; moreover, if $\hat G$ and $\hat G'$ are isomorphic then the group ${\rm Hom}(G,k^*)$
 acts {\it regularly\/} over the set of isomorphisms $\psi\,\colon \hat G \cong \hat G'$ and we denote by
 $\psi^\theta$ the $k^*\-$group isomorphism determined by $\theta\in {\rm Hom}(G,k^*)$ and $\psi\,.$
 We denote by $k^*\-\Gr$ the category of $k^*\-$groups.

 \medskip
 £2.3. Note that for any $k\-$algebra $A$ of finite dimension --- just called {\it $k\-$algebra\/} in the sequel ---
 the group $A^*$ of invertible elements has a canonical $k^*\-$group structure; we call {\it point\/} of $A$ 
 any $A^*\-$conjugacy class $\alpha$ of primitive idempotents of $A$ and denote by $A(\alpha)$ the simple
 quotient of $A$ determined by~$\alpha\,,$ and by $\P(A)$ the set of  {\it points\/} of $A\,.$
 If $S$ is a simple algebra then ${\rm Aut}_k (S)$ coincides with the $k^*\-$quotient of $S^*\,;$
 in particular, any finite group $G$ acting on $S$ determines --- by {\it pull-back\/} --- a $k^*\-$group $\hat G$
 of $k^*\-$quotient $G\,,$ together with a $k^*\-$group homomorphism [11,~5.7]
 $$\rho : \hat G\too S^*
 \eqno £2.3.1.$$

\medskip
£2.4. If $\hat G$ is a finite $k^*\-$group, we call {\it $\hat G\-$interior algebra\/} any $k\-$algebra $A$
endowed with a $k^*\-$group homomorphism
$$\rho : \hat G\too A^*
\eqno £2.4.1$$
and, as usual, we write $\hat x\. a$ and $a\.\hat x$ instead of $\rho(\hat x)a$ and $a\rho (\hat x)$ for any
$\hat x\in\hat G$ and any $a\in A\,;$ we say that $A$ is {\it primitive\/} whenever the unity element is 
primitive in $A^G\,.$ A {\it $\hat G\-$interior algebra homomorphism\/} from $A$ to another $\hat G\-$interior
algebra $A'$ is a {\it not necessarily unitary\/} algebra homomorphism\break
 \eject
 \noindent
 $f\,\colon A\to A'$ fulfilling
$f(\hat x\. a) = \hat x\. f(a)$ and $f(a\.\hat x) = f(a)\.\hat x\,;$ we say that $f$ is an {\it embedding\/}
 whenever ${\rm Ker}(f) = \{0\}$ and ${\rm Im}(f) = f(1)A'f(1)\,.$ Occasionally, it is handy to consider the 
 $(A'^{\hat G})^*\-$conjugacy class of $f$ that we denote by $\tilde f$ and call {\it exterior homomorphism\/}
 from $A$ to $A'\,;$ note that the {\it exterior homomorphisms\/} can be composed [9,~Definition~3.1].
 For a $k^*\-$group homomorphism $\varphi\,\colon \hat G'\to \hat G\,,$ we denote by ${\rm Res}_\varphi (A)$
 the $\hat G'\-$interior algebra defined by $\rho\circ\varphi\,.$ Note that the conjgation induces an action 
 of the $k^*\-$quotient $G$ of $\hat G$ on $A\,,$ so that $A$ becomes an ordinary $G\-$algebra; thus, all the {\it pointed group\/} language developed in [9] applies to $\hat G\-$interior algebras .

 \medskip
 £2.5. Namely, for any $k^*\-$subgroup $\hat H$ of $\hat G\,,$ a {\it point $\alpha$\/} of $\hat H$ on $A$
 is just a point of the $k\-$algebra $A^H\,,$ and the pair $\hat H_\alpha$ is a {\it pointed $k^*\-$group\/} on $A\,;$
we denote by $A(\hat H_\alpha)$ the simple quotient $A^H(\alpha)$ and, setting
$$\bar N_G (\hat H_\alpha) = N_G (\hat H_\alpha)/H\qq  A(\hat H_\alpha) = {\rm End}_k (V_\alpha)
\eqno £2.5.1,$$
by $\skew3\hat{\bar N}_G(\hat H_\alpha)$ the $k^*\-$group determined by the action of $\bar N_G (\hat H_\alpha)$
on $A(\hat H_\alpha)\,,$ so that $V_\alpha$ becomes a $\skew3\hat{\bar N}_G(\hat H_\alpha)\-$module called 
the {\it multiplicity $\skew3\hat{\bar N}_G(\hat H_\alpha)\-$module\/} of~$\hat H_\alpha$ [11,~6.4].
For any $i\in \alpha\,,$ $iAi$ has an evident structure of $\hat H\-$interior algebra mapping $\hat x\in \hat H$
on $\hat x\.i = i\.\hat x$ and we denote by $A_\alpha$ one of these mutually $(A^H)^*\-$conjugate $\hat H\-$interior
algebras. If $A'$ is another $\hat G\-$interior algebra and $f\,\colon A\to A'$ a {\it $\hat G\-$interior algebra embedding,\/}
$f(\alpha)$ is contained in a unique point $\alpha'$ of $\hat H$ on $A'\,,$ usually identified with $\alpha\,,$ and $f$ 
induces a $k\-$algebra embedding, a $k^*\-$group isomorphism and an $\hat H\-$interior algebra isomorphism
$$A(\hat H_\alpha)\too A'(\hat H_{\alpha'})\quad ,\quad \skew3\hat{\bar N}_G(\hat H_\alpha)\cong
\skew3\hat{\bar N}_G(\hat H_{\alpha'})\qq  A_\alpha \buildrel f_\alpha^{\alpha'}\over\cong A'_{\alpha'}
\eqno £2.5.2.$$

\medskip
£2.6. A second pointed $k^*\-$group $\hat K_\beta$ on $A$ is {\it contained\/} in $\hat H_\alpha$ if
$\hat K$ is a $k^*\-$subgroup of $\hat H$ and, for any $i\in \alpha\,,$ there is $j\in \beta$ such that 
$ij = j = ji\,;$ then, it is quite clear that the $(A^K)^*\-$conjugation induces a $\hat K\-$interior algebra 
embedding
$$f_\beta^\alpha : A_\beta \too {\rm Res}_{\hat K}^{\hat H}(A_\alpha)
\eqno £2.6.1.$$
More generally, we say that an injective $k^*\-$group homomorphism $\varphi\,\colon \hat K\to \hat H$
is an {\it $A\-$fusion from $\hat K_\beta$ to $\hat H_\alpha$\/} whenever there is a $\hat K\-$interior 
algebra embed-ding
$$f_\varphi : A_\beta \too {\rm Res}_\varphi(A_\alpha)
\eqno £2.6.2\phantom{.}$$
such that the inclusion $A_\beta\i A$ and the composition of $f_\varphi$ with the inclusion $A_\alpha\i A$
are $A^*\-$conjugate; then, the {\it exterior embedding\/} $\tilde f_\varphi$ is uniquely determined [10,~2.8].
We denote by $F_A (\hat K_\beta,\hat H_\alpha)$ the set of {\it $\hat H\-$conjugacy
classes\/} of $A\-$fusions from $\hat K_\beta$ to $\hat H_\alpha$ [12,~Definition~2.5] and we simply set
$$F_A (\hat H_\alpha) =F_A (\hat H_\alpha,\hat H_\alpha)
\eqno £2.6.3;$$
note that the conjugation in $\hat G$ induces a canonical group homomorphism
$$\bar N_G (\hat H_\alpha) \too F_A (\hat H_\alpha) 
\eqno £2.6.4.$$
If $A'$ is another $\hat G\-$interior algebra and $f\,\colon A \to A'$ a $\hat G\-$interior algebra embedding,
it follows from [10,~Proposition~2.14] that we have
$$F_A (\hat K_\beta,\hat H_\alpha) = F_{A'} (\hat K_\beta,\hat H_\alpha)
\eqno £2.6.5.$$

\medskip
£2.7. Note that any $p\-$subgroup $P$ of $\hat G$ can be identified with its image in~$G$ and
determines the $k^*\-$subgroup $k^*\.P \cong k^*\times P$ of $\hat G\,;$ as usual, we consider 
the {\it Brauer quotient\/} and the {\it Brauer algebra homomorphism\/}
$${\rm Br}_P : A^P\too A(P) = A^P\big/\sum_Q A_Q^P
\eqno £2.7.1,$$
where $Q$ runs over the set of proper subgroups of $P\,,$ and call {\it local\/} any point $\gamma$ of $P$
on $A$ not contained in ${\rm Ker}({\rm Br}_P)\,;$ recall that all the {\it maximal local pointed groups $P_\gamma$\/}
on $A$ contained in $\hat H_\alpha$ --- called {\it defect pointe groups of $\hat H_\alpha$\/} --- are
mutually $H\-$conjugate  [9,~Theorem~1.2], and that the $k\-$algebras $A_\alpha$ and $A_\gamma$ 
are {\it Morita equivalent\/} [9,~Corollary~3.5]. If $A_\gamma = iAi$ for $i\in \gamma\,,$ it follows from
[10,~Corollary~2.13] that we have a group homomorphism
$$F_A (P_\gamma) \too N_{A^{^*}_\gamma}(P\.i)\big/ P\.(A_\gamma^P)^*
\eqno £2.7.2\phantom{.}$$
and we consider the $k^*\-$group $\hat F_A (P_\gamma)$ defined by the {\it pull-back\/}
$$\matrix{F_A (P_\gamma)&\too& N_{A_\gamma^{^*}}(P\.i)/P\.(A_\gamma^P)^*\cr 
\uparrow&\phantom{\big\uparrow}&\uparrow\cr
\hat F_A (P_\gamma) &\too &N_{A_\gamma^{^*}}(P\.i)\big/P\.\big(i + J(A_\gamma^P)\big)\cr}
\eqno £2.7.3.$$

\medskip
£2.8. Then, from [11,~Proposition~6.12] suitably extended to $k^*\-$groups,
it follows that the group homomorphism~£2.6.4 can be lifted to a canonical $k^*\-$group
homomorphism
$$\skew3\hat{\bar N}_G(P_\gamma) * \bar N_{\hat G} (P_\gamma)^\circ \too \hat F_A (P_\gamma)^\circ
\eqno £2.8.1\phantom{.}$$
which, for any $\hat x\in \bar N_{\hat G} (P_\gamma) = N_{\hat G} (P_\gamma)/P$ and any $a\in (A^P)^*$ having the 
same action on $A(P_\gamma)\,,$ maps the element $(x,\bar a) * \hat x^{-1}$ of $\skew3\hat{\bar N}_G(P_\gamma) * \bar N_{\hat G} (P_\gamma)^\circ$ on the~pair [11,~Proposition~6.10]
$$\big( \tilde x^{-1},\overline{i( \hat x^{-1}\.a)i}\,\big)\in \hat F_A (P_\gamma)
\eqno £2.8.2,$$
where $x$ denotes the image of $\hat x$ in $\bar N_G (P_\gamma)\,,$ $\bar a$ the image of $a$ in $A(P_\gamma)\,,$
$\tilde x$ the image of $x$ in $F_A (P_\gamma)$ {\it via\/} homomorphism~£2.6.4 and $\overline{i( \hat x^{-1}\.a)i}$
the image of $i( \hat x^{-1}\.a)i$ in the right-hand bottom of diagram~£2.7.3.
\eject

\medskip
£2.9. If $A'$ is another $\hat G\-$interior algebra and $f\,\colon A\to A'$ a $\hat G\-$interior algebra embedding, it follows from [11,~Proposition~6.8] that, denoting by $\gamma'$ the point of $P$ on $A'$ containing $f(\gamma)\,,$ we have a canonical
$k^*\-$group isomorphism
$$\hat F_{\tilde f} (P_\gamma) : \hat F_A (P_\gamma)\cong \hat F_{A'} (P_{\gamma'})
\eqno £2.9.1\phantom{.}$$
which, according to [11,~Proposition~6.21], is compatible with the corresponding $k^*\-$group homomorphisms~£2.8.1 
and~£2.5.2. More precisely, let $Q_\delta$ be another local pointed group on $A$ and denote by $\delta'$
the point of $Q$ on $A'$ containing $f(\delta)\,;$ if there is a group 
isomorphism $\varphi\,\colon Q\cong P$ which is an {\it $A\-$fusion\/} from $Q_\delta$ to $P_\gamma$ then,
according to equality~£2.6.5 above, $\varphi$ is also  an {\it $A'\-$fusion\/} from $Q_{\delta'}$ to $P_{\gamma'}\,,$
so that we have two $Q\-$interior algebra isomorphisms
$$f_\varphi : A_\delta\cong {\rm Res}_\varphi (A_\gamma)\qq f'_\varphi : A'_{\delta'}\cong {\rm Res}_\varphi (A'_{\gamma'})
\eqno £2.9.2$$
and the uniqueness of the {\it exterior isomorphisms\/} $\tilde f_\varphi$ and $\tilde f'_\varphi$ forces 
the equality
$$\tilde f'_\varphi\circ \tilde f_\delta^{\delta'} = {\rm Res}_\varphi( \tilde f_\gamma^{\gamma'})\circ \tilde f_\varphi
\eqno £2.9.3\,.$$
In particular, since by the very definition we have 
$$\hat F_{{\rm Res}_\varphi (A_\gamma)}(Q_\delta) = \hat F_A (P_\gamma)\qq \hat F_{{\rm Res}_\varphi (A'_{\gamma'})}(Q_{\delta'}) = \hat F_{A'} (P_{\gamma'})
\eqno £2.9.4,$$
we get the following commutative diagram of $k^*\-$group isomorphisms
$$\matrix{\hat F_A (Q_\delta) &\buildrel \hat F_{\tilde f_\varphi} (Q_\delta)\over \cong& \hat F_A (P_\gamma)\cr
\hskip-40pt {\scriptstyle \hat F_{\tilde f} (Q_\delta)}\hskip4pt\wr\!\Vert&\phantom{\Big\uparrow}
&\wr\Vert\hskip4pt{\scriptstyle \hat F_{\tilde f} (P_\gamma)}\hskip-40pt\cr
\hat F_{A'} (Q_{\delta'}) &\buildrel \hat F_{\tilde f'_\varphi} (Q_{\delta'})\over \cong& \hat F_{A' }(P_{\gamma'})\cr}
\eqno £2.9.5.$$

\medskip
£2.10. Il is clear that the inclusion $k^*\i k$ determines a $k\-$algebra homomorphism to $k$ from the group 
algebra $kk^*$ of the group $k^*\,,$ so that $k$ becomes a $kk^*\-$algebra; for any finite $k^*\-$group
$\hat G\,,$ it is clear that the group algebra $k\hat G$ of the group $\hat G$ is also a $kk^*\-$algebra and 
then, we call {\it $k^*\-$group algebra\/} of $\hat G$ the algebra
$$k_*\hat G = k\otimes_{kk^*} k\hat G
\eqno £2.10.1;$$
note that the dimension of $k_*\hat G$ is equal to $\vert G\vert\,.$ Coherently, a {\it block\/} of $\hat G$
is a primitive idempotent $b$ of the center $Z(k_*\hat G)\,,$ so that $\alpha = \{b\}$ is a point of~$\hat G$ 
on $k_*\hat G\,;$ as usual, we denote by ${\rm Irr}_k (\hat G,b)$ the set of {\it Brauer characters\/} of all
 the simple $k_*\hat G\,b\-$modules, which corresponds bijectively with the set of points $\P (k_*\hat G\,b)\,.$

\medskip
£2.11. Recall that for any $p\-$subgroup $P$ of $\hat G$ we have [11,~£2.10.2 and~Pro-position~5.15]
$$(k_*\hat G)(P)\cong k_*C_{\hat G} (P)
\eqno £2.11.1;$$
\eject
\noindent
in particular, if $P$ is normal in $G\,,$ since the kernel of the obvious $k\-$algebra homomorphism
$k_*\hat G\to k_*(\hat G/P)$ is contained in the {\it radical\/} $J(k_*\hat G)$ and contains 
${\rm Ker}({\rm Br}_P)\,,$ this isomorphism implies that {\it any point of $P$ on $k_*\hat G$ is local\/}.
Moreover, it follows from [10,~Theorem~3.1] that we have
\smallskip
\noindent
£2.11.2\quad {\it  For any pair of local pointed groups $P_\gamma$
and $Q_\delta$ on $k_*\hat G\,,$ a $k_*\hat G\-$fusion from $Q_\delta$ to $P_\gamma$ coincides
with the conjugation by an element $x\in G$ such that~$Q_\delta\i (P_\gamma)^x\,.$\/}

\medskip
£2.12. If $\hat G$ is a finite $k^*\-$group, $A$ a $\hat G\-$interior algebra and $\hat H$ a $k^*\-$sub-group
of $\hat G\,,$ as usual we denote by ${\rm Res}_{\hat H}^{\hat G} (A)$ the corresponding $\hat H\-$interior
algebra. Conversely, for  any $\hat H\-$interior algebra~$B\,,$ we consider the {\it induced $G\-$interior
algebra\/}
$${\rm Ind}_{\hat H}^{\hat G}(B) = k_*\hat G\otimes_{k_*\hat H} B \otimes_{k_*\hat  H} k_*\hat  G
\eqno £2.12.1,$$
where the distributive product is defined by the {\it formula\/}
$$(\hat x\otimes b\otimes \hat y)(\hat x'\otimes b'\otimes \hat y')
=\cases{\hat x\otimes b.\hat y\hat x'.b'\otimes \hat y'&if $\hat y\hat x'\in
\hat H$\cr  {}&{}\cr 
0 &otherwise\cr}
\eqno £2.12.2\phantom{.}$$
for any $\hat x,\hat y,\hat x',\hat y'\in \hat G$ and any $b,b'\in B\,,$ and where we map $\hat x\in \hat G$ on the element
$$\sum_{\hat y}\hat x\hat y\otimes 1_B\otimes \hat y^{-1} = \sum_{\hat y}\hat y\otimes 1_B
\otimes \hat y^{-1}\hat x
\eqno £2.12.3,$$
$\hat y\in \hat G$ running  over a set of representatives for 
$\hat G/\hat H\,.$

\medskip
£2.13. For a finite $p\-$group $P\,,$ we call {\it Dade $P\-$algebra\/} [13,~1.3] a simple algebra $S$ endowed
with an action of $P$ which stibilizes a basis of $S$ containing the unity element; actually, the action of $P$ 
on $S$ can be lifted to a unique group homomorphism $P\to S^*$ and usually we consider $S$ as a
$P\-$interior algebra; moreover, the {\it Brauer quotient\/} $S(P)$ is also a simple $k\-$algebra [13,~1.8]
which implies that $P$ has a unique local point $\rho$ on $S$ that very often we omit, respectively writing
$F_S (P)$ and $\hat F_S (P)$ instead of $F_S (P_\rho)$ and $\hat F_S (P_\rho)\,.$ Recall that two Dade $P\-$algebras $S$ and $S'$ are {\it similar\/} if $S$ can be {\it
embedded\/} (cf.~£2.4) in the {\it tensor product\/} ${\rm End}(N)\otimes_k
S'$ for a suitable $kP\-$module $N$ with a $P\-$stable basis~[13,~1.5
and~2.5.1]; we denote by $\D_k (P)$ the set of {\it similarity classes\/} and the {\it tensor product\/} induces a {\it group structure\/} on $\D_k (P)$ --- called the {\it
Dade group\/} of $P$ --- where the opposite
$P\-$algebra $S^\circ$ determines the inverse of the {\it similarity\/} class  of $S\,.$

\medskip
£2.14. As in [15,~9.3], it is handy to consider the category $\frak D_k$ where the objects are the pairs $(P,S)$ formed by a finite $p\-$group $P$ and by a Dade $P\-$algebra $S\,,$ and where a morphism from $(P,S)$ to a second $\frak D_k\-$object
$(P',S')$ are the pairs $(\pi,f)$ formed by a surjective group homomorphism $\pi\,\colon P\to P'$ such that
${\rm Ker} (\pi)$ is $F_S(P)\-$stable, and by a $P\-$interior algebra embedding
$$f : {\rm Res}_\pi (S')\too S
\eqno £2.14.1.$$
Then, we have functors $\frak f$ and $\hat \frak f$ mapping $(P,S)$ on $F_S(P)$ and $\hat F_S (P)$ [15,~9.5],
together a {\it natural map\/} $\hat \frak f\to \frak f$ mapping $(P,S)$ on the structural homomorphism
$$\hat F_S (P)\too F_S (P)
\eqno £2.14.2.$$

\medskip
£2.15. As in [15,~9.5], we call {\it polarization\/} any {\it natural map $\omega$\/} from the functor
$\hat\frak f\,\colon \frak D_k\too k^*\-\Gr$ above to the {\it trivial\/} one --- namely, to the functor
mapping $(P,S)$ on $k^*$ and $(\pi, f)$ on ${\rm id}_{k^*}$ --- such that if $T$ is a $P\-$algebra with
trivial $P\-$action then $\omega$ maps $(P,T)$ on the first projection in the isomorphism
$$\hat F_T (P)\cong k^*\times F_T (P)
\eqno £2.15.1\phantom{.}$$
obtained  from the corresponding pull-back~£2.7.3. The point is that, according to [15,~Theorem~9.21], there exists
such a {\it natural map\/}, and we will construct a bijection as announced above from any choice of a 
{\it polarization\/} $\omega\,,$ namely from any choice, in a {\it coherent\/} way, of a $k^*\-$group homomorphism
$$\omega_{(P,S)} : \hat F_S (P)\too k^*
\eqno £2.15.2$$
for any $\frak D_k\-$object $(P,S)\,.$
A first application of this existence concerns the {\it multiplicity modules\/} of the indecomposable $k_*\hat G\-$modules $M$ having a vertex~$P$ and a $P\-$source $N$ such that ${\rm End}_k (N)$ is a Dade $P\-$algebra.

\bigskip
\noindent
{\bf Lemma~£2.16.} {\it Let $\hat G$ be a finite $k^*\-$group, $M$ an indecomposable $k_*\hat G\-$module,
$P$ a vertex and $N$ a $P\-$source of $M\,;$ let us denote by $P_N$ the local pointed group on the $\hat G\-$interior
algebra ${\rm End}_k (M)$ determined by the pair $(P,N)\,.$ If  ${\rm End}_k (N)$ is a Dade $P\-$algebra then the action
of $\bar N_{G}(P_N) $ on the simple quotient $\big({\rm End}_k (M)\big)(P_N)$ can be lifted to a $k^*\-$group
homomorphism
$$\bar N_{\hat G}(P_N)\too \big({\rm End}_k (M)\big)(P_N)^*
\eqno £2.16.1.$$\/}
\par
\noindent
{\bf Proof:} In any case, this action determines a $k^*\-$group $\skew3\hat{\bar N}_{G}(P_N) $ and we have a canonical
$k^*\-$group homomorphism (cf.~£2.8.1)
$$\skew3\hat{\bar N}_G(P_N) * \bar N_{\hat G} (P_N)^\circ \too \hat F_{{\rm End}_k (N)} (P_N)^\circ
\eqno £2.16.2;$$
but, if ${\rm End}_k (N)$ is a Dade $P\-$algebra, the existence of a {\it polarization\/} implies that, in particular,
we have 
$$\hat F_{{\rm End}_k (N)} (P)\cong k^*\times F_{{\rm End}_k (N)} (P)
\eqno £2.16.3;$$
consequently, we get $\skew3\hat{\bar N}_G(P_N) \cong \bar N_{\hat G} (P_N)\,.$ We are done.
\eject

\medskip
£2.17. More generally, if $S$ is a Dade $P\-$algebra and $A$ a $P\-$interior algebra, it follows from
[12,~Theorem~5.3] that, for any subgroup $Q$ of $P\,,$ we have a canonical bijection between
the sets of local points of $Q$ on $A$ and on $S\otimes_k A\,;$ moreover, if $A$ admits a $P\times P\-$stable 
basis by the multiplication on both sides, where $P\times \{1\}$ and $\{1\}\times P$ act freely, it follows 
from [6,~Lemma~1.17] that, for any pair of local pointed groups $Q_\delta$ and $R_\varepsilon$ 
on $A\,,$ we have
$$F_{S\otimes_k A} (R_{S\times\varepsilon},Q_{S\times \delta}) = F_S (R,Q)\cap F_A (R_\varepsilon,Q_\delta)
\eqno £2.17.1\phantom{.}$$
where $S\times \varepsilon$ and $S\times \delta$ denote the corresponding local points of $R$ and 
$Q$ on $S\otimes_k A\,;$ in this case, since the choice of a {\it polarization\/} $\omega$ determines a $k^*\-$group homomorphism
$$\omega_{(Q,{\rm Res}_Q^P(S))} : \hat F_S (Q)\too k^*
\eqno £2.17.2,$$
it follows from [12,~Proposition~5.11] that the inclusion of $F_{S\otimes_k A} (Q_{S\times \delta})$
in $F_A (Q_\delta)$ can be lifted  to a  $k^*\-$group homomorphism determined by $\omega$
$$\Phi^\omega_S(Q_\delta) : \hat F_{S\otimes_k A} (Q_{S\times \delta})\too \hat F_A (Q_\delta)
\eqno £2.17.3.$$
More precisely, as in £2.9 above, if $A'$ is a  $P\-$interior algebra and $f\,\colon A\to A'$ a  $P\-$interior algebra
embedding, denoting by $\delta'$ the point of $Q$ on $A'$ containing $f(\delta)\,,$ from [12,~Proposition~5.11] 
we still get the following commutative diagram of $k^*\-$group homomorphism
$$\matrix{\hat F_{S\otimes_k A} (Q_{S\times \delta})&\buildrel \Phi^\omega_S(Q_\delta)\over{\hbox to 30pt{\rightarrowfill}}  
&\hat F_A (Q_\delta)\cr
\hskip-60pt {\scriptstyle \hat F_{\tilde{\rm id}_S\otimes\tilde f} (Q_{Q\times\delta})}\hskip4pt\wr\!\Vert
&\phantom{\Big\uparrow}
&\wr\Vert\hskip4pt{\scriptstyle  \hat F_{\tilde f} (Q_\delta)}\hskip-30pt\cr
\hat F_{S\otimes_k A'} (Q_{S\times \delta'})&\buildrel \Phi^\omega_S(Q_{\delta'})\over{\hbox to 30pt{\rightarrowfill}}  &\hat F_{A'} (Q_{\delta'})\cr}\eqno £2.17.4.$$

\bigskip
\noindent
{\bf £3. The weights revisited}
\bigskip
£3.1. Let $\hat G$ be a finite $k^*\-$group; we say that a local pointed group  $Q_\delta$ on $k_*\hat G$
is {\it selfcentralizing\/} if $C_P (Q) = Z(Q)$ for any local pointed group $P_\gamma$ on $k_*\hat G$
containing $Q_\delta\,,$ and that it is a {\it radical\/} whenever it is selfcentralizing and we have
$${\Bbb O}_p \big(F_{k_*\hat G} (Q_\delta)\big) = \{1\}
\eqno 
£3.1.1.$$
Recall that, according to [15,~4.8 and~Corollary~7.3], $Q_\delta$ is {\it selfcentralizing\/} if, denoting by $f$ the block of $C_{\hat G}(Q)$ determined by $\delta\,,$ the image $\bar f$ of $f$ in the $k^*\-$group algebra of $\bar C_{\hat G}(Q)
= C_{\hat G}(Q)/Z(Q)$ is a block of {\it defect zero\/}; note that, in this case, $\delta$ is the unique local point
of $Q$ on $k_*\hat G$ determining the block $f\,.$

\medskip
£3.2. As mentioned in~£1.2 above, a {\it weight\/} $(R,Y)$ of $\hat G$ is formed by a $p\-$subgroup
$R$ of $\hat G$ and by the isomorphism class $Y$ of the restriction to $N_{\hat G} (R)$ of a {\it simple projective\/}
$k_*\bar N_{\hat G} (R)\-$module $V\,,$ where we set $\bar N_{\hat G}(R) = N_{\hat G}(R)/R\,;$\break
\eject
\noindent
 let us~denote by ${\rm Wgt}_k (\hat G)$ the set of $G\-$conjugacy classes of {\it weights\/} of $\hat G\,.$ 
 Then,~the restriction of $V$ to $\bar C_{\hat G} (R)\triangleleft \bar N_{\hat G} (R)$ is a {\it semisimple
projective\/} $k_*\bar C_{\hat G}(R)\-$module and thus any simple direct summand $W$ of
${\rm Res}_{\bar C_{\hat G}(R)}^{\bar N_{\hat G}(R)} (V)$ is also projective, so that it determines the unique
{\it local point\/} $\varepsilon$ of $R$ on $k_*\hat G$ (cf.~£2.11.1) in a block $\bar g$ of defect zero
of $\bar C_{\hat G} (R)\,;$ that is to say, $W$ determines a selfcentralizing pointed group $R_\varepsilon$
on $k_*\hat G$ and the stabilizer of the isomorphism class of $W$ in $N_{\hat G} (R)$ coincides with $N_{\hat G} (R_\varepsilon)\,.$

\medskip
£3.3. Moreover, it follows from isomorphism~£2.11.1 that we have
$$(k_*\hat G)(R_\varepsilon)\cong k_*\bar C_{\hat G} (R)\,\bar g\cong {\rm End}_k (W)
\eqno £3.3.1\phantom{.}$$
and from £2.5 we know that $W$ becomes an $\skew3\hat{\bar N}_G (R_\varepsilon)\-$module;
then, since the $\bar N_{\hat G} (R)\-$interior algebra ${\rm End}_k (V)$ is isomorphic to a suitable
{\it block algebra\/} of~$\bar N_{\hat G} (R)\,,$ and since we have (cf.~£2.10)
$$N_{\hat G} (R_\varepsilon)/ C_{\hat G} (R) \cong F_{k_*\hat G} (R_\varepsilon)
\eqno £3.3.2,$$
if follows from [16,~Theorem~£3.7] and from~£2.8 above that, for a suitable {\it simple projective\/}
$k_*\hat F_{k_*\hat G}(R_\varepsilon)\-$module $U$ restricted to 
$\skew3\hat{\bar N}_G(R_\varepsilon)^\circ * \bar N_{\hat G} (R_\varepsilon)$ {\it via\/} homomorphism~£2.8.1, we obtain
$$V\cong {\rm Ind}_{\bar N_{\hat G} (R_\varepsilon)}^{\bar N_{\hat G} (R)} (W\otimes_k U)
\eqno £3.3.3;$$
in particular, we get
$${\Bbb O}_p \big(F_{k_*\hat G}(R_\varepsilon)\big) = \{1\}
\eqno £3.3.4,$$
so that $R_\varepsilon$ is a {\it radical\/} pointed group.

\medskip
£3.4. Conversely, if $R_\varepsilon$ is a radical pointed group on $k_*\hat G$ and $U$ a {\it simple projective\/}
$k_*\hat F_{k_*\hat G}(R_\varepsilon)\-$module, it is easily checked that the restriction of $U$ to 
$\skew3\hat{\bar N}_G(R_\varepsilon)^\circ * \bar N_{\hat G} (R_\varepsilon)$ throughout homomorphism~£2.8.1, 
together with a {\it multiplicity $\skew3\hat{\bar N}_G(R_\varepsilon)\-$module $W$\/} of $R_\varepsilon$ define
a {\it simple projective \/} $k_*\bar N_{\hat G}(R)\-$module {\it via\/} the tensor product and the induction as in~£3.3.3.
In conclusion, we have proved that
\smallskip
\noindent
£3.4.1.\quad {\it The above correspondence between the sets of $G\-$conjugacy classes of weights  $(R,Y)$ 
of  $\hat G$ and of pairs $(R_\varepsilon,X)$ formed by a radical pointed group $R_\varepsilon$ on $k_*\hat G$ and by an 
isomorphism class $X$ of simple projective $k_*\hat F_{k_*\hat G}(R_\varepsilon)\-$mo-dules is bijective.\/}
\smallskip
\noindent
Let us call {\it $b\-$weight\/} of $\hat G$ any pair $(R_\varepsilon,X)$ formed by a radical pointed group
$R_\varepsilon$ on $k_*\hat G\,b$ and by an isomorphism class $X$ of simple projective 
$k_*\hat F_{k_*\hat G}(R_\varepsilon)\-$modules, and let us denote by ${\rm Wgt}_k(\hat G,b)$ the set of
$G\-$conjugacy classes of {\it $b\-$weights\/} of $\hat G\,;$ thus, statement~£3.4.1 affirms that we have
a ca-nonical bijection
$${\rm Wgt}_k(\hat G)\cong \bigsqcup_b {\rm Wgt}_k(\hat G,b)
\eqno £3.4.2\phantom{.}$$
where $b$ runs over the set of blocks of $\hat G\,;$ in particular, any {\it weight\/} of $\hat G$ determines a block.

\bigskip
\noindent
{\bf £4. Fitting pointed groups}
\bigskip
£4.1. Let us say that a finite $k^*\-$group $\hat G$ is {\it $p\-$solvable\/} if the $k^*\-$quotient $G$
of $\hat G$ is so; it is in this case that the following definition is actually useful. We call {\it Fitting pointed group\/}
of $\hat G$ any radical pointed group $Q_\delta$ on $k_*\hat G$ fulfilling the following condition
\smallskip
\noindent
£4.1.1.\quad {\it For any local pointed groups $P_\gamma$ and $R_\varepsilon$ on $k_*\hat G$ such that 
$P_\gamma$ contains $Q_\delta$ and $R_\varepsilon\,,$ any $k_*\hat G\-$fusion from $R_\varepsilon$
 to $P_\gamma$ coincides with the conjugation by an element $x\in N_G (Q_\delta)$ fulfilling 
 $R_\varepsilon\i (P_\gamma)^x\,.$\/}
\smallskip
\noindent
Note that this condition implies that a {\it Fitting pointed group\/} $Q_\delta$ of $\hat G$ is normal
in any local pointed group $P_\gamma$ containing $Q_\delta\,.$

\bigskip
\noindent
{\bf Proposition~£4.2.} {\it Let $\hat G$ be a finite $k^*\-$group and $Q_\delta$ a Fitting pointed group
of $\hat G\,.$ If a local pointed group $P_\gamma$ on $k_*\hat G$ contains both $Q_\delta$ and a radical pointed
group $R_\varepsilon$ on $k_*\hat G\,,$ then $R_\varepsilon$ contains $Q_\delta\,.$
In particular, $Q_\delta$ is the unique Fitting pointed group of $\hat G$ contained in $P_\gamma\,.$\/}

\medskip
\noindent
{\bf Proof:} We already know that $Q_\delta$ is normal in $P_\gamma$ and therefore the product $Q\.R$ is a subgroup
of $P\,;$ but, any element $y\in N_G (R_\varepsilon)$ induces by conjugation a $k_*\hat G\-$fusion from 
$R_\varepsilon$ to $P_\gamma$ and therefore, according to condition~£4.1.1, this $k_*\hat G\-$fusion
is also induced by an element $x\in N_G (Q_\delta)\,;$ in particular, the image of $N_Q (R_\varepsilon)$ in
$F_{k_*\hat G} (R_\varepsilon)$ is a normal $p\-$subgroup and therefore it is trivial.
\smallskip
On the other hand, it follows from [2,~Theorem~1.8] and from~£3.1 above that $\varepsilon$ is the unique local point 
of $R$ on $k_*\hat G$ such that $R_\varepsilon\i P_\gamma\,,$ and thus we have $N_Q (R_\varepsilon) = N_Q (R)\,;$
moreover, since $R_\varepsilon$ is selfcentralizing, we still have $C_P (R) =Z(R)$ and therefore $\bar N_{Q\.R} (R)$
maps injectively into the group of {\it outer\/} automorphisms of $R\,.$

\smallskip
Consequently, we get $\bar N_{Q\.R} (R) = \{1\}$ which implies that $Q\.R = R\,,$ so that $Q\i R\,;$
finally, once again it follows from [2,~Theorem~1.8] and from £3.1 above that $Q_\delta\i R_\varepsilon\,.$ Since
any Fitting pointed group is a radical, the last statement is now clear. We are done.

\bigskip
\noindent
{\bf Corollary~£4.3.} {\it Let $\hat G$ be a finite $k^*\-$group and $P_\gamma$ a maximal local pointed group
on $k_*\hat G\,.$ A radical pointed group $Q_\delta$ on $k_*\hat G$ contained in $P_\gamma$ is a Fitting 
pointed group of $\hat G$ if and only if it is contained in each radical pointed group on $k_*\hat G$ contained in
$P_\gamma\,.$\/}

\medskip
\noindent
{\bf Proof:} It follows from Proposition~£4.2 that this condition  is necessary. Conversely, if $Q_\delta$ is contained in any
radical pointed group on $k_*\hat G$ contained in~$P_\gamma\,,$ it follows from [14,~Theorem~A.9] that,
in particular, $Q_\delta$ is contained in each {\it essential\/} pointed group $R_\varepsilon$ contained in $P_\gamma\,;$
moreover, for any $x\in G$ normalizing either $R_\varepsilon$ or $P_\gamma\,,$ $(Q_\delta)^x$ is a Fitting
pointed group on $k_*\hat G$ contained in $P_\gamma$ and therefore it coincides with $Q_\delta\,;$
hence, $N_{\hat G} (Q_\delta)$ contains $N_{\hat G} (P_\gamma)$ and $N_{\hat G} (R_\varepsilon)$ for each
{\it essential\/} pointed group $R_\varepsilon$  contained in $P_\gamma\,.$ At this point,
condition~£4.1.1 follows from [14,~Corollary~A.12].

\medskip 
£4.4. From now on, we assume that {\it $\hat G$ is a $p\-$solvable finite $k^*\-$group\/} and let $b$ be a block
of $\hat G$ and $P_\gamma$ a maximal local pointed group on $k_*\hat G\,b\,;$ it follows from [16,~Theorem~4.6]
that there exists a {\it $P\-$source\/} pair $(S,\hat L)\,,$ unique up to isomorphisms, formed by a {\it primitive
Dade $P\-$algebra\/} $S$ and by a $p\-$solvable finite $k^*\-$group $\hat L$ containing $P\,,$ which fulfills the following
two conditions
\smallskip
\noindent
£4.4.1.\quad {\it $C_L\big({\Bbb O}_p(L)\big) = Z\big({\Bbb O}_p (L)\big)$ where $L$ denotes the $k^*\-$quotient
of $\hat L\,.$\/}
\smallskip
\noindent
£4.4.2.\quad {\it There is a $P\-$interior algebra embedding $\,e_\gamma\,\colon (k_*\hat G)_\gamma
\too S\otimes_k k_*\hat L\,.$\/}
\smallskip
\noindent
Note that, according to isomorphism~£2.11.1, any $p\-$subgroup of $L$ containing $\Bbb O_p (L)$ has a unique
local point on $k_*\hat L$ --- actually, it coincides with $\{1\}$ (cf.~£2.10). In particular, $P$ has a unique local
point $\dot\gamma= \{1\}$ on $k_*\hat L$ and therefore it follows from [12,~Theorem~5.3] that it has also a unique local point $S\times \dot\gamma$ on $S\otimes_k k_*\hat L\,;$ then, the embedding above is equivalent to the existence of
a $P\-$interior algebra isomorphism
$$(k_*\hat G)_\gamma \cong (S\otimes_k k_*\hat L)_{S\times \dot\gamma}
\eqno £4.4.3.$$

\medskip
£4.5. In particular, from isomorphism~£4.4.3 and from [12,~Theorem~5.3], any local pointed group $Q_\delta$ 
on $k_*\hat G$ contained in $P_\gamma$ determines a local pointed group $Q_\delta$ on $k_*\hat L$
 and this correspondence is bijective. Moreover, since we have  $P\-$interior algebra embeddings
 $$k_*\hat L\too S^\circ \otimes_k S\otimes_k k_*\hat L \longleftarrow S^\circ\otimes_k (k_*\hat G)_\gamma
 \eqno £4.5.1\phantom{.}$$
 and $(S^\circ\otimes_k S)\times \dot\gamma$ is the unique local point of $P$ on 
 $S^\circ\otimes_k S\otimes_k k_*\hat L\,,$ we still have a $P\-$interior algebra embedding
 $$e_\gamma^\circ : k_*\hat L\too S^\circ \otimes_k (k_*\hat G)_\gamma
 \eqno £4.5.2\phantom{.}$$
inducing the same bijection between the sets of local pointed groups on $(k_*\hat G)_\gamma$ and on $k_*\hat L\,;$
then, since $(k_*\hat G)_\gamma$ and $k_*\hat L$ admit $P\times P\-$stable bases\break 
\eject
\noindent
by the multiplication on both sides, where $P\times \{1\}$ and $\{1\}\times P$ act freely, it follows from £2.17 above 
applied twice that we have
$$F_{k_*\hat G} (R_\varepsilon,Q_\delta) = F_{k_*\hat L} (R_{\dot\varepsilon},Q_{\dot\delta}) \i F_S (R,Q)
\eqno £4.5.3\phantom{.}$$
for any pair of local pointed groups $Q_\delta$ and $R_\varepsilon$ on $k_*\hat G$ contained in $P_\gamma\,,$
and that the choice of a {\it polarization\/} $\omega$ and of the embedding $e_\gamma$ determine  $k^*\-$isomorphisms (cf.~£2.8.3 and~£2.17.3)
$$\hat F_{k_*\hat G} (Q_\delta)\buildrel \hat F_{\tilde e_\gamma}(Q_\delta)\over\cong 
\hat F_{S\otimes_k k_*\hat L}(Q_{S\times \dot\delta})\buildrel \Phi^\omega_S(Q_{\dot\delta})\over 
\cong \hat F_{k_*\hat L} (Q_{\dot\delta}) 
\eqno £4.5.4.$$

\medskip
£4.6. Set $O = \Bbb O_p (L)$ and denote by $\dot\eta$ and by $\eta$ the respective unique local points
of $O$ on $k_*\hat L$ and on $(k_*\hat G)_\gamma\,;$ since $O_{\dot\eta}$ is clearly a {\it Fitting pointed group\/}
of $\hat L\,,$ it follows from £4.5 above that {\it $O_\eta$ is a Fitting pointed group of\/} $\hat G\,.$
Moreover, from the $k^*\-$group homomorphism~£2.8.1 and from the last statement in~£2.10,
we get the $k^*\-$group isomorphism
$$\hat L/O\cong \hat F_{k_*\hat L} (O_{\dot\eta})
\eqno £4.6.1\phantom{.}$$
and therefore the choice of a {\it polarization\/} $\omega$ determines a $k^*\-$isomorphism
$$\hat L/O \cong \hat F_{k_*\hat G} (O_\eta)
\eqno £4.6.2.$$

\bigskip
\noindent
{\bf Remark £4.7.} It follows from~[16,~4.7] that the Dade $P\-$algebra $S$ above
always come from a suitable {\it nilpotent block\/} admitting $P$ as a {\it defect group\/} 
and therefore, according to~[15,~Theorem~7.8], the {\it similarity\/} class of $S$ in the {\it
Dade group\/} $\D_k (P)$ is a torsion element (cf.~£2.13). In particular, we can restrict our
{\it polarizations\/} to the {\it full\/} subalgebra $\frak D_k^{\rm tor}$ of $\frak D_k$
over the objects $(P,S)$ fulfilling this condition.

\bigskip
\noindent
{\bf £5. The key parameterizations}
\bigskip
£5.1. Let $\hat G$ be a $p\-$solvable finite $k^*\-$group, $b$ a block of $\hat G$ and 
$P_\gamma$ a maximal local pointed group on $k_*\hat G\,b\,,$ and denote by $(S,\hat L)$ a {\it $P\-$source\/} pair of this block and by $O_\eta$ the {\it Fitting pointed group\/} of $\hat G$ contained in $P_\gamma\,;$ in this section, our purpose is to show that the choice of a 
{\it polarization\/} $\omega$ determines
two bijections
$$\eqalign{\Gamma^\omega_{(\hat G,b)} &: {\rm Irr}_k (\hat G,b) \cong {\rm Irr}_k 
\big(\hat F_{k_*\hat G} (O_\eta)\big)\cr
\Delta^\omega_{(\hat G,b)} &: {\rm Wgt}_k (\hat G,b) \cong {\rm Wgt}_k 
\big(\hat F_{k_*\hat G} (O_\eta)\big)\cr}
\eqno £5.1.1\phantom{.}$$
which are {\it natural\/} with respect to the isomorphisms between blocks. We first need to know the group
of {\it exterior automorphisms\/} ${\rm Out}_P\big((k_*\hat G)_\gamma\big)$ (cf.~£2.4) of the $P\-$interior algebra 
$(k_*\hat G)_\gamma\,;$ recall that, according to [11,~Proposition~14.9],\break 
\eject
\noindent
 we have an injective group homomorphism
$${\rm Out}_P\big((k_*\hat G)_\gamma\big)\too {\rm Hom} \big(F_{k_*\hat G} (P_\gamma),k^*\big)
\eqno £5.1.2\phantom{.}$$
and therefore ${\rm Out}_P\big((k_*\hat G)_\gamma\big)$ is Abelian.

\bigskip
\noindent
{\bf Proposition~£5.2.} {\it With the notation above, there are group isomorphisms
$${\rm Out}_P\big((k_*\hat G)_\gamma\big)\cong {\rm Out}_P(k_*\hat L)\cong {\rm Hom} (L,k^*)
\eqno £5.2.1\phantom{.}$$
mapping $\tilde\sigma\in {\rm Out}_P\big((k_*\hat G)_\gamma\big)$ on an element 
$\dot{\tilde\sigma}\in {\rm Out}_P(k_*\hat L)$ such that, for any $P\-$interior algebra embedding 
$e_\gamma\,\colon (k_*\hat G)_\gamma\to S\otimes_k k_*\hat L$ we have
$$\tilde e_\gamma\circ \tilde\sigma = (\widetilde{\rm id}_S \otimes \dot{\tilde\sigma}) \circ \tilde e_\gamma
\eqno £5.2.2,$$
and mapping $\zeta\in {\rm Hom} (L,k^*)$ on the exterior class of the $P\-$interior algebra automorphism
of $k_*\hat L$ sending $\hat y\in \hat L$ to $\zeta (y)\.\hat y$ where $y$ is the image of $\hat y$ in $L\,.$
Moreover, ${\rm Out}_P\big((k_*\hat G)_\gamma\big)$ acts regularly over the set of exterior embeddings from
$(k_*\hat G)_\gamma$ to $S\otimes_k k_*\hat L\,.$\/}

\medskip
\noindent
{\bf Proof:} Since $S^\circ\times \gamma$ is the unique local point of $P$ on $S^\circ \otimes_k (k_*\hat G)_\gamma\,,$
embedding~£4.5.2 induces a $P\-$interior algebra isomorphism
$$k_*\hat L\cong \big(S^\circ \otimes_k (k_*\hat G)_\gamma\big)_{S^\circ\times \gamma}
\eqno £5.2.3\phantom{.}$$
and therefore, for a representative $\sigma$ of $\tilde\sigma\,,$ the automorphism ${\rm id}_S\otimes \sigma$
of $S^\circ \otimes_k (k_*\hat G)_\gamma\,,$ composed with a suitable inner automorphism, induces an
automorphism $\dot\sigma$ of $k_*\hat L$ and it is quite clear that the {\it exterior class $\dot{\tilde\sigma}$\/} of 
$\dot\sigma$ does not depend on our choices, and fulfills 
$$\tilde e_\gamma^\circ \circ \dot{\tilde\sigma} = (\widetilde{\rm id}_S \otimes \tilde\sigma) \circ \tilde e_\gamma^\circ
\eqno £5.2.4.$$

\smallskip
Tensoring embedding~£4.5.2 by $S$ and arguing as in~£4.5 above, it is not difficult to prove that equality~£5.2.2
also holds. Similarly, since this correspondence comes from ``conjugation'' {\it via\/} the {\it exterior\/} class
of isomorphisms~£4.4.3 and~£5.2.3, it is clear that it is a group isomorphism; actually, this argument also proves the 
last statement.

\smallskip
 On the other hand, for any $\zeta\in {\rm Hom} (L,k^*)\,,$ it is clear that the map sending $\hat y\in \hat L$
 to $\zeta(y)\.\hat y$ defines an automorphism of the $k^*\-$group $\hat L$ inducing the identity on $P$ and thus,
 it determines a $P\-$interior algebra automorphism of~$k_*\hat L\,;$ moreover, since $y$ is also the image of
$\zeta(y)\.\hat y$ in $L\,,$  we clearly get a group homomorphism
 $${\rm Hom} (L,k^*)\too {\rm Aut}_P(k_*\hat L) 
 \eqno £5.2.5.$$
 \eject
\noindent 
Conversely, any  $P\-$interior algebra automorphism $\dot\sigma$ of~$k_*\hat L$ stabilizes the {\it Fitting pointed group\/} $O_{\dot\eta}\,,$ acting trivially on $O\,;$ hence, it acts on the $k^*\-$group $\hat F_{k_*\hat L} (O_{\dot\eta})$
acting trivially on its $k^*\-$quotient $F_{k_*\hat L} (O_{\dot\eta})\i {\rm Out}(O)$ and therefore, according
to isomorphism~£4.6.1 above, it determines an element of 
$${\rm Hom} (L/O,k^*) ={\rm Hom} (L,k^*)
\eqno £5.2.6;$$
clearly, any inner  $P\-$interior algebra automorphism of~$k_*\hat L$ determines the trivial element  of ${\rm Hom} (L,k^*)$ and thus, we easily get the second isomorphism in~£5.2.1.

\medskip
£5.3. We are ready to define the first bijection in~£5.1.1. Since the restriction determines a {\it Morita equivalence\/} 
between the $k\-$algebras $k_*\hat G\,b$ and~$(k_*\hat G)_\gamma$  (cf.~£2.7),  we certainly have a {\it natural\/} bijection (cf.~£2.10)
$${\rm Irr}_k (\hat G,b)\cong \P \big((k_*\hat G)_\gamma\big)
\eqno £5.3.1\phantom{.}$$
and  any embedding $e_\gamma\,\colon (k_*\hat G)_\gamma\to S\otimes_k k_*\hat L$ induces an injective
map and a $k^*\-$group isomorphism (cf.~£2.5 and~£2.8.3)
$$\eqalign{\P(\tilde e_\gamma) &: \P\big((k_*\hat G)_\gamma\big)\too \P (S\otimes_k k_*\hat L)\cr 
\hat F_{\tilde e_\gamma}(O_\eta) &:\hat F_{k_*\hat G} (O_\eta)\cong 
\hat F_{S\otimes_k k_*\hat L}(O_{S\times\dot\eta})\cr}
\eqno £5.3.2;$$
then, the existence of embedding~£4.5.2 proves that the map $\P(e_\gamma) $ is actually bijective.
On the other hand, the choice of a {\it polarization\/} $\omega$ determines a $k^*\-$group isomorphism
$$\Phi^\omega_S (O_{\dot\eta}) : \hat F_{S\otimes_k k_*\hat L}(O_{S\times\dot\eta})
\cong \hat F_{k_*\hat L} (O_{\dot\eta}) 
\eqno £5.3.3.$$
Finally, isomorphism~£4.6.1 determines a canonical bijection
$$\Gamma_{\hat L} : {\rm Irr}_k (\hat L)\cong {\rm Irr}_k \big(\hat F_{k_*\hat L} (O_{\dot\eta})\big)
\eqno £5.3.4.$$

\bigskip
\noindent
{\bf Corollary~£5.4.} {\it With the notation and the choice above, there is a bijection
$$\Gamma^\omega_{(\hat G,b)} : {\rm Irr}_k (\hat G,b) \cong {\rm Irr}_k  \big(\hat F_{k_*\hat G} (O_\eta)\big)
\eqno £5.4.1\phantom{.}$$
such that, for any  embedding $e_\gamma\,\colon (k_*\hat G)_\gamma\to S\otimes_k k_*\hat L\,,$ 
we have the commutative diagram
$$\matrix{{\rm Irr}_k (\hat G,b)&\cong &\P (S\otimes_k k_*\hat L)
&\cong &{\rm Irr}_k (\hat L)\cr
\hskip-30pt{\scriptstyle \Gamma^\omega_{(\hat G,b)}}\hskip3pt\wr\!\Vert&&\phantom{\Big\uparrow}
&&\hskip-15pt{\scriptstyle \Gamma_{\hat L}}\hskip3pt\wr\!\Vert\cr
{\rm Irr}_k  \big(\hat F_{k_*\hat G} (O_\eta)\big)&\cong &
{\rm Irr}_k  \big(\hat F_{S\otimes_k k_*\hat L}(O_{S\times\dot\eta})\big)&\cong 
&{\rm Irr}_k \big(\hat F_{k_*\hat L} (O_{\dot\eta})\big)\cr}
\eqno £5.4.2.$$\/}

\medskip
\noindent
{\bf Proof:} It is clear that, for a choice of an  embedding 
$$e_\gamma : (k_*\hat G)_\gamma\too  S\otimes_k k_*\hat L
\eqno £5.4.3,$$
\eject
\noindent
 the  bijections~£5.3.1 and $\P(\tilde e_\gamma)\,,$ and the $k^*\-$group isomorphism
$\hat F_{\tilde e_\gamma}(O_\eta)$ above determine the horizontal left-hand bijections in diagram~£5.4.2;
the top hori-zontal right-hand bijection follow from~£2.10 and~£2.17, and the bottom hori-zontal right-hand bijection
from isomorphism~£5.3.3 up to the choice of $\omega\,;$ then, the bijection $\Gamma_{\hat L}$ and
the commutativity of the diagram define the bijection~$\Gamma^\omega_{(\hat G,b)}\,.$

\smallskip
We claim that this bijection does not depend on the choice of $e_\gamma\,;$ indeed, for another choice $e'_\gamma$
 of this embedding, it follows from Proposition~£5.2 that there is $\tilde\sigma\in 
 {\rm Out}_P\big((k_*\hat G)_\gamma\big)$ fulfilling 
 $$\tilde e'_\gamma = \tilde e_\gamma\circ \tilde\sigma = (\widetilde{\rm id}_S \otimes \dot{\tilde\sigma}) 
 \circ \tilde e_\gamma
 \eqno £5.4.4\phantom{.}$$
 and therefore, with obvious notation, we get the following commutative diagrams
 $$\hskip-30pt\matrix{{\rm Irr}_k (\hat G,b)&\cong &\P \big((k_*\hat G)_\gamma\big)&\buildrel \P(\tilde e'_\gamma)\over\cong &\P (S\otimes_k k_*\hat L)\cr
\Vert&\phantom{\big\uparrow}&\Vert &&\wr\Vert\hskip4pt{\scriptstyle \P(\tilde{\rm id}_S \otimes 
\dot{\tilde\sigma})}\hskip-40pt\cr
{\rm Irr}_k (\hat G,b)&\cong &\P \big((k_*\hat G)_\gamma\big)&\buildrel \P(\tilde e_\gamma)\over\cong &\P (S\otimes_k k_*\hat L)\cr}
\eqno £5.4.5\phantom{.}$$
\smallskip
$$\hskip-30pt\matrix{{\rm Irr}_k  \big(\hat F_{k_*\hat G} (O_\eta)\big)&\buildrel {\rm Irr}_k (\hat F_{\tilde e'_\gamma}(O_\eta))\over \cong &
{\rm Irr}_k  \big(\hat F_{S\otimes_k k_*\hat L}(O_{S\times\dot\eta})\big)\cr
\Vert&\phantom{\Big\uparrow}&\wr\Vert\hskip4pt{\scriptstyle {\rm Irr}_k (\hat F_{\tilde{\rm id}_S \otimes 
\dot{\tilde\sigma}}(O_{S\times\dot\eta}))}\hskip-70pt \cr
{\rm Irr}_k  \big(\hat F_{k_*\hat G} (O_\eta)\big)&\buildrel {\rm Irr}_k (\hat F_{\tilde e_\gamma}(O_\eta))\over\cong &
{\rm Irr}_k  \big(\hat F_{S\otimes_k k_*\hat L}(O_{S\times\dot\eta})\big)\cr}
\eqno £5.4.6.$$

\smallskip
Moreover, we have the evident commutative diagram
$$\matrix{\P (S\otimes_k k_*\hat L)
&\cong &{\rm Irr}_k (\hat L)\cr
\hskip-50pt{\scriptstyle \P(\tilde{\rm id}_S \otimes \dot{\tilde\sigma})}\hskip4pt\wr\!\Vert
&\phantom{\Big\uparrow}&\wr\Vert\hskip4pt{\scriptstyle {\rm Irr}_k (\dot{\tilde\sigma})}\hskip-30pt\cr
\P (S\otimes_k k_*\hat L)
&\cong &{\rm Irr}_k (\hat L)\cr}
\eqno £5.4.7;$$
on the other hand, since the groups of $k^*\-$group automorphisms of $\hat F_{k_*\hat L} (O_{\dot\eta})$ and 
$\hat F_{S\otimes_k k_*\hat L}(O_{S\times\dot\eta})$ which induce the identity over (cf.~£2.17 applied twice)
$$F_{k_*\hat L} (O_{\dot\eta}) = F_{S\otimes_k k_*\hat L}(O_{S\times\dot\eta})
\eqno £5.4.8\phantom{.}$$
are both canonically isomorphic to the Abelian group ${\rm Hom}\big(F_{k_*\hat L} (O_{\dot\eta}),k^*\big)\,,$
 we still have the commutative diagram
$$\matrix{\
{\rm Irr}_k  \big(\hat F_{S\otimes_k k_*\hat L}(O_{S\times\dot\eta})\big)&\cong 
&{\rm Irr}_k \big(\hat F_{k_*\hat L} (O_{\dot\eta})\big)\cr
\hskip-70pt{\scriptstyle {\rm Irr}_k (\hat F_{\tilde{\rm id}_S \otimes 
\dot{\tilde\sigma}}(O_{S\times\dot\eta}))}\hskip4pt\wr\!\Vert&\phantom{\Big\uparrow}
&\wr\Vert\hskip4pt{\scriptstyle {\rm Irr}_k (\hat F_{\dot{\tilde\sigma}}(O_{\dot\eta}))}\hskip-50pt  \cr
{\rm Irr}_k  \big(\hat F_{S\otimes_k k_*\hat L}(O_{S\times\dot\eta})\big)&\cong 
&{\rm Irr}_k \big(\hat F_{k_*\hat L} (O_{\dot\eta})\big)\cr}
\eqno £5.4.9.$$

\smallskip
Finally, from isomorphism~£4.6.2 we obviously   get the following commutative diagram
$$\matrix{{\rm Irr}_k (\hat L)&\buildrel {\rm Irr}_k (\dot{\tilde\sigma})\over \cong& {\rm Irr}_k (\hat L)\cr
\hskip-15pt{\scriptstyle \Gamma_{\hat L}}\hskip3pt\wr\!\Vert&&\wr\Vert\hskip4pt
{\scriptstyle \Gamma_{\hat L}}\hskip-15pt\cr
{\rm Irr}_k \big(\hat F_{k_*\hat L} (O_{\dot\eta})\big)&\buildrel {\rm Irr}_k (\hat F_{\dot{\tilde\sigma}}(O_\eta))\over \cong &{\rm Irr}_k \big(\hat F_{k_*\hat L} (O_{\dot\eta})\big)\cr}
\eqno £5.4.10;$$
now, our claim follows from putting together all these commutative diagrams.

\medskip
£5.5. In order to define the second bijection in~£5.1.1, let $(R_\varepsilon,X)$ be a {\it $b\-$weight\/} of $\hat G\,;$ 
for our purposes, we may assume that 
$P_\gamma$ contains $R_\varepsilon\,;$ then, $R_\varepsilon$ and $R_{\dot\varepsilon}$ respectively contain 
$O_\eta$ and $O_{\dot\eta}\,;$ recall that, with the notation and the choice in~£4.5 above, we have  $k^*\-$group
isomorphisms
$$\hat F_{k_*\hat G} (R_\varepsilon)\buildrel \hat F_{\tilde e_\gamma}(R_\varepsilon)\over\cong 
\hat F_{S\otimes_k k_*\hat L}(R_{S\times \dot\varepsilon})\buildrel \Phi^\omega_S(R_{\dot\varepsilon})\over 
\cong \hat F_{k_*\hat L} (R_{\dot\varepsilon}) 
\eqno £5.5.1\phantom{.}$$
and, in this case, $X$ determines an isomorphism class $\dot X$ of {\it simple projective\/} 
$k_*\hat F_{k_*\hat L} (R_{\dot\varepsilon})\-$modules; moreover, we clearly have $N_{\hat L}(R_{\dot\varepsilon}) 
= N_{\hat L} (R)$ and from~£2.11.2 it is easily checked that the $k^*\-$group homomorphism~£2.8.3 induces a $k^*\-$group
isomorphism
$$\bar N_{\hat L} (R_{\dot\varepsilon})\cong \hat F_{k_*\hat L} (R_{\dot\varepsilon})
\eqno £5.5.2;$$
consequently, the pair $(R,\dot X)$ is a {\it weight\/} of $\hat L\,.$

\bigskip
\noindent
{\bf Proposition~£5.6.} {\it  With the notation and the choice above, let $(R_\varepsilon,X)$ and $(R'_{\varepsilon'},X')$ 
be $b\-$weights of $\hat G$ such that $P_\gamma$ contains  $R_\varepsilon$ and $R'_{\varepsilon'}\,.$
If $(R_\varepsilon,X)$ and $(R'_{\varepsilon'},X')$ are $G\-$conjugate then  the corresponding weights 
$(R,\dot X)$ and $(R',\dot X')$  of $\hat L$ are $L\-$conjugate. In particular, this correspondence induces  a bijection
$${\rm Wgt}_k^\omega (e_\gamma) : {\rm Wgt}_k (\hat G,b)\cong {\rm Wgt}_k (\hat L)
\eqno £5.6.1.$$\/}
\par
\noindent
{\bf Proof:} Assume that  $(R'_{\varepsilon'},X')^x = (R_{\varepsilon},X)$ for some $x\in G\,;$
then, the conjugation by $x$ determines a {\it $(k_*\hat G)_\gamma\-$fusion\/} $\varphi$
from $R_\varepsilon$ to $R'_{\varepsilon'}$ (cf.~£2.6), and the corresponding $R\-$interior algebra isomorphism
$$f_\varphi : (k_*\hat G)_\varepsilon \cong {\rm Res}_\varphi \big((k_*\hat G)_{\varepsilon'}\big)
\eqno £5.6.2\phantom{.}$$
induces a $k^*\-$group isomorphism (cf~£2.8.3)
$$\hat F_{\tilde f_\varphi}(R_\varepsilon) : \hat F_{(k_*\hat G)_\varepsilon }(R_\varepsilon)\cong \hat F_{(k_*\hat G)_{\varepsilon'}}(R'_{\varepsilon'})
\eqno £5.6.3;$$
actually, we have  $X = {\rm Res}_{\hat F_{\tilde f_\varphi}(R_\varepsilon)}(X')\,.$ 
\eject

\smallskip
But, according to equality~£4.5.3,  the group homomorphism $\varphi$ is also a {\it $k_*\hat L\-$fusion\/} 
from $R_{\dot\varepsilon}$ to  $R'_{\dot\varepsilon'}\,,$ so that $\varphi\,\colon R\cong R'$ is also induced by
some element~$\dot x\in L$ (cf.~statement~£2.11.2);
moreover, we have the corresponding $R\-$interior algebra isomorphism
$$\dot f_\varphi : (k_*\hat L)_{\dot\varepsilon} \cong {\rm Res}_\varphi \big((k_*\hat L)_{\dot\varepsilon'}\big)
\eqno £5.6.4\phantom{.}$$
inducing  a $k^*\-$group isomorphism (cf~£2.8.3)
$$\hat F_{\dot{\tilde f}_\varphi}(R_{\dot\varepsilon}) : \hat F_{(k_*\hat L)_{\dot\varepsilon} }(R_{\dot\varepsilon})
\cong \hat F_{(k_*\hat L)_{\dot\varepsilon'}}(R'_{\dot\varepsilon'})
\eqno £5.6.5.$$

\smallskip
Then, the commutativity of diagrams~£2.9.5 and~£2.17.4 applied here yields the following commutative diagrams
of  $k^*\-$group isomorphisms
$$\matrix{\hat F_{(k_*\hat G)_\gamma }(R_\varepsilon)\hskip-5pt
&\buildrel \hat F_{\tilde e_\gamma}(R_\varepsilon)\over\cong \hskip-5pt
&\hat F_{S \otimes_k k_*\hat L }(R_{S\times \dot\varepsilon})\hskip-5pt
&\buildrel \Phi^\omega_S(R_{\dot\varepsilon})\over \cong \hskip-5pt 
&\hskip-5pt \hat F_{k_*\hat L }(R_{\dot\varepsilon})\cr
\hskip-30pt{\scriptstyle \hat F_{\tilde f_\varphi}(R_\varepsilon)}\hskip2pt\wr\!\Vert
&\phantom{\Big\uparrow}&\wr\Vert&
&\wr\Vert\hskip2pt {\scriptstyle \hat F_{\dot{\tilde f}_\varphi}(R_{\dot\varepsilon})}\hskip-30pt\cr
\hat F_{(k_*\hat G)_\gamma }(R'_{\varepsilon'})\hskip-5pt
 &\buildrel \hat F_{\tilde e_\gamma}(R'_{\varepsilon'})\over\cong \hskip-5pt
 &\hat F_{S \otimes_k k_*\hat L }(R'_{S\times \dot\varepsilon'})\hskip-5pt
 &\buildrel \Phi^\omega_S(R'_{\dot\varepsilon'})\over \cong \hskip-5pt
 &\hat F_{k_*\hat L }(R'_{\dot\varepsilon'})\cr }
 \eqno £5.6.6.$$
Consequently, we also have   $\dot X = {\rm Res}_{\hat F_{\dot{\tilde f}_\varphi}(R_{\dot\varepsilon})}(\dot X')$ 
and therefore we get 
$$(R'_{\dot\varepsilon'},\dot X')^{\dot x} = (R_{\dot\varepsilon},\dot X)
\eqno £5.6.7;$$
that is to say, the correspondence above induces a map
$${\rm Wgt}_k^\omega (e_\gamma) : {\rm Wgt}_k (\hat G,b)\too {\rm Wgt}_k (\hat L)
\eqno £5.6.8.$$
which is quite clear that it is a bijection. We are done.

\bigskip
\noindent
{\bf Proposition~£5.7.} {\it With the the notation above, the canonical $k^*\-$group isomorphism  $\hat L/O\cong 
\hat F_{k_*\hat L} (O_{\dot\eta})$ induces a bijection
$$\Delta_{\hat L} : {\rm Wgt}_k (\hat L) \cong {\rm Wgt}_k \big(\hat F_{k_*\hat L} (O_{\dot\eta})\big)
\eqno £5.7.1.$$\/}

\par
\noindent
{\bf Proof:} Let $(R,Y)$ be a {\it weight\/} of $\hat L\,;$ since the unity element in $k_*\hat L$ is a block of~$\hat L$
and condition~£4.4.1 holds,  $R$ has a unique local point $\dot\varepsilon$ on $k_*\hat L$ and $R_{\dot\varepsilon}$
is a radical pointed group which contains $O_{\dot\eta}$ (cf.~Corollary~£4.3); moreover, since~we have the $k^*\-$group isomorphism $\hat L/O\cong \hat F_{k_*\hat L} (O_{\dot\eta})$ (cf.~£4.6.2),  setting $\bar R = R/O$ and identifying
$\bar R$ with its image in $\hat F_{k_*\hat L} (O_{\dot\eta})\,,$ the normalizer $N_{\hat L}(R_{\dot\varepsilon}) = N_{\hat L} (R)$ is just the converse image in $\hat L$ of $N_{\hat F_{k_*\hat L} (O_{\dot\eta})} (\bar R)$ 
and therefore we have the canonical $k^*\-$group isomorphism
$$\bar N_{\hat L}(R) \cong \bar N_{\hat F_{k_*\hat L} (O_{\dot\eta})} (\bar R)
\eqno £5.7.2\,;$$
\eject
\noindent
in particular, $Y$ determines an isomorphism class $\bar Y$ of simple $N_{\hat F_{k_*\hat L} (O_{\dot\eta})} (\bar R)\-$ modules of vertex $\bar R\,,$ so that the pair $(\bar R,\bar Y)$ is a {\it weight\/} of $\hat F_{k_*\hat L} (O_{\dot\eta})\,.$

\smallskip
Conversely, if we start with a {\it weight\/} $(\bar R,\bar Y)$ of $\hat F_{k_*\hat L} (O_{\dot\eta})\,,$
it is clear that, for the converse image $R$ of $\bar R$ in $\hat L\,,$ isomorphism~£5.6.2 still holds and therefore
$\bar Y$  determines an isomorphism class $Y$ of simple $k_*\hat L\-$modules of vertex $R\,,$ so that
 the pair $(R,Y)$ is a {\it weight\/} of $\hat L\,.$ Since this correspondence is compatible with the $L\-$conjugation,
 we get the announced bijection~£5.6.1.

\bigskip
\noindent
{\bf Corollary~£5.8.} {\it  With the notation and the choice above, there is a bijection
$$\Delta^\omega_{(\hat G,b)} : {\rm Wgt}_k (\hat G,b) \cong {\rm Wgt}_k  \big(\hat F_{k_*\hat G} (O_\eta)\big)
\eqno £5.8.1\phantom{.}$$
such that, for any  embedding $e_\gamma\,\colon (k_*\hat G)_\gamma\to S\otimes_k k_*\hat L\,,$ 
we have the commutative diagram
$$\matrix{{\rm Wgt}_k (\hat G,b)&\buildrel {\rm Wgt}_k^\omega (e_\gamma)\over\cong &{\rm Wgt}_k (\hat L)\cr
\hskip-30pt{\scriptstyle \Delta^\omega_{(\hat G,b)}}\hskip3pt\wr\!\Vert&\phantom{\Big\uparrow}
&\hskip-15pt{\scriptstyle \Delta_{\hat L}}\hskip3pt\wr\!\Vert\cr
{\rm Wgt}_k  \big(\hat F_{k_*\hat G} (O_\eta)\big)&\cong 
&{\rm Wgt}_k \big(\hat F_{k_*\hat L} (O_{\dot\eta})\big)\cr}
\eqno £5.8.2.$$
where the bottom bijection is induced by the $k^*\-$group isomorphisms
$$\hat F_{k_*\hat G} (O_\eta)\buildrel \hat F_{\tilde e_\gamma}(O_\eta)\over\cong 
\hat F_{S\otimes_k k_*\hat L}(O_{S\times \dot\eta})\buildrel \Phi^\omega_S(O_{\dot\eta})\over 
\cong \hat F_{k_*\hat L} (O_{\dot\eta}) 
\eqno £5.8.3.$$\/}

\par
\noindent
{\bf Proof:} It is clear that, for a choice of $e_\gamma\,,$ Propositions~£5.6 and~£5.7, and the commutativity of the diagram define the bijection $\Delta^\omega_{(\hat G,b)}\,.$ We claim that this bijection does not depend on this choice; indeed, for another choice $e'_\gamma$
 of this embedding, it follows from Proposition~£5.2 that there is $\tilde\sigma\in 
 {\rm Out}_P\big((k_*\hat G)_\gamma\big)$ fulfilling 
 $$\tilde e'_\gamma = \tilde e_\gamma\circ \tilde\sigma = (\widetilde{\rm id}_S \otimes \dot{\tilde\sigma}) 
 \circ \tilde e_\gamma
 \eqno £5.8.4;$$
in particular, if $(R_\varepsilon,X)$ is a {\it $b\-$weight\/} of $\hat G$ and $(R_{\dot\varepsilon},\dot X)$ the corresponding {\it weight\/} of $\hat L$ in~£5.5 above --- $\dot X$ is the isomorphism class of a {\it simple projective\/} 
$k_*\hat F_{k_*\hat L} (R_{\dot\varepsilon})\-$module $V$  restricted to $N_{\hat L}(R)$ ---  
then ${\rm Wgt}_k^\omega (e'_\gamma)$ sends~the $G\-$conjugacy class of $(R_\varepsilon,X)$ to the $L\-$conjugacy class of $(R_{\dot\varepsilon},\dot X')$ where~$\dot X' $ is the isomorphism class of corresponding the {\it simple projective\/} 
$k_*\hat F_{k_*\hat L} (R_{\dot\varepsilon})\-$mo-dule 
${\rm Res}_{\hat F_{\dot{\tilde\sigma}^{-1}}(R_{\dot\varepsilon})} (V)\,,$ since  ${\rm Hom}(L,k^*)$ clearly
acts trivially on the set of local pointed groups on $k_*\hat L$ and we have the commutative diagram~(cf.~£2.17.4)
$$\matrix{\hat F_{S \otimes_k k_*\hat L }(R_{S\times \dot\varepsilon})&\buildrel \Phi^\omega_S(R_{\dot\varepsilon})\over \cong &\hat F_{k_*\hat L }(R_{\dot\varepsilon})\cr
\hskip-50pt{\scriptstyle \hat F_{\tilde{\rm id}_S\otimes\dot{\tilde\sigma}}(R_{S\times \dot\varepsilon})}\hskip2pt\wr\!\Vert&\phantom{\Big\uparrow}
&\wr\Vert\hskip2pt {\scriptstyle \hat F_{\dot{\tilde\sigma}}(R_{\dot\varepsilon})}\hskip-30pt\cr
\hat F_{S \otimes_k k_*\hat L }(R_{S\times \dot\varepsilon})&\buildrel \Phi^\omega_S(R_{\dot\varepsilon})\over \cong &\hat F_{k_*\hat L }(R_{\dot\varepsilon})\cr }
 \eqno £5.8.5.$$
 \eject

 \smallskip
 Now, setting $\bar R = R/O\,,$ since we have (cf.~isomorphisms~£4.6.1 and~£5.5.2)
 $$\hat F_{k_*\hat L} (R_{\dot\varepsilon})\cong \bar N_{\hat L} (R_{\dot\varepsilon}) = \bar N_{\hat L} (R)
 \cong  \bar N_{\hat F_{k_*\hat L} (O_{\dot\eta}) } (\bar R)
 \eqno £5.8.6,$$
$V$ determines a  {\it simple projective\/}  $\bar N_{\hat F_{k_*\hat L} (O_{\dot\eta}) } (\bar R)\-$module $\bar V\,;$
moreover, since   ${\rm Hom}(L,k^*)$ clearly stabilizes $\hat L$ and it acts trivially on $P\,,$ the corresponding representative $\dot\sigma$ of $\dot{\tilde\sigma}$ induces a $k^*\-$group automorphism $\dot{\bar\sigma}$ 
of~$\hat F_{k_*\hat L} (O_{\dot\eta}) $ (cf.~isomorphism~£4.6.1) which stabilizes $\bar N_{\hat F_{k_*\hat L} 
(O_{\dot\eta}) } (\bar R)\,,$ and it is quite clear that, with obvious notation, we get the following commutative diagram
 $$\matrix{\hat F_{k_*\hat L} (R_{\dot\varepsilon})&\cong &  \bar N_{\hat F_{k_*\hat L} (O_{\dot\eta}) } (\bar R)\cr
\hskip-30pt{\scriptstyle \hat F_{\dot{\tilde\sigma}}(R_{\dot\varepsilon})}\hskip2pt\wr\!\Vert  
&\phantom{\Big\uparrow}&\wr\Vert\hskip2pt {\scriptstyle \bar N_{\dot{\bar\sigma}}(\bar R)}\hskip-30pt\cr
 \hat F_{k_*\hat L} (R_{\dot\varepsilon})&\cong &  \bar N_{\hat F_{k_*\hat L} (O_{\dot\eta}) } (\bar R)\cr}
\eqno £5.8.7;$$
hence, {\it via\/} isomorphisms~£5.6.7, ${\rm Res}_{\hat F_{\dot{\tilde\sigma}^{-1}}(R_{\dot\varepsilon})} (V)$
 determines the  {\it simple projective\/}  $\bar N_{\hat F_{k_*\hat L} (O_{\dot\eta}) } (\bar R)\-$module 
 ${\rm Res}_{\bar N_{\dot{\bar\sigma}^{-1}}(\bar R)} (\bar V)\,.$

 \smallskip
 At this point, denoting by $\skew3\dot{\bar X}$ and $\dot{\bar X'}$ the respective isomorphism classes 
 of the $\bar N_{\hat F_{k_*\hat L} (O_{\dot\eta}) } (\bar R)\-$modules $\bar V$ and 
 ${\rm Res}_{\bar N_{\dot{\bar\sigma}^{-1}}(\bar R)} (\bar V)\,,$ it follows from Proposition~£5.7 that
 $\Delta_{\hat L}$ maps $(R_{\dot\varepsilon},\dot X)$ on $(\bar R,\skew3\dot{\bar X})\,,$ and 
 $(R_{\dot\varepsilon},\dot X')$ on $(\bar R,\skew3\dot{\bar X'})\,.$ But, we also have the
 following commutative diagram (cf.~£2.17.4)
$$\matrix{\hat F_{S \otimes_k k_*\hat L }(O_{S\times \dot\eta})&\buildrel \Phi^\omega_S(O_{\dot\eta})\over \cong &\hat F_{k_*\hat L }(O_{\dot\eta})\cr
\hskip-50pt{\scriptstyle \hat F_{\tilde{\rm id}_S\otimes\dot{\tilde\sigma}}(O_{S\times \dot\eta})}\hskip2pt\wr\!\Vert&\phantom{\Big\uparrow}
&\wr\Vert\hskip2pt {\scriptstyle \hat F_{\dot{\tilde\sigma}}(O_{\dot\eta})}\hskip-30pt\cr
\hat F_{S \otimes_k k_*\hat L }(O_{S\times \dot\eta})&\buildrel \Phi^\omega_S(O_{\dot\eta})\over \cong 
&\hat F_{k_*\hat L }(O_{\dot\eta})\cr }
 \eqno £5.8.8\phantom{.}$$
and we consider its restriction to all the normalizers of $\bar R\,.$ Consequently, since we have (cf.~£5.8.4)
$$\hat F_{\tilde e'_\gamma}(O_\eta) = \hat F_{\tilde{\rm id}_S\otimes\dot{\tilde\sigma}}(O_{S\times \dot\eta})\circ \hat F_{\tilde e_\gamma}(O_\eta)
\eqno £5.8.9,$$
the corresponding bottom bijections in diagram~£5.8.2 maps the {\it weights\/} $(\bar R,\skew3\dot{\bar X})$ and 
$(\bar R,\skew3\dot{\bar X'})$ of $\hat F_{k_*\hat L }(O_{\dot\eta})$ on the {\it same weight\/}
of $\hat F_{k_*\hat L} (O_{\dot\eta})\,.$ We are done.

\bigskip
\noindent
{\bf £6. The Fitting block sequences}
\bigskip
£6.1. In order to exhibit bijections between the sets of  isomorphism classes of simple $k_*\hat G\-$modules
and of  $G\-$conjugacy classes of {\it weights\/} of~$\hat G\,,$ we need a third set, namely the set of $G\-$conjugacy classes of {\it Fitting block sequences\/} of~$\hat G\,.$ We call  {\it Fitting block sequence\/} of~$\hat G$
any sequence $\B = \{(\hat G_n,b_n)\}_{n\in \Bbb N}$ of pairs formed by a $k^*\-$group $\hat G_n$ and by 
a block $b_n$ of $\hat G_n\,,$ such that $\hat G_0 = \hat G$ and that, for any $n\in \Bbb N\,,$ 
we have $\hat G_{n+1} = \hat F_{k_*\hat G_{n}}(O^n_{\,\eta_n})$ for some {\it  Fitting pointed\break
\eject
\noindent
 group\/}~$O^n_{\,\eta_n}$ of~$\hat G_n\,.$ Note that, since clearly $\vert G_{n+1} \vert\le \vert G_{n} \vert\,,$ 
 such a sequence stabilizes, and actually we have $\vert G_{n+1} \vert= \vert G_{n} \vert$ if and only if $b_n$ 
 is a block of {\it defect zero\/} of~$\hat G_n$ (cf.~statement~£4.4.1). Moreover, for any $h\in \Bbb N\,,$ the sequence 
$\B_h = \{(\hat G_{h +n},b_{h+ n})\}_{n\in \Bbb N}$ is clearly a   {\it Fitting block sequence\/} of~$\hat G_h\,.$

\medskip
£6.2. If  $\hat G'$ is a $k^*\-$group isomorphic to $\hat G$ and $\theta\,\colon \hat G\cong \hat G'$
a $k^*\-$group isomorphism of $\hat G\,,$ it is quite clear that, from
any  {\it Fitting block sequence\/}  $\B = \{(\hat G_n,b_n)\}_{n\in \Bbb N}$ of~$\hat G\,,$ we
are able to construct a {\it Fitting block sequence\/} $\B' = \{(\hat G'_n,b'_n)\}_{n\in \Bbb N}$ of~$\hat G'$ inductively defining a sequence of $k^*\-$group isomorphisms~$\theta_n\,\colon \hat G_n\cong \hat G'_n$ by $\theta_0 = \theta$
and, for any  $n\in \Bbb N\,,$ by (cf.~£2.9)
$$\theta_{n+1} = \hat F_{\tilde\theta_n} (O^n_{\,\eta_n}) :  \hat F_{k_*\hat G_{n}}(O^n_{\,\eta_n})\cong 
\hat F_{k_*\hat G'_{n}}\big(\theta_n(O^n)_{\,\theta_n(\eta_n)}\big)
\eqno £6.2.1,$$
where we sill denote by $\theta_n\,\colon k_*\hat G_n\cong k_*\hat G'_n$ the corresponding $k\-$algebra isomorphism,
and  setting 
$$b'_n = \theta_n (b_n)\qq \hat G_{n+1} = \hat F_{k_*\hat G'_{n}}\big(\theta_n(O^n)_{\,\theta_n(\eta_n)}\big)
\eqno £6.2.2\phantom{.}$$
 for any  $n\in \Bbb N\,.$ In particular, the group of inner automorphisms of $G$
acts on the set of {\it Fitting block sequences\/} of~$\hat G$ and then we denote by  ${\rm Fbs}_k(\hat G)$ the set 
of ``$G\-$conjugacy classes'' of the {\it Fitting block sequences\/} of~$\hat G\,,$ and by $N_G (\B)$ 
the stabilizer of $\B$ in $G\,.$

\medskip
£6.3. In this section, our purpose is to show that the choice of a {\it polarization\/} $\omega$ determines two 
bijections
$${\rm Fbs}_k(\hat G)\cong {\rm Irr}_k (\hat G)\qq {\rm Fbs}_k(\hat G)\cong {\rm Wgt}_k(\hat G)
\eqno £6.3.1\phantom{.}$$
which are {\it natural\/} with respect to the $k^*\-$group isomorphisms, the composition of the inverse of the first one
with the second one being our announced parameterization.

\medskip 
£6.4. Let $\B = \{(\hat G_n,b_n)\}_{n\in \Bbb N}$ be a  {\it Fitting block sequence\/} of~$\hat G\,,$
so that we have $\hat G_{n+1} = \hat F_{k_*\hat G_{n}}(O^n_{\,\eta_n})$ for some Fitting pointed 
group~$O^n_{\,\eta_n}$ of~$\hat G_n$ and, choosing a {\it polarization\/} $\omega\,,$ we denote by 
$$\eqalign{\Gamma_{(\hat G_n,b_n)}^\omega &: {\rm Irr}_k (\hat G_n,b_n)\cong {\rm Irr}_k (\hat G_{n+1})\cr
\Delta_{(\hat G_n,b_n)}^\omega &: {\rm Wgt}_k (\hat G_n,b_n)\cong {\rm Wgt}_k (\hat G_{n+1})\cr}
\eqno £6.4.1\phantom{.}$$
the bijections coming from Corollaries~£5.4 and~£5.8. Let us call {\it character   sequence $\omega\-$associated to $\B$\/} any sequence  $\{\varphi_n\}_{n\in \Bbb N}$ where $\varphi_n$ belongs to ${\rm Irr}_k (\hat G_n,b_n)$ in such a way that we have
$$\Gamma_{(\hat G_n,b_n)}^\omega (\varphi_n) = \varphi_{n+1}
\eqno £6.4.2\phantom{.}$$
 for any $n\in \Bbb N\,.$ Similarly, let us call {\it weight sequence $\omega\-$associated to $\B$\/} any sequence 
 $\{\,\overline{\!(R^n,Y^n)\!}\,\}_{n\in \Bbb N}$ where $\,\overline{\!(R^n,Y^n)\!}\,$ is the $G_n\-$conjugacy class of a {\it weight\/} $(R^n,Y^n)$ of $\hat G_n\,,$ determining a $G_n\-$conjugacy class 
 $\,\overline{\!(R^n_{\,\varepsilon_n},X^n)\!}\,$ of {\it $b_n\-$weights\/} of~$\hat G_n$ (cf.~statement~£3.4.1) 
in such a way that we have
$$\Delta_{(\hat G_n,b_n)}^\omega \big(\,\overline{\!(R^n_{\,\varepsilon_n},X^n)\!}\,\big) = 
\,\overline{\!(R^{n+1},Y^{n+1})\!}\,
\eqno £6.4.3\phantom{.}$$
 for any $n\in \Bbb N\,.$

 \bigskip
 \noindent
 {\bf Theorem~£6.5.} {\it With the notation and the choice above, any Fitting block sequence 
 $\B = \{(\hat G_n,b_n)\}_{n\in \Bbb N}$ of~$\hat G$ admits a unique  character   sequence  
 $\{\varphi_n\}_{n\in \Bbb N}$  and a unique  weight sequence $\{\,\overline{\!(R^n,Y^n)\!}\,\}_{n\in \Bbb N}$  
 $\omega\-$associated to $\B\,.$ Moreover, the correspondences mapping $\B$ to $\varphi_0$ and to 
 $\,\overline{\!(R^0,Y^0)\!}\,$ induce two natural bijections
$${\rm Fbs}_k(\hat G)\cong {\rm Irr}_k (\hat G)\quad and\quad {\rm Fbs}_k(\hat G)\cong {\rm Wgt}_k(\hat G)
\eqno £6.5.1.$$ \/}

\par
\noindent
{\bf Proof:} Since the sequence $\B$ stabilizes, we can argue by induction on the ``length to stabilization''. If this length is
zero then the block $b_0$ is already of {\it defect zero\/} and therefore ${\rm Irr}_k (\hat G_0,b_0)$ has a  unique element
$\varphi_0$ and, setting $\varphi_n = \varphi_0$ for any $n\in \Bbb N\,,$ we get  a {\it character sequence\/}
$\omega\-$associated to $\B\,;$ similarly, ${\rm Wgt}_k (\hat G_0,b_0)$ has a  unique element and the corresponding constant sequence defines a {\it weight sequence\/}
$\omega\-$associated to $\B\,.$

\smallskip
If the  ``length to stabilization'' is not zero then the  Fitting block sequence 
 $\B_1 = \{(\hat G_{1 +n},b_{1 +n})\}_{n\in \Bbb N}$ of~$\hat G_1$ already admits a  character   sequence  
 $\{\varphi_{1 +n}\}_{n\in \Bbb N}$  and a   weight sequence $\{\,\overline{\!(R^{1 +n},Y^{1 +n})\!}\,\}_{n\in \Bbb N}$  
 $\omega\-$associated to $\B_1\,;$ then, in order to get a  character   sequence and a weight sequence  
 $\omega\-$associated to $\B\,,$ it suffices to define (cf.~£6.4.1)
 $$\eqalign{\varphi_0 &= (\Gamma_{(\hat G_0,b_0)}^\omega)^{-1}(\varphi_{1})\cr
 \,\overline{\!(R^0_{\,\varepsilon_0},X^0)\!}\, &=   (\Delta_{(\hat G_0,b_0)}^\omega)^{-1} 
 \big(\,\overline{\!(R^{1},Y^{1})\!}\,\big)\cr}
 \eqno £6.5.2\phantom{.}$$
and to consider the $G\-$conjugacy class $\,\overline{\!(R^{0},Y^{0})\!}\,$ of {\it weights\/} of $\hat G$ determined 
by~$\,\overline{\!(R^0_{\,\varepsilon_0},X^0)\!}\,$ (cf.~statement~£3.4.1).

\medskip
On the other hand, since the maps $\Gamma_{(\hat G_n,b_n)}^\omega$ and $\Delta_{(\hat G_n,b_n)}^\omega$
are bijective, equalities~£6.4.2 and~£6.4.3 show that a {\it character sequence\/} $\{\varphi_n\}_{n\in \Bbb N}$ and
a {\it weight sequence\/} $\{\,\overline{\!(R^n,Y^n)\!}\,\}_{n\in \Bbb N}$
$\omega\-$associated to $\B$ are uniquely determined by one of their terms; but, for $n$ big enough, we know that
$b_n$ is a block of {\it defect zero\/} of $\hat G_n$ and then $\varphi_n$ and $\,\overline{\!(R^n,Y^n)\!}\,$
are uniquely determined; consequently, $\{\varphi_n\}_{n\in \Bbb N}$  and $\{\,\overline{\!(R^n,Y^n)\!}\,\}_{n\in \Bbb N}$ are uniquely determined and it is quite clear that they only depend on the $G\-$conjugacy class of $\B\,;$
thus, since  $\Gamma_{(\hat G_n,b_n)}^\omega$ and $\Delta_{(\hat G_n,b_n)}^\omega$ are {\it natural\/},
 we have obtained two {\it natural\/} maps
$${\rm Fbs}_k(\hat G)\too {\rm Irr}_k (\hat G)\qq {\rm Fbs}_k(\hat G)\too {\rm Wgt}_k(\hat G)
\eqno £6.5.3\,.$$
\eject

\smallskip
We claim that they are both bijective; actually, we will define the inverse maps. For any $\varphi\in {\rm Irr}_k (\hat G)\,,$
we inductively define two sequences $\{\varphi_n\}_{n\in \Bbb N}$ and $\{(\hat G_n,b_n)\}_{n\in \Bbb N}$
by setting $\varphi_0 = \varphi\,,$ $\hat G_0 = \hat G$ and by denoting by $b_0$ the block of $\varphi\,,$ and further, for any 
$n\in \Bbb N\,,$ by setting
$$\varphi_{n+1} = \Gamma_{(\hat G_n,b_n)}^\omega (\varphi_n)\quad ,\quad \hat G_{n+1} = \hat F_{k_*\hat G_{n}}(O^n_{\,\eta_n})
\eqno £6.5.4\phantom{.}$$
 for some Fitting pointed group~$O^n_{\,\eta_n}$ on $k_*\hat G_nb_n\,,$ and by denoting by $b_{n+1}$ the 
 block of $\varphi_{n+1}\,;$ then, it is clear that  $\B =\{(\hat G_n,b_n)\}_{n\in \Bbb N}$ is a {\it Fitting block sequence\/}
 of~$\hat G$ and that  $\{\varphi_n\}_{n\in \Bbb N}$ becomes the {\it character sequence\/} $\omega\-$associated 
 to~$\B\,;$ note that, our construction only depends on the choice of the Fitting pointed group~$O^n_{\,\eta_n}$ on 
 $k_*\hat G_nb_n$ for a {\it finite  set of values of $n\,.$\/} Moreover, since all the  {\it Fitting pointed group\/} on 
 $k_*\hat G_nb_n$ are mutually $G_n\-$conjugate (cf.~Proposition~£4.2), $\varphi$ determines a unique $G\-$conjugacy
 class of {\it Fitting block sequence\/} of~$\hat G\,.$
 That is to say, we have obtained a map
 $${\rm Irr}_k (\hat G)\too {\rm Fbs}_k(\hat G)
 \eqno £6.5.5$$
 and it is easily checked that it is the inverse of the left-hand map in~£6.5.3.

 \smallskip
 Analogously, for any $\,\overline{\!(R,Y)\!}\,\in {\rm Wgt}_k(\hat G)\,,$ we inductively define two sequences 
 $\{\,\overline{\!(R^n,Y^n)\!}\,\}_{n\in \Bbb N}$  and $\{(\hat G_n,b_n)\}_{n\in \Bbb N}$
by setting $\,\overline{\!(R^0,Y^0)\!}\, = \,\overline{\!(R,Y)\!}\,,$ $\hat G_0 = \hat G$ and by~denoting by $b_0$ 
the block of $\hat G_0$ determined by $\,\overline{\!(R,Y)\!}\,$ (cf.~bijection~£3.4.2), and further, for any 
$n\in \Bbb N\,,$ by setting
$$\overline{\!(R^{n+1},Y^{n+1})\!}\, = \Delta_{(\hat G_n,b_n)}^\omega 
\big(\,\overline{\!(R^n_{\,\varepsilon_n},X^n)\!}\,\big) \quad ,\quad \hat G_{n+1} 
= \hat F_{k_*\hat G_{n}}(O^n_{\,\eta_n})
\eqno £6.5.6\phantom{.}$$
where $\,\overline{\!(R^n_{\,\varepsilon_n},X^n)\!}\,$ is the {\it $b_n\-$weight\/} of $\hat G_n$ determined 
by $\,\overline{\!(R^n,Y^n)\!}\,$ and $O^n_{\,\eta_n}$ a Fitting pointed group  on $k_*\hat G_nb_n\,,$
and by denoting by $b_{n+1}$ the block of determined by the {\it weight\/} $\,\overline{\!(R^{n+1},Y^{n+1})\!}\,$
 (cf.~bijection~£3.4.2); then, it is clear that  $\B =\{(\hat G_n,b_n)\}_{n\in \Bbb N}$ is a {\it Fitting block sequence\/}
 of~$\hat G$ and that $\{\,\overline{\!(R^n,Y^n)\!}\,\}_{n\in \Bbb N}$ becomes the {\it weight sequence\/}
 $\omega\-$associated to $\B\,.$ As above, our construction only depends on the choice of the Fitting pointed group~$O^n_{\,\eta_n}$ on  $k_*\hat G_nb_n$ for a {\it finite  set of values of $n$\/} and therefore we have obtained a map
 $${\rm Wgt}_k (\hat G)\too {\rm Fbs}_k(\hat G)
 \eqno £6.5.7\phantom{.}$$
which is the inverse of the right-hand map in~£6.5.3.

\bigskip
\noindent
{\bf £7. Vertex, sources and multiplicity modules}
\bigskip
£7.1. Let $\hat G$ be again a $p\-$solvable finite $k^*\-$group and choose a {\it polarization\/}~$\omega\,;$
then, it follows from Theorem~£6.5 above that any simple $k_*\hat G\-$mo-dule~$M$ determines a $G\-$conjugacy
class $\,\overline{\!(R,Y)\!}\,$ of {\it weights\/} of $\hat G$ and in this section we discuss the relationship
between this $G\-$conjugacy class  $\,\overline{\!(R,Y)\!}\,$ and the $G\-$conjugacy class  of the triples formed by 
a {\it vertex\/} $Q\,,$ a {\it $Q\-$source\/} $E$ and a {\it multiplicity $\skew3\hat{\bar N}_G (Q)_E\-$module\/}
$V$ of $M$ (cf.~£2.5). 
\eject

\medskip
£7.2. Actually, $M$ also determines a $G\-$conjugacy class of {\it Fitting block sequences\/} 
$\B =\{(\hat G_n,b_n)\}_{n\in \Bbb N}\,;$ let us denote by $\{O^n_{\,\eta_n}\}_{n\in \Bbb N}$ the corresponding
sequence of Fitting pointed groups $O^n_{\,\eta_n}$ on $k_*\hat G_n b_n\,,$ so that for any $n\in \Bbb N$ we have
$$\hat G_{n+1} = \hat F_{k^*\hat G}(O^n_{\,\eta_n})
\eqno £7.2.1;$$
let $P^n_{\gamma_n}$ be a maximal local pointed group on $k_*\hat G_n$ containing $O^n_{\,\eta_n}$ and 
$(S_n, \hat L_n)$ a {\it $P^n\-$source\/} pair for $k_*\hat G_n b_n$ (cf.~£4.4); note that, according to statement £2.11.2
above and to [12,~Lemma~3.10], up to a suitable identification, $\bar P^n = P^n/O^n$ is a Sylow $p\-$subgroup of~$G_{n+1}$ and therefore there is $\bar x\in G_{n+1}$ such that $(\bar P^n)^{\bar x}$ contains $P_{n+1}\,;$ thus, since the sequence $\B$ ``stabilizes'', up to finite number of choices we may assume that $\bar P_n$ contains $P_{n+1}$ for any~$n\in \Bbb N\,.$

\medskip
£7.3. Moreover, from Theorem~£6.5 we still obtain a {\it weight sequence\/} 
$\{\,\overline{\!(R^n,Y^n)\!}\,\}_{n\in \Bbb N}$  $\omega\-$associated to $\B$ starting on  
$\,\overline{\!(R,Y)\!} = \,\overline{\!(R^0,Y^0)\!}\,\,,$ and from Corollary~£5.4 we get a 
{\it simple sequence\/}  $\{M_n\}_{n\in \Bbb N}$ $\omega\-$associated to $M$ of simple   $k_*\hat G_n\-$modules $M_n$ inductively defined by $M_0 = M$ and, denoting by $\varphi_n$ the Brauer character of~$M_n\,,$ by $\varphi_{n+1} = \Gamma_{(\hat G_n,b_n)}^\omega(\varphi_n)$ for any $n\in \Bbb N\,;$
explicitly, the Morita equivalence between $k_*\hat G_nb_n$ and $(k_*\hat G_n)_{\gamma_n}$ determines a simple $(k_*\hat G_n)_{\gamma_n}\-$module  $(M_n)_{\gamma_n}$ and let us
set 
$${\rm End}_k(M_n)_{\gamma_n} = {\rm End}_k\big((M_n)_{\gamma_n}\big)
\eqno £7.3.1;$$
then, choosing an embedding (cf. statement~£4.4.2)
$$e_{\gamma_n} : (k_*\hat G_n)_{\gamma_n}\too S_n\otimes_k k_*\hat L_n
\eqno £7.3.2,$$
the restriction {\it via\/} the embedding~£4.5.2 determines a simple $k_*\hat L_n\-$module 
$\dot M_n$ which becomes a simple 
$k_*\hat F_{k_*\hat L}(O^n_{\,\dot\eta_n})\-$module (cf.~isomorphism~£4.6.1); finally, we may assume that we have (cf.~isomorphism~£4.5.4)
$$M_{n+1}=  {\rm Res}_{\Phi^\omega_{S_n}(O^n_{\,\dot\eta_n})\circ 
\hat F_{\tilde e_{\gamma_n}} (O^n_{\,\eta_n})} (\dot M_n)
\eqno £7.3.3.$$

\medskip
£7.4. For any $n\in \Bbb N\,,$ let $Q^n$ be a {\it vertex\/} and $E_n$ a {\it $Q^n\-$source\/} of $M_n\,;$ denoting by 
$Q^n_{\,E_n}$ the corresponding local pointed group on ${\rm End}_k (M_n)\,,$ it is clear that there is a local
point $\delta_n$ of $Q^n$ on $k_*\hat G_n b_n$ which has a nonzero image in 
$\big({\rm End}_k (M_n)\big)(Q^n_{\,E_n})\,;$ thus, we may assume that $P^n_{\,\gamma_n}$
contains $Q^n_{\,\delta_n}$ and it follows easily from [9~Proposition~1.6] applied to ${\rm End}_k (M)$ that 
 $Q^n_{\,\delta_n}$ is a radical pointed group on $k_*\hat G_n b_n\,,$ so that  $Q^n_{\,\delta_n}$ contains $O^n_{\,\eta_n}$ (cf.~Proposition~£4.2).

\bigskip
\noindent
{\bf Lemma~£7.5.} {\it With the notation above and up to a suitable identification, the quotient $\bar Q^n = Q^n/O^n$ is a vertex of $M_{n+1}\,.$ In particular, there is  $\bar x\in G_{n+1}$ such that $(\bar Q^n)^{\bar x} = Q^{n+1}\,.$\/}
\eject
\medskip
\noindent
{\bf Proof:} Since $Q^n$ is a vertex of 
$M_n$ and  $P^n_{\,\gamma_n}$ contains $Q^n_{\,\delta_n}\,,$ we have 
$${\rm End}_k (M_n)_{\gamma_n}(Q^n)\not= \{0\}
\eqno £7.5.1\phantom{.}$$
and therefore we still have
$$\big({\rm End}_k (M_{n+1})\big)(\bar Q^n)\cong \big({\rm End}_k (\dot M_{n})\big)(\bar Q^n)\not= \{0\}
\eqno £7.5.2,$$
so that  there is  $\bar x\in G_{n+1}$ such that $(\bar Q^n)^{\bar x}\i  Q^{n+1}\,.$

\smallskip
 Conversely, denoting by  $\check Q^{n+1}$ the converse image of $Q^{n +1}$ in $P^n$ and by $x$ a lifting of $\bar x$ to $N_{G_n} (O^n_{\,\eta_n})\,,$ we have $(Q^n)^{x}\i \check Q^{n+1}\,;$ moreover, since $Q^{n+1}\i P^n\,,$ it is clear that $S_n(\check Q^{n+1})\not= \{0\}$ and therefore we get
$$\eqalign{\big(S_n\otimes_k {\rm End}_k (\dot M_{n})&\big)(\check Q^{n+1})\cr
&\cong S_n(\check Q^{n+1})\otimes_k \big({\rm End}_k (M_{n+1})\big)(Q^{n +1})\not= \{0\}\cr}
\eqno £7.5.3\phantom{.}$$
which implies that ${\rm End}_k (M_n)_{\gamma_n}(\check Q^{n +1})\not= \{0\}$ and 
{\it a frotiori\/} that
$$\big({\rm End}_k (M_n)\big)(\check Q^{n +1})\not= \{0\}
\eqno £7.5.4;$$
thus, $\check Q^{n+1}$ is contained in a  vertex of $M_n$ and thus we have $(\bar Q^n)^{\bar x} = Q^{n+1}\,.$

\medskip
£7.6.  Once again, since the sequence $\B$ ``stabilizes'', up to finite number of choices we may assume that $\bar Q^n = Q^{n+1}$ 
for any $n\in \Bbb N\,;$ at this point, setting $Q= Q^0\,,$ these equalities determine  group homomorphisms
$$\rho_n : Q\too Q^n\i P^n
\eqno £7.6.1\phantom{.}$$
and therefore we have a Dade $Q\-$algebra ${\rm Res}_{\rho_n}(S_n)$ for any $n\in \Bbb N\,;$ then, since all but a finite number of these Dade $Q\-$algebras are isomorphic to $k\,,$ it makes sense to define the Dade $Q\-$algebra
$$T = \bigotimes_{n\in \Bbb N} {\rm Res}_{\rho_n}(S_n)
\eqno £7.6.2;$$
more generally, we denote by $T_h$ the Dade $Q^h\-$algebra obtained from the tensor product
$\bigotimes_{n\in \Bbb N} {\rm Res}_{\rho_{h +n}}(S_{h +n})$ for any $h\in \Bbb N\,.$
We are ready to describe a vertex $Q$ and a $Q\-$source $E = E_0$ of $M\,;$ as it could be expected, our parameterizations agree with the correspondence exhibited by Okuyama in~[8].

\bigskip
\noindent
{\bf Proposition~£7.7.} {\it With the notation and the choice above, $R$ is a vertex of~$M$ and, assuming
that $Q = R\,,$ an $R\-$source $E$ of $M$ is determined by an $R\-$interior algebra embedding
$\,{\rm End}_k (E)\to T\,.$\/}

\medskip
\noindent
{\bf Proof:} We argue by induction on the ``length to stabilization'' of $\B\,;$ if this length is zero then
we have $Q = \{1\} = R$ and $T\cong k\,,$ so that everything is clear. Otherwise, considering the $k^*\-$group
$\hat G_1\,,$ the simple $k_*\hat G_1\-$module~$M_1\,,$ the {\it Fitting block
sequence\/} $\B_1 = \{(\hat G_{1 +n},b_{1 +n})\}_{n\in \Bbb N}$ and the $G_1\-$conjugacy\break
\eject
\noindent
class $\,\overline{\!(R^1,Y^1)\!}\,$ of {\it weights\/} of $\hat G_1$ determined by $M_1\,,$ it follows from the induction hypothesis
that we may assume that $Q^1 = R^1$ and that an $R^1\-$source $E_1$ of $M_1$ is determined by an $R^1\-$interior algebra embedding $\,{\rm End}_k (E_1)\to T_1\,.$

\smallskip
But, it follows from Lemma~£7.5 that we may assume that $Q$ is the converse image of $Q^1$ in $P\,,$ and from 
Corollary~£5.7 that $R$ is $G\-$conjugate to the converse image of $R^1\,;$ consequently, $R$ is also a vertex of 
$M$ and we may assume that $R = Q\,.$ Moreover, we clearly have a $P^0\-$interior algebra embedding (cf.~£7.3)
$${\rm End}_k(M_0)_{\gamma_0}\too S_0\otimes_k {\rm End}_k (\dot M_{0})
\eqno £7.7.2\phantom{.}$$
and therefore, since we have $Q^0_{\delta_0}\i P^0_{\,\gamma_0}$ (cf.~£7.4), we can choose an $R\-$source 
$E$ of~$M_0 = M$ such that embedding~£7.7.2 determines an $R\-$interior algebra embedding
$${\rm End}_k (E)\too {\rm Res}_{\rho_0} (S_0)\otimes_k {\rm End}_k (E_1) \too 
{\rm Res}_{\rho_0} (S_0)\otimes_k T_1 = T
\eqno £7.7.3.$$
We are done.

\medskip
£7.8. From now on, we assume that $Q^n = R^n$ for any $n\in \Bbb N$ and, as in~Lemma~£2.16 above, we denote by
$R^n_{\, E_n}$ the corresponding local pointed group on~${\rm End}_k (M_n)\,;$   
let us consider a {\it multiplicity  $k_*\skew3\hat{\bar N}_{G_n} (R^n_{\, E_n})\-$module\/} $V_n$ of $M_n\,;$ since
 by Proposition~£7.7 we already know that  ${\rm End}_k (E_n)$ is a Dade $R^n\-$algebra, it follows from Lemma~£2.16
 that there exists a $k^*\-$group isomorphism
$$\skew3\hat{\bar N}_{G_n} (R^n_{\, E_n})\cong \bar  N_{\hat G_n}(R^n_{\, E_n})
\eqno £7.8.1\phantom{.}$$
which, according to the $k^*\-$group homomorphism £2.8.1, depends on the choice of a splitting for the $k^*\-$group
(cf.~£2.9 and~Proposition~£7.7)
$$\hat F_{{\rm End}_k (M_n)}(R^n_{\, E_n}) \cong \hat F_{{\rm End}_k (E_n)}(R^n_{\, E_n}) \cong \hat F_{T_n}(R^n) 
\eqno £7.8.2\phantom{.}$$
and indeed, from our choice of the {\it polarization\/} $\omega\,,$ we have the splitting
$$\omega_{(R^n,T_n)} : \hat F_{T_n}(R^n) \too k^*
\eqno £7.8.3.$$

\medskip
£7.9. On the other hand, for any $n\in \Bbb N\,,$ let $W_n$ be the  restriction to the stabilizer 
$\skew3\hat {\bar N}_{G_n}  (R^n_{\,\delta_n})_{E_n}$  in $\skew3\hat {\bar N}_{G_n} (R^n_{\,\delta_n})$ of the isomorphism class of~$E_n\,,$ of a {\it multiplicity $k_*\skew3\hat {\bar N}_{G_n} 
(R^n_{\,\delta_n})\-$module\/} of $R^n_{\,\delta_n}$ (cf.~£2.5);
more explicitly, denoting by $\bar b(\delta_n)$ the block of $\bar C_{\hat G_n} (R^n)$ determined by $\delta_n\,,$
since $R^n_{\,\delta_n}$ is a radical pointed group on $k_*\hat G_n\,,$ we have (cf.~£2.11.1 and~£3.1)
$$k_*\bar C_{\hat G_n} (R^n)\bar b(\delta_n)\cong {\rm End}_k (W_n)
\eqno £7.9.1.$$ 
Consequently, since $\big({\rm End}_k (M_n)\big)(R^n_{\,E_n}) \cong {\rm End}_k (V_n)$ 
 has a $\bar C_{\hat G}(R^n)\-$interior algebra structure, it makes sense to consider
 $\bar b (\delta_n)\.\big({\rm End}_k (M_n)\big)(R^n_{\,E_n})\.\bar b (\delta_n)$
 as a $\skew3\hat {\bar N}_{G_n}(R^n_{\,E_n})_{\delta_n}\-$interior algebra and, from
  the structural homomorphism above, we get an injective $k\-$algebra homomorphism
 $$(k_*\hat G)(R^n_{\,\delta_n})\too \bar b (\delta_n)\.\big({\rm End}_k (M_n)\big)(R^n_{\,E_n})\.\bar b (\delta_n)
 \eqno £7.9.2;$$
  then,  denoting by  $\skew3\hat {\bar N}_{G_n}  (R^n_{\,E_n})_{\delta_n}$ the stabilizer of $\delta_n$ in $\skew3\hat {\bar N}_{G_n}  (R^n_{\,E_n})$ and setting
$$\skew8\widehat{\bar N}^{^{E_n}}_n = \skew3\hat {\bar N}_{G_n}(R^n_{\,E_n})_{\delta_n} *                            
\big(\skew3\hat {\bar N}_{G_n}  (R^n_{\,\delta_n})_{E_n}\big)^\circ 
\eqno £7.9.3,$$
 it follows from [9,~Proposition~2.1] that, for a suitable {\it primitive\/} $\skew8\widehat{\bar N}^{^{E_n}}_n\-$interior algebra $B_n\,,$ we have a $\skew3\hat {\bar N}_{G_n}(R^n_{\,E_n})\-$interior algebra isomorphism
$${\rm End}_k (V_n)\cong {\rm Ind}_{\skew3\hat {\bar N}_{G_n}(R^n_{\,E_n})_{\delta_n}}^{\skew3\hat {\bar N}_{G_n}(R^n_{\,E_n})} \big(k_*\bar C_{\hat G_n} (R^n) \bar b(\delta_n)\otimes_k B_n\big)
\eqno £7.9.4;$$
actually,  it is easily checked that the subgroup
 $$\bar C_G(R^n)\cong \bar C_{\hat G}(R^n) * \bar C_{\hat G}(R^n)^\circ \i 
 \skew8\widehat{\bar N}^{^{E_n}}_n
 \eqno £7.9.5\phantom{.}$$
 has a trivial image in $B_n$ so that, up to an obvious identification, $B_n$ becomes an 
 $\skew8\widehat{\bar N}^{^{E_n}}_n/\bar C_G (R^n) \-$interior algebra.

\medskip
£7.10. Similarly, denoting by $U_n$ a {\it simple projective $k_*\bar N_{\hat G_n} (R^n)\-$module\/}
which restricted to $N_{\hat G_n} (R^n)$ belongs to the isomorphism class $Y^n\,,$ it follows from [16,~Proposition~3.2] applied to the {\it primitive\/} $N_{\hat G_n}(R^n)\-$interior algebra ${\rm End}_k (U_n)$ that, setting
$$\skew3\hat{\bar N}_n = \bar N_{\hat G_n} (R^n_{\,\delta_n}) * \skew3\hat {\bar N}_{G_n}  (R^n_{\,\delta_n})^\circ
\eqno £7.10.1,$$
for a suitable {\it primitive\/} $\skew3\hat{\bar N}_n\-$interior algebra $D_n$  we have 
$${\rm End}_k (U_n)\cong {\rm Ind}_{\bar N_{\hat G_n} (R^n_{\,\delta_n})}^{ N_{\hat G_n} 
(R^n)} \big(k_*\bar C_{\hat G_n} (R^n)\bar b(\delta_n)\otimes_k D_n\big)
\eqno £7.10.2;$$
actually, it is clear from its very definition that $D_n$ becomes a $\skew3\hat{\bar N}_n/\bar C_{G_n} (R^n)\-$ interior algebra and note that, according to homomorphism~£2.8.1, we have a {\it canonical\/} $k^*\-$group isomorphism
$$\skew3\hat{\bar N}_n/\bar C_{G_n} (R^n)\cong \hat F_{k_*\hat G_n} (R^n_{\,\delta_n})
\eqno £7.10.3.$$
In order to relate $D_n$ with $U_{n+1}\,,$ we have to consider the following $k^*\-$group isomorphism.
\eject

\bigskip
\noindent
{\bf Proposition~£7.11.} {\it  With the notation and the choice above, there is a $k^*\-$group isomorphism
$${}^\omega\phi_n : \hat F_{k_*\hat G_n} (R^n_{\,\delta_n})\cong \bar N_{\hat G_{n+1}}(R^{n +1})
\eqno £7.11.1\phantom{.}$$
such that, for any $P^n\-$interior algebra  embedding 
$$e_{\gamma_n} : (k_*\hat G_n)_{\gamma_n}\too S_n \otimes_k k_*\hat L_n
\eqno £7.11.2\,,$$
  we have the commutative diagram
$$\matrix{\hat F_{k_*\hat G_n} (R^n_{\,\delta_n}) \hskip-5pt &\cong &\hskip-5pt
\hat F_{S_n\otimes_k k_*\hat L_n}(R^n_{\,S_n\times \dot\delta_n})\hskip-5pt
 &\cong &\hskip-5pt \hat F_{k_*\hat L_n} (R^n_{\,\dot\delta_n}) \cr
&&\phantom{\Big\uparrow}&&\wr\Vert\cr
\hskip-30pt{\scriptstyle {}^\omega\phi_n}\hskip3pt\wr\!\Vert&&&&\bar N_{\hat L_n} (R^n)\cr
&&\phantom{\Big\uparrow} &&\wr\Vert\cr
\bar N_{\hat G_{n+1}}(R^{n +1})\hskip-5pt &\cong &
\hskip-5pt\bar N_{\hat F_{S_n\otimes_k k_*\hat L_n}(O^n_{S_n\times \dot\eta_n})}(\bar R^n)\hskip-5pt
&\cong &\hskip-5pt \bar N_{\hat F_{k_*\hat L_n}(O^n_{\,\dot\eta_n})}(\bar R^n)
\cr}
\eqno £7.11.3.$$\/}

\par
\noindent
{\bf Proof:} Choosing a $P^n\-$interior algebra  embedding $e_\gamma\,\colon (k_*\hat G)_\gamma\to S
\otimes_k k_*\hat L\,,$ 
we have $k^*\-$group isomorphisms (cf.~£5.5.1)
$$\hat F_{k_*\hat G_n} (R^n_{\,\delta_n})\buildrel \hat F_{\tilde e_{\gamma_n}}(R^n_{\,\delta_n})\over\cong 
\hat F_{S_n\otimes_k k_*\hat L_n}(R^n_{\,S_n\times \dot\delta_n})\buildrel \Phi^\omega_{S_n}(R^n_{\,\dot\delta_n})\over  \cong \hat F_{k_*\hat L_n} (R^n_{\,\dot\delta_n}) 
\eqno £7.11.4;$$
similarly, the $k^*\-$group isomorphisms~£4.5.4 applied to $O^n_{\,\eta_n}$ yield
$$\hat F_{k_*\hat G_n} (O^n_{\,\eta_n})\buildrel \hat F_{\tilde e_{\gamma_n}}(O^n_{\,\eta_n})\over\cong 
\hat F_{S_n\otimes_k k_*\hat L_n}(O^n_{S_n\times \dot\eta_n})\buildrel \Phi^\omega_{S_n}(O^n_{\,\dot\eta_n})\over 
\cong \hat F_{k_*\hat L_n} (O^n_{\,\dot\eta_n}) 
\eqno £7.11.5;$$
furthermore, setting $\bar R^n = R^n/O^n\,,$  we have {\it canonical\/} 
$k^*\-$group isomorphisms (cf.~isomorphisms £4.6.2 and~£5.5.2)
$$\hat F_{k_*\hat L_n} (R^n_{\,\dot\delta_n})\cong \bar N_{\hat L_n} (R^n)\cong 
\bar N_{\hat F_{k_*\hat L_n}(O^n_{\,\dot\eta_n})}(\bar R^n)
\eqno £7.11.6.$$
Now, it is clear that the commutativity of the corresponding diagram above defines a $k^*\-$group isomorphism
 ${}^\omega\phi_n\,.$

 \smallskip
We claim that this $k^*\-$group isomorphism does not depend on the choice of~$e_{\gamma_n}\,;$ indeed, for another choice $e'_{\gamma_n}$  of this embedding, it follows from Proposition~£5.2 that there is $\tilde\sigma_n\in 
 {\rm Out}_{P^n}\big((k_*\hat G_n)_{\gamma_n}\big)$ fulfilling 
 $$\tilde e'_{\gamma_n} = \tilde e_{\gamma_n}\circ \tilde\sigma_n = (\widetilde{\rm id}_{S_n} \otimes \dot{\tilde\sigma_n}) 
 \circ \tilde e_{\gamma_n}
 \eqno £7.11.7\phantom{.}$$
 and therefore, with obvious notation, we get the following commutative diagrams
 $$\matrix{\hat F_{k_*\hat G_n} (R^n_{\,\delta_n})\hskip-10pt
 &\buildrel \hat F_{\tilde e'_{\gamma_n}}(R^n_{\,\delta_n})\over\cong 
 &\hat F_{S_n\otimes_k k_*\hat L_n}(R^n_{\,S_n\times \dot\delta_n})\hskip-10pt
 &\buildrel \Phi^\omega_{S_n}  (R^n_{\,\dot\delta_n})\over \cong 
 &\hat F_{k_*\hat L_n} (R^n_{\,\dot\delta_n}) \cr
\Vert&\phantom{\big\uparrow}&\hskip-40pt{\scriptstyle \hat F_{\tilde{\rm id}_S \otimes 
\dot{\tilde\sigma}}(R^n_{\,S_n\times \dot\delta_n})}\hskip4pt\wr\!\Vert&&\hskip-40pt{\scriptstyle \hat F_{\dot{\tilde\sigma}}(R^n_{\,\dot\delta_n})}
\hskip4pt\wr\!\Vert\cr
\hat F_{k_*\hat G_n} (R^n_{\,\delta_n})\hskip-10pt&\buildrel \hat F_{\tilde e_{\gamma_n}}(R^n_{\,\delta_n})\over\cong &\hat F_{S_n\otimes_k k_*\hat L_n}(R^n_{\,S_n\times \dot\delta_n})\hskip-10pt&\buildrel \Phi^\omega_{S_n}(R^n_{\,\dot\delta_n})\over  \cong 
&\hat F_{k_*\hat L_n} (R^n_{\,\dot\delta_n}) \cr}
\eqno £7.11.8\phantom{.}$$

\smallskip
$$\matrix{\hat F_{k_*\hat G_n} (O^n_{\,\eta_n})\hskip-10pt&\buildrel \hat F_{\tilde e'_{\gamma_n}}(O^n_{\,\eta_n})\over\cong &\hat F_{S_n\otimes_k k_*\hat L_n}(O^n_{S_n\times \dot\eta_n})\hskip-10pt&\buildrel \Phi^\omega_{S_n}(O^n_{\,\dot\eta_n})\over \cong &\hat F_{k_*\hat L_n} (O^n_{\,\dot\eta_n})\cr
\Vert&\phantom{\Big\uparrow}&\hskip-40pt{\scriptstyle \hat F_{\tilde{\rm id}_S \otimes 
\dot{\tilde\sigma}}(O^n_{\,S_n\times \dot\eta_n})}\hskip4pt\wr\!\Vert&&\hskip-40pt{\scriptstyle 
\hat F_{\dot{\tilde\sigma}}(O^n_{\,\dot\eta_n})}\hskip4pt\wr\!\Vert\cr
\hat F_{k_*\hat G_n} (O^n_{\,\eta_n})\hskip-10pt&\buildrel \hat F_{\tilde e_{\gamma_n}}(O^n_{\,\eta_n})\over\cong &
\hat F_{S_n\otimes_k k_*\hat L_n}(O^n_{S_n\times \dot\eta_n})\hskip-10pt&\buildrel \Phi^\omega_{S_n}
(O^n_{\,\dot\eta_n})\over  \cong &\hat F_{k_*\hat L_n} (O^n_{\,\dot\eta_n}) \cr}
\eqno £7.11.9.$$

\smallskip
Now, the commutativity of the corresponding diagram above follows from these commutative diagrams and from the  {\it naturality\/} of the right-hand vertical isomorphisms in the diagram. We are done.

\bigskip
\noindent
{\bf Corollary~£7.12.} {\it  With the notation and the choice above, we have an\break
$\hat F_{k_*\hat G_n} (R^n_{\,\delta_n})\-$interior algebra isomorphism
$$D_n\cong  {\rm Res}_{\,{}^\omega\phi_n}\big({\rm End}_k (U_{n +1})\big)
\eqno £7.12.1.$$\/}

\par
\noindent
{\bf Proof:} Since ${\rm End}_k (U_n)$ is actually isomorphic to a {\it block of defect zero\/}
of the $k^*\-$group  $\bar N_{\hat G_n}  (R^n)$ (cf.~£7.10), it follows from [16,~Theorem~3.7] and from  isomorphism~£7.10.2 above that $D_n$ is isomorphic  to a {\it block of defect zero\/} of the  $k^*\-$group  
$\hat F_{k_*\hat G_n} (R^n_{\,\delta_n})$ and then isomorphism~£7.12.1 easily follows from the commutativity of 
diagram~£7.11.3 above.

\medskip
£7.13. On the other hand, according to Proposition~£7.7, we have $R^n\-$ and $R^{n+1}\-$interior algebra embeddings
$${\rm End}_k (E_n) \too T_n\qq {\rm End}_k (E_{n+1}) \too T_{n+1}
\eqno £7.13.2;$$
but, denoting by $t_n\,,$ $t_{n+1}$ and $s_n$ the respective {\it similarity\/} classes in the {\it Dade group\/}
$\D_k (R^n)$ of $T_n\,,$ of the restriction to $R^n$ of $T_{n+1}\,,$ and of $S_n$ (cf.~£2.13), we clearly have 
$t_n = s_n + t_{n+1}$ (cf.~£7.6) and therefore any automorphism of $R^n$ stabilizing $s_n$ stabilizes $t_n$ 
if and only if it stabilizes  $t_{n+1}\,;$ moreover, it follows from the inclusion in~£4.5.3 that 
$F_{k_*\hat G_n} (R^n_{\,\delta_n})$ stabilizes $s_n\,.$ Consequently, since $t_n$ and $t_{n+1}$ 
respectively determine the isomorphism classes of $E_n$
and of the restriction to $R^n$ of $E_{n+1}\,,$ the stabilizers in $F_{k_*\hat G_n} (R^n_{\,\delta_n})$ of these
isomorphism classes coincide with each other, and therefore we have a  {\it canonical\/}  surjective homomorphism
$$\nu_n : \bar N_{G_n}   (R^n_{\,\delta_n})_{E_n} = \bar N_{G_n}  (R^n_{\,E_n})_{\delta_n}\too \bar N_{G_{n+1}}   
(R^{n+1}_{\; E_{n+1}})
\eqno £7.13.3.$$
\eject

\bigskip
\noindent
{\bf Proposition~£7.14.} {\it With the notation and the choice above,  there are a $k^*\-$group homomorphism and 
 a $\skew8\widehat{\bar N}^{^{E_n}}_n\-$interior algebra isomorphism
$$\eqalign{\hat\nu_n &: \skew8\widehat{\bar N}^{^{E_n}}_n\too \skew3\hat {\bar N}_{G_{n+1}}   (R^{n+1}_{\; E_{n+1}})\cr
f_n &:  B_n \cong {\rm Res}_{\nu_n}\Big(\big({\rm End}_k (M_{n+1})\big)(R^{n+1}_{\; E_{n+1}})\Big)\cr}
\eqno £7.14.1\phantom{.}$$
 such that, for any $P^n\-$interior algebra  embedding 
 $$e_{\gamma_n} : (k_*\hat G_n)_{\gamma_n} \too  S_n \otimes_k k_*\hat L_n
 \eqno £7.14.2,$$
   we have the commutative diagram 
$$\matrix{B_n &\cong &\big({\rm End}_k (\dot M_n)\big)(R^n_{\,\dot E_n})\cr
\hskip-20pt {\scriptstyle f_n}\hskip4pt\wr\!\Vert&\phantom{\big\downarrow}&\wr\Vert\cr
\big({\rm End}_k (M_{n+1})\big)(R^{n+1}_{\; E_{n+1}})&\cong &\big({\rm End}_k 
(\dot{\bar M}_n)\big)
(\bar R^n_{\,\dot{\bar E}_n})\cr}
\eqno £7.14.3.$$\/}

\par
\noindent
{\bf Proof:} We have a structural injective $k\-$algebra homomorphism (cf.~£7.3.1)
 $$(k_*\hat G)_{\gamma_n}(R^n_{\,\delta_n})\too \big({\rm End}_k (M_n)\big)(R^n_{\,E_n})_{\gamma_n}
 \eqno £7.14.4\phantom{.}$$
 and, as in~£7.9 above, denoting by $C_n$ the centralizer in ${\rm End}_k (M_n)_{\gamma_n}(R^n_{\,E_n})$ of the image of~$(k_*\hat G_n )_{\gamma_n}
 (R^n_{\,\delta_n})\,,$ it follows from [9,~Proposition~2.1] that
 we have a $k\-$algebra isomorphism
 $${\rm End}_k (M_n)_{\gamma_n}(R^n_{\,E_n})
 \cong (k_*\hat G)_{\gamma_n}(R^n)\otimes_k C_n
 \eqno £7.14.5.$$
 Moreover, always according to [9,~Proposition~2.1], the obvious commutative diagram 
 $$\matrix{(k_*\hat G)_{\gamma_n} (R^n_{\,\delta_n})&
 \too& (k_*\hat G)(R^n_{\,\delta_n})\cr
\big\downarrow&\phantom{\Big\downarrow}&\big\downarrow \cr
 {\rm End}_k (M_n)_{\gamma_n}(R^n_{\,E_n})&
 \too&\bar b (\delta_n)\.\big({\rm End}_k (M_n)\big)(R^n_{\,E_n})\.\bar b (\delta_n)\cr}
 \eqno £7.14.6\phantom{.}$$
induce a {\it canonical\/} $k\-$algebra isomorphism $C_n\cong B_n$ which allows us to identify both centralizers.

 \smallskip
 Choosing a  $P^n\-$interior algebra  embedding (cf.~statement~£4.4.2)
 $$e_{\gamma_n} : (k_*\hat G_n)_{\gamma_n}\too S_n \otimes_k k_*\hat L_n
 \eqno £7.14.7,$$
 note that $(k_*\hat L_n) (R^n_{\,\dot\delta_n})\cong k$ since $R^n_{\,\delta_n}$ and
 $R^n_{\,\dot\delta_n}$ are {\it radical\/} and therefore  they are {\it selfcentralizing\/};
 then, the corresponding commutative diagram
 $$\matrix{(k_*\hat G_n)_{\gamma_n} (R^n_{\,\delta_n})&\too
 & S_n(R^n)  \cr
 \big\downarrow&\phantom{\Big\downarrow}&\big\downarrow  \cr
  \big({\rm End}_k (M_n)\big)_{\gamma_n}(R^n_{\,E_n})&
 \too&  S_n(R^n)  \otimes_k  \big({\rm End}_k (\dot M_n)\big)(R^n_{\,\dot E_n})\cr}
 \eqno £7.14.8\phantom{.}$$
 \eject
 \noindent
 and the argument above yield the top $k\-$algebra isomorphism
 $$B_n\cong \big({\rm End}_k (\dot M_n)\big)(R^n_{\,\dot E_n})
 \eqno £7.14.9\phantom{.}$$
 in the diagram £7.14.3 above, and the $k^*\-$group isomorphism
$$\skew8\widehat{\bar N}^{^{E_n}}_n/C_G (R^n)\cong 
\skew3\hat{\bar N}_{L_n} (R^n_{\,\dot E_n})
 \eqno £7.14.10.$$

 \smallskip
Moreover, according to~£7.3, we get
 $${\rm End}_k (\dot M_n) \cong {\rm End}_k (\skew3\dot{\bar M}_n) \cong 
 {\rm End}_k (M_{n+1})
 \eqno £7.14.11\phantom{.}$$
 with the interior structures coming from the $k^*\-$group isomorphisms (cf.~£4.5.4
 and~£4.6.1)
 $$\hat L_n/O^n\cong \hat F_{k_*\hat L_n}(O^n_{\,\dot\eta_n})\cong  \hat G_{n+1}
\eqno £7.14.12;$$
consequently, we still get
$$\eqalign{\big({\rm End}_k (\dot M_n)\big)(R^n_{\,\dot E_n}) \cong 
\big({\rm End}_k (\skew3\dot{\bar M}_n)\big)&(R^n_{\,\skew3\dot{\bar E}_n})\cr
& \cong 
\big( {\rm End}_k (M_{n+1})\big)(R^{n+1}_{\; E_{n+1}})\cr}
\eqno £7.14.13\phantom{.}$$
 with the interior structures coming from the $k\-$group isomorphisms 
  $$\skew3\hat{\bar N}_{L_n} (R^n_{\,\dot E_n})\cong \skew3\hat{\bar N}_{ F_{k_*\hat L_n}(O^n_{\,\dot\eta_n})}(\bar R^n_{\,\skew3\dot{\bar E}_n})\cong  
  \skew3\hat{\bar N}_{G_{n+1}} (R^{n+1}_{\; E_{n+1}})
\eqno £7.14.14.$$

\smallskip
Finally, for a particular choice of $e_{\gamma_n}\,,$ the commutativity of diagram £7.14.3 induce the isomorphism $f_n\,,$  and the $k^*\-$group isomorphisms~£7.14.10  and £7.14.14 determine the $k^*\-$group homomorphism $\hat\nu_n\,.$ Once again $f_n$ and~$\hat\nu_n$ do not depend on the choice of~$e_{\gamma_n}\,;$ indeed, for another choice $e'_{\gamma_n}$
 of this embedding, it follows from Proposition~£5.2 that there is $\tilde\sigma_n\in 
 {\rm Out}_{P^n}\big((k_*\hat G_n)_{\gamma_n}\big)$ fulfilling 
 $$\tilde e'_{\gamma_n} = \tilde e_{\gamma_n}\circ \tilde\sigma_n = (\widetilde{\rm id}_{S_n} \otimes \dot{\tilde\sigma_n}) 
 \circ \tilde e_{\gamma_n}
 \eqno £7.14.15\phantom{.}$$
 and  we get commutative diagrams as above.

\medskip
£7.15. We are ready to describe the {\it multiplicity 
$\skew3\hat{\bar N}_G   (R_{\, E})\-$module\/} of $M\,;$ first of all, from the very definition of
$N_G(\B)$ (cf.~£6.2), we get a sequence of $k^*\-$groups 
$\skew6\widehat{N}_G^{^n} (\B)\,,$ with the same $k^*\-$quotient $N_G (\B)\,,$ and of 
$k^*\-$group homomorphisms $\hat\mu_n\,\colon\skew6\widehat{N}^{^n}_G (\B)
\to \hat G_n$ inductively defined as follows; the $k^*\-$group 
$\skew6\widehat{N}_G^{^0} (\B)$ is just the converse image of $N_G (\B)$ in $\hat G$
 and $\hat\mu_0$ the inclusion map;  then, for any $n\ge 1\,,$ arguing by induction on $n$ it is easily checked that the image of $\hat\mu_{n -1}$ normalizes the pointed group
  $O^{n-1}_{\;\eta_{n-1}}$ on $k_*\hat G_{n-1}$ and therefore $\hat\mu_{n-1}$
  induces a group homomorphism $\mu_{n-1}$ from $N_G (\B)$ to 
  $N_{G_{n-1}} (O^{n-1}_{\;\eta_{n-1}})\,;$\break 
  \eject
  \noindent
   since we have (cf.~statement~£2.11.2)
  $$N_{G_{n-1}} (O^{n-1}_{\;\eta_{n-1}})/O^{n-1}\.C_{G_{n-1}}(O^{n-1})
  \cong F_{k_*\hat G_{n-1}}(O^{n-1}_{\;\eta_{n-1}})
\eqno £7.15.1,$$
we define $\skew6\widehat{N}_G^{^n} (\B)$ and $\hat\mu_n$ by the following
{\it pull-back\/}
$$\matrix{N_G(\B) &\too &F_{k_*\hat G_{n-1}}(O^{n-1}_{\;\eta_{n-1}})\cr
\big\uparrow&\phantom{\Big\uparrow}&\big\uparrow\cr
\skew6\widehat{N}_{\!G}^{^n} (\B)&\buildrel \hat\mu_n\over\too &\hat G_n\cr}
\eqno £7.15.2.$$

\medskip
£7.16. Now, consider the {\it pointed vertex sequence\/} $\R =
\{R^n_{\,\delta_n}\}_{n\in \Bbb N}$ of $M$ associated  to~$\omega$  (cf.~£7.4 and~£7.8), and denote by  $N_G(\R)$ the stabilizer of~$\R$ in~$G\,;$ clearly, $N_G (\B)$ contains 
$N_G (\R)$ and we denote by $\skew6\widehat{N}_G^{^n} (\R)$ the corresponding 
$k^*\-$subgroup of $\skew6\widehat{N}_G^{^n} (\B)\,.$ Moreover, arguing by induction on $n\,,$ it is easily checked that the subgroup
$\mu_n \big(N_G(\R)\big)\i G_n$ normalizes the pointed group $R^n_{\,\delta_n}$ 
on~$k_*\hat G_n\,;$ thus, setting $\bar N_G (\R) = N_G (\R)/R\,,$ we get  
a $k^*\-$group $\skew6\widehat{\bar N}^{^{\delta_n}}_G (\R)$ from the following {\it pull-back\/}
$$\matrix{\bar N_G(\R) &\buildrel \bar\mu_n\over\too &\bar N_{G_n} (R^n_{\,\delta_n})\cr
\big\uparrow&\phantom{\Big\uparrow}&\big\uparrow\cr
\skew6\widehat{\bar N}^{^{\delta_n}} _{\!G}(\R)&\buildrel \skew3\hat{\bar\mu}^{\delta_n}_n\over\too &\skew3\hat{\bar N}_{G_n} (R^n_{\,\delta_n})\cr}
\eqno £7.16.1\phantom{.}$$
and then we have the $k_*\skew6\widehat{\bar N}^{^{\delta_n}} _{\!G}(\R)\-$module
${\rm Res}_{\skew3\hat{\bar\mu}^{\delta_n}_n}(W_n)$ for any $n\in \Bbb N$  (cf.~£7.9.1).

\medskip
£7.17. On the other hand, from the $k^*\-$group homomorphism~£2.8.1 and from
 Proposition~£7.11, for any $n\in \Bbb N$ we get the  $k^*\-$group homomorphisms
$$\matrix{\bar N_{\hat G_n}(R^n_{\,\delta_n})\star \skew3\hat{\bar N}_{G_n} 
(R^n_{\,\delta_n})^\circ &\too &\hat F_{k_*\hat G_n}(R^n_{\,\delta_n})&
\buildrel  {}^\omega\phi_n\over \cong& \bar N_{\hat G_{n+1}}(R^{n+1})\cr
\hskip-35pt{\scriptstyle \skew2\hat{\bar\mu}_n \star \skew3\hat{\bar\mu}^{\delta_n}_n}
\hskip5pt\big\uparrow&\phantom{\Big\uparrow}&&&\hskip-25pt{\scriptstyle \skew2\hat{\bar\mu}_{n+1}} \hskip5pt\big\uparrow\cr 
\skew6\widehat{\bar N}_G^{^n} (\R) \star \skew6\widehat{\bar N}^{^{\delta_n}} _{\!G}
(\R)^\circ &&&&\skew6\widehat{\bar N}_G^{^{n+1}} (\R)\cr}
 \eqno £7.17.1\phantom{.}$$
and therefore, since all the bottom $k^*\-$groups admit the same $k^*\-$quotient, we still get a $k^*\-$group isomorphism
$${}^\omega\Psi_n : \skew6\widehat{\bar N}^{^{\delta_n}} _{\!G}(\R)\cong \skew6\widehat{\bar N}_G^{^n} (\R) \star \skew6\widehat{\bar N}_G^{^{n+1}} (\R)^\circ 
\eqno £7.17.2.$$
But, for~$n$ big enough we have 
$$W_n \cong k \qq\skew6\widehat{N}_G^{^n} (\B)\cong k^*\times N _G(\B)
\eqno £7.17.3;$$
moreover, note that $\skew6\widehat{N}_G^{^0} (\B)$ coincides with the converse image
$N_{\hat G}(\B)$ of $N_G (\B)\,$ in~$\hat G$ and similarly we set $N_{\hat G}(\R) = 
\skew6\widehat{N}_G^{^0} (\R)$ and $\bar N_{\hat G}(\R) = N_{\hat G}(\R)/R\,;$\break
\eject
\noindent
in particular, the following tensor product
$${}^\omega W = \bigotimes_{n\in \Bbb N}  {\rm Res}_{\skew3\hat{\bar\mu}^{\delta_n}_n
\circ ({}^\omega\Psi_n)^{-1}}(W_n)
\eqno £7.17.4\phantom{.}$$
 makes sense and it is clearly a $k_*\bar N_{\hat G} (\R)\-$module.

\medskip
£7.18.  Finally, it follows from~£4.5.3 that $F_{k_*\hat G_n}(R^n_{\,\delta_n})$ stabilizes
the isomorphism class of ${\rm Res}_{R^n}^{P^n}(S_n)$ and therefore $\bar N_G(\R)$ stabilizes 
the {\it similarity\/} class of $T_n$ and, in particular,  the isomorphism class of $E_n$ 
for any $n\in \Bbb N$ (cf.~£7.6 and~Proposition~£7.7); hence,  we get  again
a $k^*\-$group $\skew6\widehat{\bar N}^{^{E_n}}_G (\R)$ from the following {\it pull-back\/}
$$\matrix{\bar N_G(\R) &\buildrel \bar\mu_n\over\too &\bar N_{G_n} (R^n_{\,\delta_n})_{E_n}\cr
\big\uparrow&\phantom{\Big\uparrow}&\big\uparrow\cr
\skew6\widehat{\bar N}^{^{E_n}} _{\!G}(\R)&\buildrel \skew3\hat{\bar\mu}^{^{E_n}}_n\over\too &\skew3\hat{\bar N}_{G_n} (R^n_{\,E_n})_{\delta_n}\cr}
\eqno £7.18.1.$$
 Similarly, from Proposition~£7.14, for any $n\in \Bbb N$ we get the  $k^*\-$group homomorphisms
$$\matrix{\skew3\hat{\bar N}_{G_n}(R^n_{\,E_n})_{\delta_n}\star \big(\skew3\hat{\bar N}_{G_n} 
(R^n_{\,\delta_n})_{E_n}\big)^\circ &
\buildrel  \hat\nu_n\over \too& \skew3\hat{\bar N}_{G_{n+1}}(R^{n+1}_{\,E_{n+1}})\cr
\hskip-50pt{\scriptstyle \skew3\hat{\bar\mu}^{^{E_n}}_n \star \skew3\hat{\bar\mu}^{\delta_n}_n}
\hskip5pt\big\uparrow&\phantom{\Big\uparrow}&\hskip-35pt{\scriptstyle \skew3\hat{\bar\mu}^{^{E_{n+1}}}_{n+1}} \hskip5pt\big\uparrow\cr 
\skew6\widehat{\bar N}_G^{^{E_n}} (\R) \star \skew6\widehat{\bar N}^{^{\delta_n}} _{\!G}(\R)^\circ 
&&\skew6\widehat{\bar N}_G^{^{^{E_{n+1}}}} (\R)\cr}
 \eqno £7.18.2\phantom{.}$$
and therefore, since all the bottom $k^*\-$groups admit the same $k^*\-$quotient, we still get a $k^*\-$group isomorphism
$$\Psi_n : \skew6\widehat{\bar N}^{^{\delta_n}} _{\!G}(\R)\cong \skew6\widehat{\bar N}_G^{^{E_n}} (\R) \star \skew6\widehat{\bar N}_G^{^{E_{n+1}}} (\R)^\circ 
\eqno £7.18.3;$$
thus,  the following tensor product
$$W = \bigotimes_{n\in \Bbb N}  {\rm Res}_{\skew3\hat{\bar\mu}^{\delta_n}_n
\circ (\Psi_n)^{-1}}(W_n)
\eqno £7.18.4\phantom{.}$$
 makes sense and it is clearly a $k_*\skew6\widehat{\bar N}_G^{^{E_0}} (\R)\-$module.
 As above, we set $R = R^0\,,$ $E = E_0\,,$ $V = V_0$ and $U = U_0\,.$

\bigskip
\noindent
{\bf Theorem~£7.19.} {\it With the notation and the choice above, we have natural 
$k_*\bar N_{\hat G} (R)\-$ and $k_*\skew3\hat{\bar N}_G (R_E)\-$module isomorphisms
$$U\cong {\rm Ind}_{\bar N_{\hat G}(\R)}^{\bar N_{\hat G}(R)}({}^\omega W)\quad and\quad
V\cong {\rm Ind}_{\skew6\widehat{\bar N}^{^{E}}_G(\R)}^{\skew3\hat{\bar N}_G (R_E)} (W)
\eqno £7.19.1.$$\/}

\par 
\noindent
{\bf Proof:} Once again, we can argue by induction on the ``length to stabilization'' of $\B\,.$
If this length is zero then the block $b_0$ is already of {\it defect zero\/} and therefore 
everything is trivial so that isomorphisms~£7.19.1 above are trivially true.

\smallskip
If the  ``length to stabilization'' is not zero then we consider the  {\it Fitting block sequence\/} 
 $\B_1 = \{(\hat G_{1 +n},b_{1 +n})\}_{n\in \Bbb N}$ of~$\hat G_1$ and the corresponding
{\it weight sequence\/} $\{\,\overline{\!(R^{1 +n},Y^{1 +n})\!}\,\}_{n\in \Bbb N}$ and  
{\it simple sequence\/}  $\{M_{1 +n}\}_{n\in \Bbb N}$ $\omega\-$associated to~$\B_1\,;$
{\it mutatis mutandis\/}, we consider the corresponding {\it pointed vertex sequence\/} 
$\R_1 =\{R^{n+1}_{\,\delta_{n+1}}\}_{n\in \Bbb N}$ of $M_1$ associated  to~$\omega\,,$ 
and denote by $N_{G_1}(\R_1)$ and $N_{\hat G_1}(\R_1)$ the respective stabilizers of 
$\R_1$ in $G_1$ and in $\hat G_1\,.$

\smallskip
Moreover, from the corresponding {\it pull-back\/}~£7.16.1, for any $n\ge 1$ we get a $k^*\-$group 
$\skew6\widehat{\bar N}^{^{\delta_n}} _{\!G_1}(\R_1) $ of $k^*\-$quotient 
$\bar N _{G_1}(\R_1)$ and a $k^*\-$group homomorphism
$$\skew6\widehat{\bar N}^{^{\delta_n}} _{\!G_1}(\R_1) 
\buildrel \skew3\hat{\bar\mu}^{\delta_n}_{1,n} \over\too 
\skew3\hat{\bar N}_{G_n} (R^n_{\,\delta_n})
\eqno £7.19.2,$$
so that we still get the $k_*\skew6\widehat{\bar N}^{^{\delta_n}} _{\!G_1}(\R_1)\-$module
${\rm Res}_{\skew3\hat{\bar\mu}^{\delta_n}_{1,n}}(W_n)\,.$ Analogously, for any 
$n\ge 1$ we still get the corresponding $k^*\-$group isomorphisms~£7.17.2 and~£7.18.3
$$\eqalign{{}^\omega\Psi_{1,n} : \skew6\widehat{\bar N}^{^{\delta_n}} _{\!G_1}(\R_1)&\cong \skew6\widehat{\bar N}_{G_1}^{^n} (\R_1) \star \skew6\widehat{\bar N}_{G_1}^{^{n+1}} (\R_1)^\circ\cr
\Psi_{1,n} : \skew6\widehat{\bar N}^{^{\delta_n}} _{\!G_1}(\R_1)&\cong \skew6\widehat{\bar N}_{G_1}^{^{E_n}} (\R_1) \star \skew6\widehat{\bar N}_{G_1}^{^{E_{n+1}}} (\R_1)^\circ\cr} 
\eqno £7.19.3,$$
and, once again, the following tensor products
$$\eqalign{{}^\omega W^1 &= \bigotimes_{n\ge 1} {\rm Res}_{\skew3\hat{\bar\mu}^{\delta_n}_{1,n}
\circ ({}^\omega\Psi_{1,n})^{-1}}(W_n)\cr
W^1 &= \bigotimes_{n\ge 1} {\rm Res}_{\skew3\hat{\bar\mu}^{\delta_n}_{1,n}
\circ (\Psi_{1,n})^{-1}}(W_n)\cr}
\eqno £7.19.4\phantom{.}$$
 make sense and respectively become $k_*\bar N_{\hat G_1} (\R_1)\-$ and 
 $k_*\skew6\widehat{\bar N}_{G_1}^{^{E_1}} (\R_1)\-$modules.

\smallskip
At this point, it follows from the induction hypothesis that we have {\it natural\/} 
$k_*\bar N_{\hat G_1} (R^1)\-$ and $k_*\bar N_{\hat G_1} (R^1_{\,E_1})\-$module isomorphisms
$$U_1\cong {\rm Ind}_{\bar N_{\hat G_1}(\R_1)}^{\bar N_{\hat G_1}(R^1)}({}^\omega W^1)\qq  
V_1\cong {\rm Ind}_{\skew6\widehat{\bar N}_{G_1}^{^{E_1}} (\R_1)}^{\skew3\hat{\bar N}_{G_1}(R^1_{\,E_1})} (W^1)
\eqno £7.19.5.$$
But, it follows from isomorphisms~£7.9.4 and~£7.10.2, and from Corollary~£7.12 and~Proposition~£7.14
 that, considering the surjective $k^*\-$group homomorphism (cf.~£2.8.1 and~Propositions~£7.11 and~£7.14)
$$\eqalign{\skew3\hat{\bar N}_0 = \bar N_{\hat G}(R_\delta)\star \skew3\hat{\bar N}_{G} 
(R_\delta)^\circ &\too \hat F_{k_*\hat G}(R_\delta)\buildrel  {}^\omega\phi_0\over 
\cong \bar N_{\hat G_1}(R^1)\cr
\skew6\widehat{\bar N}^{^{E}}_0 = \skew3\hat {\bar N}_{G}(R_{\,E})_{\delta} \star                            
\big(\skew3\hat {\bar N}_{G}  (R_{\,\delta})_{E}\big)^\circ &\too \skew3\hat {\bar N}_{G_{1}}   (R^{1}_{\; E_{1}})\cr}
\eqno £7.19.6\phantom{.}$$
and denoting by $\check U_1$ and $\check V_1$ the corresponding restricitions of 
$U_1$ to $\skew3\hat{\bar N}_0$ and of $V_1$ to
 $\skew6\widehat{\bar N}^{^{E}}_0\,,$ we have
$$U\cong {\rm Ind}_{\bar N_{\hat G} (R_{\,\delta})}^{ N_{\hat G} (R)} (W_0\otimes_k \check U_1)
\qq V\cong {\rm Ind}_{\skew3\hat{\bar N}_{G}  (R_{E})_{\delta}}^{\skew3\hat{\bar N}_{G} (R_{E})}  (W_0\otimes_k \check V_1)
\eqno £7.19.7.$$

\smallskip
Furthermore, it is easily checked that the  image of $\bar N_G (\R)\i \bar N_G(R_\delta)$
in~$\bar N_{\hat G_1}(R^1)$ throughout the $k^*\-$quotient of homomorphism~£7.19.2
is contained in $\bar N_{G_1}(\R_1)$ and then that this $k^*\-$group homomorphism induces
$k^*\-$group homomorphisms (cf.~£7.17.1 and~£7.18.2)
$$\eqalign{\bar N_{\hat G}(\R) \star \skew6\widehat{\bar N}^{^\delta}_{G} (\R)^\circ
\cong \skew6\widehat{\bar N}_G^{^1} (\R) &\too \bar N_{\hat G_1}(\R_1)\cr
 \skew6\widehat{\bar N}^{^{E}}_{G}(\R) \star \skew6\widehat{\bar N}^{^\delta}_{G} (\R)^\circ \cong 
 \skew6\widehat{\bar N}_G^{^{E_1}} (\R) &\too \skew6\widehat{\bar N}_{G_1}^{^{E_1}}(\R_1)\cr} 
\eqno £7.19.8;$$
thus, denoting by ${}^\omega\check W^1$ and by $\check W^1$ the corresponding restrictions  of ${}^\omega W^1$ 
to~$\skew6\widehat{\bar N}_G^{^1} (\R)$ and of $W^1$ to $\skew6\widehat{\bar N}_G^{^{E_1}} (\R)$ we have
(cf.~£7.19.5)
$$\check U_1 = {\rm Ind}_{\skew6\widehat{\bar N}_G^{^1} (\R) }^{\skew3\hat{\bar N}_0}({}^\omega \check W^1)\qq 
\check V_1 = {\rm Ind}_{\skew6\widehat{\bar N}_G^{^{E_1}} (\R) }
^{\skew6\widehat{\bar N}^{^E}_0}(\check W^1)
\eqno £7.19.9.$$

\smallskip
More explicitly, for any $n\ge 1$ the following diagrams of $k^*\-$group homomorphisms
$$\matrix{\skew6\widehat{\bar N}_{G_1}^{^n} (\R_1) \star \skew6\widehat{\bar N}_{G_1}^{^{n+1}} (\R_1)^\circ &\buildrel ({}^\omega\Psi_{1,n})^{-1}\over\cong & \skew6\widehat{\bar N}^{^{\delta_n}} _{\!G_1}(\R_1) &
\buildrel \skew3\hat{\bar\mu}^{\delta_n}_{1,n} \over\too &
\skew3\hat{\bar N}_{G_n} (R^n_{\,\delta_n})\cr
\big\uparrow&\phantom{\Big\uparrow}&\big\uparrow&&\Vert\cr
\skew6\widehat{\bar N}_{G}^{^n} (\R) \star \skew6\widehat{\bar N}_{G}^{^{n+1}} (\R)^\circ 
&\buildrel  ({}^\omega\Psi_{n})^{-1}\over\cong & \skew6\widehat{\bar N}^{^{\delta_n}} _{\!G}(\R) 
&\buildrel \skew3\hat{\bar\mu}^{\delta_n}_{n} \over\too 
&\skew3\hat{\bar N}_{G_n} (R^n_{\,\delta_n})\cr}
\eqno £7.19.10\phantom{.}$$
$$\matrix{\skew6\widehat{\bar N}_{G_1}^{^{E_n}} (\R_1) \star \skew6\widehat{\bar N}_{G_1}^{^{E_{n+1}}} (\R_1)^\circ &\buildrel (\Psi_{1,n})^{-1}\over\cong & \skew6\widehat{\bar N}^{^{\delta_n}} _{\!G_1}(\R_1) &
\buildrel \skew3\hat{\bar\mu}^{\delta_n}_{1,n} \over\too &
\skew3\hat{\bar N}_{G_n} (R^n_{\,\delta_n})\cr
\big\uparrow&\phantom{\Big\uparrow}&\big\uparrow&&\Vert\cr
\skew6\widehat{\bar N}_{G}^{^{E_n}} (\R) \star \skew6\widehat{\bar N}_{G}^{^{E_{n+1}}} (\R)^\circ 
&\buildrel  (\Psi_{n})^{-1}\over\cong & \skew6\widehat{\bar N}^{^{\delta_n}} _{\!G}(\R) 
&\buildrel \skew3\hat{\bar\mu}^{\delta_n}_{n} \over\too 
&\skew3\hat{\bar N}_{G_n} (R^n_{\,\delta_n})\cr}
\eqno £7.19.11\phantom{.}$$
are commutative since all the vertical arrows are defined by {\it pull-back via\/}  the group 
homomorphism $N_G(\R)\to N_{G_1} (\R_1)$ determined by the $k^*\-$quotient of homomorphism~£7.19.2; 
hence, we actually get a  $k_*\skew6\widehat{\bar N}_{G}^{^1} (\R)\-$ and a 
$k_*\skew6\widehat{\bar N}_{G}^{^{E_1}} (\R)\-$ module isomorphisms
$$\eqalign{{}^\omega \check W^1&\cong \bigotimes_{n\ge 1}  {\rm Res}_{\skew3\hat{\bar\mu}^{\delta_n}_n
\circ ({}^\omega\Psi_n)^{-1}}(W_n)\cr 
\check W^1&\cong \bigotimes_{n\ge 1}  {\rm Res}_{\skew3\hat{\bar\mu}^{\delta_n}_n
\circ (\Psi_n)^{-1}}(W_n)\cr }
\eqno £7.19.12\phantom{.}$$

\smallskip
Consequently, from the {\it Frobenius property\/}, we get a 
$k_*\bar N_{\hat G}  (R_{\,\delta})\-$module isomorphism
$$\eqalign{W_0\otimes_k \check U_1&\cong 
{\rm Ind}_{\bar N_{\hat G}(\R)}^{\bar N_{\hat G}  (R_{\,\delta})}
\big({\rm Res}_{\,\skew3\hat{\bar\mu}^{\delta}_0\circ ({}^\omega\Psi_0)^{-1}}(W_0)
\otimes_k {}^\omega\check W^1\big)\cr
&\cong {\rm Ind}_{\bar N_{\hat G}(\R)}^{\bar N_{\hat G} 
(R_{\,\delta})}(W)\cr}
\eqno £7.19.13\phantom{.}$$
and therefore from the left-hand isomorphism in~£7.19.7 we obtain the left-hand isomorphism in~£7.19.1. Similarly, we get a $k_*\bar N_{\hat G}  (R_{\,\delta})_E\-$module isomorphism
$$\eqalign{W_0\otimes_k \check V_1&\cong  {\rm Ind}_{\skew6\widehat{\bar N}^{^E}_{G}(\R)}^{\skew3\hat{\bar N}_{G}  (R_{E})_\delta}
\big({\rm Res}_{\,\skew3\hat{\bar\mu}^{\delta}_0\circ (\Psi_0)^{-1}}(W_0)
\otimes_k \check W^1\big)\cr
&\cong {\rm Ind}_{\skew6\widehat{\bar N}^{^E}_{G}(\R)}^{\bar N_{\hat G} 
(R_{\,\delta})_E}(W)\cr}
\eqno £7.19.14\phantom{.}$$
and therefore from the right-hand isomorphism in~£7.19.7 we obtain the right-hand isomorphism in~£7.19.1. We are done.

\medskip
£7.20. In order to compare both isomorphisms in~£7.19.1, note that from homomorphism~£2.8.2 and from our choice of a 
{\it polarization\/} $\omega$ we have a $k^*\-$group homomorphism
$$\skew3\hat{\bar N}_G (R_E)^\circ \star \bar N_{\hat G} (R_E)\too \hat F_T (R)\buildrel \omega_{(R,T)}\over
{\hbox to 25pt{\rightarrowfill}} k^*
\eqno £7.20.1\phantom{.}$$
which determines a $k^*\-$group isomorphism $\bar N_{\hat G} (R_E)\cong \skew3\hat{\bar N}_G (R_E)\,;$
let us denote by ${}^\omega V$ the restriction of $V$ throughout  this isomorphism. Similarly, for any $n\in \Bbb N\,,$
the $k_*\big(\skew6\widehat{\bar N}_G^{^{E_n}} (\R) \star \skew6\widehat{\bar N}_G^{^{E_{n+1}}} (\R)^\circ 
\big)\-$module ${\rm Res}_{\skew3\hat{\bar\mu}^{\delta_n}_n
\circ (\Psi_n)^{-1}}(W_n)$ restricted throughout the composed  $k^*\-$group isomorphism
$$\skew6\widehat{\bar N}_G^{^{n}}(\R) \star \skew6\widehat{\bar N}_G^{^{{n+1}}} (\R)^\circ  
\buildrel ({}^\omega\Psi_n)^{-1}\over\cong\skew6\widehat{\bar N}^{^{\delta_n}} _{\!G}(\R) \buildrel \Psi_n\over
\cong \skew6\widehat{\bar N}_G^{^{E_n}} (\R) \star \skew6\widehat{\bar N}_G^{^{E_{n+1}}} (\R)^\circ 
\eqno £7.20.2\phantom{.}$$
coincides with ${}^\omega W_n\,.$ 

\medskip
£7.21. But, according to the right-hand $k^*\-$group isomorphism in~£4.5.4, the corresponding splitting
$$\big(\skew6\widehat{\bar N}_G^{^{n}}(\R) \star
\skew6\widehat{\bar N}_G^{^{{n+1}}} (\R)^\circ\big) \star \big(\skew6\widehat{\bar N}_G^{^{E_n}} (\R) \star \skew6\widehat{\bar N}_G^{^{E_{n+1}}} (\R)^\circ \big)^\circ\too k^*
\eqno £7.21.1\phantom{.}$$
comes from $\omega_{(R^n,{\rm Res}_{R^n}^{P^n}(S_n))}\,\colon \hat F_{S_n} (R^n)\to k^*$
and needs not coincide with the splitting
$$\big(\skew6\widehat{\bar N}_G^{^{n}}(\R) \star \skew6\widehat{\bar N}_G^{^{E_n}} (\R)^\circ\big) \star
\big(\skew6\widehat{\bar N}_G^{^{{n+1}}} (\R)  \star \skew6\widehat{\bar N}_G^{^{E_{n+1}}} (\R)^\circ \big)^\circ
\too  k^*
\eqno £7.21.2\phantom{.}$$
coming from (cf.~£2.8.1)
$$\omega_{(R^n,T_n)}\,\colon\hat F_{T_n}(R^n)\to k^*\qq \omega_{(R^{n+1},T_{n+1})}\,\colon \hat F_{T_{n+1}}(R^{n+1})\to k^*
\eqno £7.21.3;$$
that is to say, this splitting determines a new  $k^*\-$group isomorphism between $\skew6\widehat{\bar N}_G^{^{n}}(\R) \star \skew6\widehat{\bar N}_G^{^{{n+1}}} (\R)^\circ$ and $\skew6\widehat{\bar N}_G^{^{E_n}} (\R) \star \skew6\widehat{\bar N}_G^{^{E_{n+1}}} (\R)^\circ \,;$ thus, this isomorphism and $\psi_n\circ ({}^\omega \psi_n)^{-1}$ determine 
an automorphism ${}^\omega \theta_n$ of $\skew6\widehat{\bar N}_G^{^{n}}(\R) \star \skew6\widehat{\bar N}_G^{^{{n+1}}} 
(\R)^\circ\,.$ Then, it is clear that the product of all these automorphisms defines an automorphism~${}^\omega\theta$
of~$\skew6\widehat{\bar N}_G^{^{0}}(\R)  =\bar N_{\hat G} (\R)$ and that the right-hand isomorphism in~£7.19.1 implies the following result.

\bigskip
\noindent
{\bf Corollary~£7.22.} {\it With the notation and the choice above, we have a natural 
$k_*\bar N_{\hat G} (R_E)\-$module isomorphism
$${}^\omega V\cong {\rm Ind}_{\bar N_{\hat G}(\R)}^{\bar N_{\hat G} (R_E)}\big({\rm Res}_{\,{}^\omega\theta} ({}^\omega W)\big)
\eqno £7.22.1.$$ \/}

\bigskip
\noindent
{\bf Appendix: The odd order case}
\bigskip
£A.1. Assume that $p \not= 2$ and let $\hat G$ be a $k^*\-$group with finite $k^*\-$quotient $G$ of {\it odd\/} order.
In this case, by the fundamental Feit-Thompson Theorem [3], $G$ is {\it solvable\/}
and therefore, for any choice of a {\it polarization\/} $\omega\,,$ Theorem~£6.5 above supplies a {\it natural\/} bijection
$${\rm Irr}_k (\hat G)\cong {\rm Wgt}_k (\hat G)
\eqno £A.1.1;$$
actually, it suffices to consider $\omega$ over the {\it torsion\/} subcategory
 $\frak D_k^{\rm tor}$ (cf.~Remark~£4.7); further, the oddness of our groups only demands the choice of a splitting for the 
 $k^*\-$subgroup $\Bbb O^2\big(\hat F_S (P)\big)$ of  $\hat F_S (P)$ for any $\frak D_k^{\rm tor}\-$object $(P,S)\,.$

 \medskip
 £A.2.  That is to say, in the present situation we can replace $\frak D_k\,,$ $\hat\frak f$ and  $\omega$ (cf.~£2.5) by the 
 {\it full\/} subcategory $\frak D_k^{\rm tor}$ of $\frak D_k\,,$ by the subfunctor of $\hat\frak f$
 $${}^2\hat\frak f : \frak D_k^{\rm tor}\too k^*\-\Gr
 \eqno £A.1.2\phantom{.}$$
 mapping  any $\frak D_k^{\rm tor}\-$object $(P,S)$ on  the $k^*\-$group 
 $\Bbb O^2\big(\hat F_S (P)\big)\,,$ and finally by a {\it natural\/} map 
 ${}^2\omega\,\colon {}^2\hat\frak f\too k^*$ fulfilling the condition in~£2.15.1  --- called a {\it odd-polarization\/}.  Although any  {\it odd-polarization\/} can be easily extended to a {\it polarization\/}, the point is that  there is a unique  
 {\it odd-polarization\/} com-patible with the {\it tensor product\/}  of Dade $P\-$algebras. We borrow the notation from [15,~Chap.~9] and denote by   $\D_k^{\rm tor}(P)$ the subgroup of {\it torsion\/} elements of~$\D_k (P)\,;$ actually, it is known that all the nontrivial {\it torsion\/} elements of $\D_k (P)$ have order 2 [15,~8.16 and~Corollary~8.22] or, equivalently, that $S\cong S^\circ$  for any $\frak D_k^{\rm tor}\-$object $(P,S)\,.$
\eject

 \bigskip
 \noindent
 {\bf Theorem~£A.3.} {\it There is a unique odd-polarization ${}^2 \omega$ such that  the
following diagram is commutative
$$\matrix{\Bbb O^2\big(\hat F_S (P)\big)\,\,\hat\cap\,\,\Bbb O^2\big(\hat F_{S'} (P)\big)
&\buildrel \hat\nu_{P,S,S'}\over{\hbox to 30pt{\rightarrowfill}} &\Bbb O^2\big(\hat F_{S\otimes_k S'} (P)\big)\cr
{\scriptstyle {{}^2\omega_{(P,S)}\,\hat\times\,{}^2\omega_{(P,S')}}}\; \searrow\hskip-15pt
&\phantom{\bigg\uparrow} & \swarrow\hskip5pt
\scriptstyle {{}^2\omega_{(P,S\otimes_k S')}}\cr &\ k^* &\cr}
\eqno £A.3.1.$$ 
for any pair of $\frak D_k^{\rm tor}\-$objects $(P,S)$ and $(P,S')\,.$ 
Moreover, for any normal subgroup $Q$ of $P\,,$ setting  $T ={\rm Res}_Q^P (S)$ and $\bar P = P/Q\,,$ 
and denoting by $\hat F_S (P)_Q$ the stabilizer of $Q$ in $\hat F_S (P)\,,$
the following diagram is also commutative
$$\matrix{\Bbb O^2\big(\hat F_S (P)_Q\big)Ê\hskip-5pt&\buildrel \Delta_{P,S,Q}
\over{\hbox to 30pt{\rightarrowfill}} &\hskip-5pt \Bbb O^2\big(\hat F_T (Q)\big)
\,\hat\times\, \Bbb O^2\big(\hat F_{S(Q)}(\bar P)\big)\cr 
{\scriptstyle {}^2\omega_{(P,S)}} \searrow \hskip-20pt&\phantom{\bigg\uparrow} &
\hskip-20pt\swarrow\hskip5pt \scriptstyle
{{}^2\omega_{(Q,T)}\,\hat\times\,{}^2\omega_{(\bar P,S(Q))}}\cr  &\ k^* &\cr}
\eqno £A.3.2\phantom{.}$$\/}

\medskip
\noindent
{\bf Proof:} Let $\omega$ be a {\it polarization\/} [15,~Theorem~9.21];  for any $\frak D_k^{\rm tor}\-$object $(P,S)\,,$
it is clear that there is a group homomorphism $\beta_{(P,S)}\,\colon \Bbb O^2\big( F_S (P)\big)\to k^*$ fulfilling
$$(\omega_{(P,S)}\,\hat\times\,\omega_{(P,S)})(\hat\varphi\. \hat\varphi) = \beta_{(P,S)} 
(\varphi)\omega_{(P,S\otimes_k S)}\big(\hat\nu_{P,S,S'}(\hat\varphi\.\hat\varphi)\big)
\eqno £A.3.3\phantom{.}$$
for any $\hat\varphi\in \Bbb O^2\big(\hat F_S (P)\big)\,,$ where $\hat\varphi\. \hat\varphi$ denotes the image
of $(\hat\varphi,\hat\varphi)$ in the $k^*\-$group (cf.~£2.2)
$$\hat F_S (P)\,\,\hat\cap\,\,\hat F_{S} (P) = \hat F_S (P) \star \hat F_{S} (P)
\eqno £A.3.4\phantom{.}$$
 and $\varphi$ is the image of $\hat\varphi$ in 
$\Bbb O^2\big( F_S (P)\big)\,;$ then, there is a {\it unique\/} group homomorphism $\alpha_{(P,S)}\,\colon 
\Bbb O^2\big( F_S (P)\big)\to k^*$ fulfilling $(\alpha_{(P,S)})^2 = \beta_{(P,S)}$ and we claim that it suffices to define
$${}^2 \omega_{(P,S)}(\hat\varphi) = \alpha_{(P,S)}(\varphi)^{-1}\omega_{(P,S)}(\hat\varphi)
\eqno £A.3.5\phantom{.}$$
for any $\hat\varphi\in \Bbb O^2\big(\hat F_S (P)\big)\,.$

\smallskip
In any case, note that the uniqueness of ${}^2 \omega_{(P,S)}$ follows from the uniqueness of $\alpha_{(P,S)}\,.$
The commutativity of diagram~£A.3.1 for $S'= S$ follows from
 our very definition; otherwise, the diagrams corresponding to the pairs of Dade $P\-$algebras 
 $(S\otimes_k S',S\otimes_k S')\,,$ $(S,S)$ and $(S',S')$ are certainly commutative and then the commutativity of
 diagram~£A.3.1 follows.

 \smallskip
Moreover, once again it is clear that there is a group homomorphism 
$$\gamma_{(P,S,Q)} : \Bbb O^2\big( F_S (P)_Q\big)\too k^*
\eqno £A.3.6\phantom{.}$$
such that, for any $\hat\varphi\in \Bbb O^2\big(\hat F_S (P)_Q\big)\,,$ we have
$${}^2\omega_{(P,S)} (\hat\varphi) = \gamma_{(P,S,Q)} (\varphi) ({}^2\omega_{(Q,T)}\,\hat\times\,{}^2\omega_{(\bar P,S(Q))})\big(\Delta_{P,S,Q} (\hat\varphi)\big)
\eqno £A.3.7.$$
\eject
\noindent
But, it follows from [15, Proposition~9.16] that the diagram
$$\matrix{\hat F_S (P)_Q \starÊ\hat F_{S} (P)_Q\hskip-8pt
&\too &\hskip-8pt\big(\hat F_T (Q) \star \hat F_{T} (Q)\big)\hat\times
\big(\hat F_{S(Q)} (\bar P) \star \hat F_{S(Q)} (\bar P)\big)\cr
\hskip-20pt{\scriptstyle \hat\nu_{P,S,S}}\downarrow
&\phantom{\bigg\downarrow}&\downarrow {\scriptstyle \hat\nu_{Q,T,T}
\,\hat\times\, \hat\nu_{\bar P,S(Q),S(Q)}}\hskip-50pt\cr
\hat F_{S\otimes_k S} (P)_QÊ\hskip-8pt
& \too  &\hskip-15pt\hat F_{T\otimes_k T} (Q)\,\hat\times\, \hat F_{(S\otimes_k S)(Q)} (\bar P)\cr} 
\eqno £A.3.8\phantom{.}$$
is commutative; moreover, since the Dade $P\-$algebra 
$$S\otimes_k S\cong S\otimes_k S^\circ \cong {\rm End}_k(S)
\eqno £A.3.9\phantom{.}$$
 is {\it similar\/} to  $k\,,$ the corresponding diagram~£A.3.2 is clearly commutative. Consequently, for any $\hat\varphi\in 
\Bbb O^2\big(\hat F_S (P)_Q\big)\,,$ the element $({}^2\omega_{(P,S)}\,\hat\times\,{}^2\omega_{(P,S)}) 
(\hat\varphi\.\hat\varphi)$ coincides with the image of $\Delta_{P,S,Q} (\hat\varphi)\.\Delta_{P,S,Q} (\hat\varphi)$
throughout the map
$$ ({}^2\omega_{(Q,T)}\,\hat\times\,{}^2\omega_{(Q,T)})\,\hat\times\,({}^2\omega_{(\bar P,S(Q))}\,\hat\times\,{}^2\omega_{(\bar P,S(Q))})
\eqno £A.3.10\phantom{.}$$
and therefore we get $\gamma_{(P,S,Q)} (\varphi)^2 = 1$ which forces $\gamma_{(P,S,Q)} (\varphi) = 1\,.$
We are done.

\medskip
£A.4. Since the unique {\it odd-polarization\/} ${}^2\omega$ in Theorem~£4.3 can be easily extended to a {\it polarization\/},
if follows from Theorem~£6.5 above that it supplies a {\it natural\/} bijection
$${\rm Irr}_k (\hat G)\cong {\rm Wgt}_k (\hat G)
\eqno £A.4.1\phantom{.}$$
and we claim that this bijection coincides with the bijection defined by Gabriel Navarro in [7,~Theorem~4.3]
for $\pi = \{p\}\,.$ First of all, borrowing all the notation in \S7, suitably translated to our present situation, and choosing 
this {\it odd-polarization\/} ${}^2\omega\,,$ we claim that the corresponding automorphism~${}^{{}^2\omega}\theta$ of $\bar N_G (\R)$ in~£7.21 above is the identity map and therefore, according to Theorem~£7.19 
and Corollary~£7.22,  in tis case we have
$$U\cong {\rm Ind}_{\bar N_{\hat G} (R_E)}^{\bar N_{\hat G}(R)}(\,{}^{{}^2\omega} V)
\eqno £A.4.2.$$
That is to say, in the bijection~£A.4.1 determined by ${}^2\omega$ the image of any simple $k_*\hat G\-$module $M$
 can be directly computed from  the triple formed by  a {\it vertex\/} $R\,,$ an {\it $R\-$source\/} $E$ and a 
 {\it multiplicity $\skew3\hat{\bar N}_G (R_E)\-$module\/}
$V$ of $M\,.$

\medskip
£A.5. More precisely, for any $n\in \Bbb N\,,$ we claim that the automorphism ${}^{{}^2\omega}\theta_n$ 
of~$\skew6\widehat{\bar N}_G^{^{n}}(\R) \star \skew6\widehat{\bar N}_G^{^{{n+1}}}  (\R)^\circ$
is the identity map; indeed, since  (cf.~£7.13) 
$${\rm Res}_{\rho_n}(T_n) \cong {\rm Res}_{\rho_n} (S_n)\otimes_k {\rm Res}_{\rho_{n+1}} (T_{n+1})
\eqno £A.5.1,$$
up to a suitable identification, from Theorem~£A.3 above we get the following commutative diagram
$$\matrix{\Bbb O^2\big(\hat F_{S_n} (R^n)\big)\,\,\hat\cap\,\,\Bbb O^2\big(\hat F_{T_{n+1}} (R^{n})\big)
&\too &\Bbb O^2\big(\hat F_{T_n} (R^n)\big)\cr
{\scriptstyle {{}^2\omega_{(R^n,S_n)}\,\hat\times\,{}^2\omega_{(R^{n},T_{n+1})}}}\; \searrow\hskip-10pt
&\phantom{\bigg\uparrow} & \swarrow\hskip5pt
\scriptstyle {{}^2\omega_{(R^n,T_n)}}\cr &\ k^* &\cr}
\eqno £A.5.2\phantom{.}$$ 
which, according to the very definition of ${}^{{}^2\omega}\theta_n\,,$ proves our claim.
\eject

\medskip
£A.6. In particular, if $\hat G'$ is a $k^*\-$subgroup of $\hat G$ and $M'$ a $k_*\hat G'\-$module such that
$M\cong {\rm Ind}_{\hat G'}^{\hat G}(M')\,,$ then $M'$ is clearly a simple $k_*\hat G'\-$module and it is 
easily checked that a vertex $R'$ and an $R'\-$source $E'$ of $M'$  are also a vertex  and an $R'\-$source of $M\,;$
moreover, we claim that {\it if $V'$ is a multiplicity $\skew3\hat{\bar N}_{G'} (R'_{E'})\-$module of $M'$  then the 
$k_*\skew3\hat{\bar N}_{G} (R'_{E'})\-$module
$$V = {\rm Ind}_{\skew3\hat{\bar N}_{G'} (R'_{E'})}^{\skew3\hat{\bar N}_{G} (R'_{E'})}(V')
\eqno £A.6.1$$
is a  multiplicity $\skew3\hat{\bar N}_{G} (R'_{E'})\-$module of $M\,.$\/} Indeed, recall that we have (cf.~£2.12)
$${\rm Ind}_{\hat G'}^{\hat G}\big({\rm End}_k (M')\big)\cong {\rm End}_k (M)
\eqno £A.6.2\phantom{.}$$
and that ${\rm id}_M$ is the image of ${\rm Tr}_{\hat G'}^{\hat G}(1\otimes {\rm id}_{M'}\otimes 1)\,,$
and denote by $G'_{M'}$ the pointed group on ${\rm End}_k (M)$ determined by the group $G'$ and the idempotent
$1\otimes {\rm id}_{M'}\otimes 1\,.$  Since $R'_{E'}$ is a local pointed group on ${\rm End}_k (M)\,,$ the unity element 
in  $\big({\rm End}_k (M)\big)(R'_{E'})$ coincides with the sum $\sum_x \overline{x\otimes {\rm id}_{M'}\otimes x^{-1}}$ 
where $x$ runs over the elements fulfilling $(R'_{E'})^x\i G'_{M'}$ in a set of representatives for $\hat G/\hat G'$ in~$\hat G$ and, for such an element $x\,,$ $\overline{x\otimes {\rm id}_{M'}\otimes x^{-1}}$ denotes the image 
of~$x\otimes {\rm id}_{M'}\otimes x^{-1}$in  $\big({\rm End}_k (M)\big)(R'_{E'})$ [9,~Proposition~1.3].
But, it is clear that $(R'_{E'})^x$ is also a maximal local pointed group on ${\rm End}_k (M')$ and therefore there is 
$x'\in\hat G'$ such that $(R'_{E'})^x = (R'_{E'})^{x'}\,.$ Consequently, it follows from [6,~statement~2.13.2] that we get an $\skew3\hat{\bar N}_{G} (R'_{E'})\-$interior algebra isomorphism
$${\rm Ind}_{\skew3\hat{\bar N}_{G'} (R'_{E'})}^{\skew3\hat{\bar N}_{G} (R'_{E'})}
\Big(\big({\rm End}_k (M')\big)(R'_{E'})\Big)\cong \big({\rm End}_k (M)\big)(R'_{E'})
\eqno £A.6.3\phantom{.}$$
which proves our claim.

\medskip
£A.7.  Moreover, by the very definitions of ${}^{{}^2\omega} V$ and ${}^{{}^2\omega} V'$ in~£7.20 above,
then we still have 
$${}^{{}^2\omega} V = {\rm Ind}_{\bar N_{\hat G'} (R'_{E'})}^{\bar N_{\hat G} (R'_{E'})}({}^{{}^2\omega} V')
\eqno £A.7.1\phantom{.}$$
and therefore it follows from isomorphism~£A.4.2 that we have a $k_*\bar N_{\hat G}(R')\-$mo-dule isomorphism
$$U\cong {\rm Ind}_{\bar N_{\hat G'}(R')}^{\bar N_{\hat G}(R')}(U')
\eqno £A.7.2\phantom{.}$$
where $U'$ is a {\it simple projective\/} $\bar N_{\hat G'}(R')\-$module which, together with $R'$, determines
the $G'\-$conjugacy class of {\it weights\/} of $\hat G'$ determined by $M'$ {\it via\/} the corresponding 
bijection £A.4.1.

\medskip
£A.8. In conclusion, in order to prove that Navarro's correspondence in [7,~Theorem~4.3] also maps
$M$ on the $G\-$conjugacy class of the {\it weight\/} of~$\hat G$ determined by $R$ and $U\,,$
we may assume that $M$ is {\it primitive\/} --- namely,\break
\eject
\noindent
that it is {\it not\/} induced from any proper $k^*\-$subgroup of $\hat G\,.$ In this case,
as we mention in~£1.4 above, it follows from  [17,~Lemma~30.4] that there is a $G\-$stable finite $p'\-$subgroup $K$ 
of ${\rm End}_k (M)^*$ which generates  the $k\-$algebra~${\rm End}_k (M)\,;$ in particular, ${\rm End}_k(M)$ is actually a {\it Dade $R\-$algebra\/} [13,~1.3], $R$ is a Sylow $p\-$subgroup of $G$ and $\bar N_G (R)\-$stabilizes the isomorphism class of $E$ [13,~1.8], so that $U\cong {}^{{}^2\omega} V$ (cf.~£A.4.2). At this point, a careful inspection of the origin of Navarro's correspondence in [7, Theorem~3.1] shows that it maps $M$ on the $G\-$conjugacy class of the {\it weight\/} of~$\hat G$ determined 
by $R$ and $U$ if, for a suitable Brauer character $\psi$ over $N_{\hat G}(R)\,,$ we have
$${\rm Res}_{N_{\hat G}(R)}^{\hat G}(\varphi_M) = \varphi_U + 2\.\psi
\eqno £A.8.1\phantom{.}$$
where $\varphi_M$ and $\varphi_U$ respectively denote  the Brauer characters of $M$ and of the $k_*N_{\hat G} (R)\-$module
$U\cong {}^{{}^2\omega} V\,.$ Then, the fact that in our situation such an equality holds is more or less a consequence of [5,~Theorem~5.3] but here we give a direct proof.

\bigskip
\noindent
{\bf Proposition~£A.9.} {\it Let $M$ be a simple primitive $k_*\hat G\-$module, $R$ a vertex, $E$ an $R\-$source 
and $V$ a multiplicity $k_*\skew3\hat{\bar N}_G (R_E)\-$module of $M\,.$
Consider the unique odd-polarization ${}^2\omega$ such that diagram~{\rm £4.3.1} is commutative and denote by
${}^{{}^2\omega} V$ the restriction of $\,V$ throughout the isomorphism $\bar N_{\hat G}(R)\cong 
\skew3\hat{\bar N}_G (R_E)$ determined by ${}^2\omega\,,$ and by $\varphi_M$ and $\varphi_{\,{}^{{}^2\omega} V}$ the respective Brauer characters of~$M$ and of ${}^{{}^2\omega} V$ considered as a $k_* N_{\hat G} (R)\-$module. Then,  for a suitable Brauer character $\psi$ over $N_{\hat G}(R)\,,$ we have
$${\rm Res}_{N_{\hat G}(R)}^{\hat G}(\varphi_M) = \varphi_{\,{}^{{}^2\omega} V} + 2\.\psi
\eqno £A.9.1.$$\/}

\par
\noindent
{\bf Proof:} Arguing by induction on $\vert G\vert\,,$ we may assume that $M$ is a faithful $k_*\hat G\-$module, then
identifying $\hat G$ with its image in ${\rm End}_k (M)\,;$ moreover, the case where ${\rm dim}_k (M) = 1$ being clear,
we assume that ${\rm dim}_k (M)\not= 1\,.$ Then, a minimal normal nontrivial subgroup $K$ of $G$
is an Abelian $\ell\-$elementary group for a prime number $\ell\not= p$ and the {\it primitivity\/} of $M$ 
forces the converse image $\hat K$ of~$K$ in $\hat G$ to be the {\it central\/} product of $k^*$ by an 
{\it extra-special\/} normal $\ell\-$subgroup of~$\hat G$ [4,~Ch.~5,~\S5].

\smallskip
Let $S$ be the $k\-$subalgebra of ${\rm End}_k (M)$ generated by $\hat K\,;$ once again, the {\it primitivity\/} of $M$ 
forces $S$ to be a simple $k\-$algebra and then the $k^*\-$quotient $G$ of $\hat G$ acts on $S$
determining a $k^*\-$group $^{\skew6\hat {}}G$ together with  a $k^*\-$group homomorphism 
$\,{}^{\skew6\hat {}} G\to S^*$ (cf.~£2.3), and we set
$$ G^{\skew4\hat {}}\, = \hat G \star (\,{}^{\skew6\hat {}} G)^\circ
\eqno £A.9.2;$$
\eject
\noindent
then, it follows from [16,~Proposition~3.2] that there exists a $k_*G^{\skew4\hat {}}\-$module $\bar M$
such that we have a $\hat G$-interior algebra isomorphism
$${\rm End}_k (M)\cong  S\otimes_k {\rm End}_k (\bar M)
\eqno £A.9.3;$$ 
actually, $\hat K$ is {\it canonically\/} isomorphic to the converse image of $K$ in ${}^{\skew6\hat {}} G$
and therefore $K$ lifts to a normal subgroup of $G^{\skew4\hat {}}\,$ acting trivially on $\bar M\,;$
thus, up to suitable identifications, setting $\skew3\hat{\bar G} = G^{\skew4\hat {}}\,/K$ and $S = {\rm End}_k (N)\,,$
$\bar M$ becomes a $k_*\skew3\hat{\bar G}\-$module,  we have a $k_*\hat G\-$module isomorphism
$$M \cong N\otimes_k \bar M
\eqno £A.9.4\phantom{.}$$
and, denoting by $\varphi_N$ the Brauer character of $N$ and by $\hat\pi\,\colon G^{\skew4\hat {}}\to \skew3\hat{\bar G}$
the canonical $k^*\-$group homomorphism, we have 
$$\varphi_M = \varphi_N\. {\rm Res}_{\hat\pi}(\varphi_{\bar M})
\eqno £A.9.5.$$

\smallskip
Now, it is clear that $\bar M$ is a  simple primitive $k_*\skew3\hat{\bar G}\-$module, that the image $\bar R$ of 
$R$ in $\bar G = G/K$ is a vertex of $\bar M$ (actually, it is a Sylow $p\-$subgroup of~$\bar G$), that we have a canonical $R\-$interior algebra {\it embedding\/} (cf.~£2.4)
$${\rm End}_k (E)\too S\otimes_k {\rm End}_k (\bar E)
\eqno £A.9.6\phantom{.}$$
where $\bar E$ denotes an $\bar R\-$source of $\bar M\,,$ and that we still have a $\skew3\hat{\bar N}_G (R_E)\-$interior algebra isomorphism [12,~Proposition~5.6]
$$\big({\rm End}_k (M)\big)(R_E)\cong  S (R)\otimes_k \big({\rm End}_k (\bar M)\big)(\bar R_{\bar E})
\eqno £A.9.7,$$
together with a $k^*\-$group isomorphism [12,~Proposition~5.11]
$$\skew3\hat{\bar N}_G (R_E)\cong \skew6\widehat{\bar N}^{^S}_G (R) \star {\rm Res}_{\pi}\big(\skew3\hat{\bar N}_{\bar G} (\bar R_{\bar E})\big)
\eqno £A.9.8\phantom{.}$$
where $\skew6\widehat{\bar N}^{^S}_G (R)$ and ${\rm Res}_{\pi}\big(\skew3\hat{\bar N}_{\bar G} (\bar R_{\bar E})\big)$ respectively denote the $k^*\-$groups coming from the action of $\bar N_G (R)$ 
on the simple $k\-$algebra $S(R)$ [13,~1.8], and  obtained by {\it pull-back\/} from the canonical
group homomorphism $\pi\,\colon \bar N_G(R)\to \bar N_{\bar G}(\bar R)\,.$

\smallskip
On the one hand, denoting by $\bar V$ a {\it multiplicity $\skew3\hat{\bar N}_{\bar G} (\bar R_{\bar E})\-$module\/}
of $\bar M\,,$ so that we have
$${\rm End}_k (\bar V) \cong \big({\rm End}_k (\bar M)\big)(\bar R_{\bar E})
\eqno £A.9.9,$$
it follows from the induction hypothesis that,  for a suitable Brauer character~$\bar\psi$ over $N_{\skew3\hat{\bar G}}(\bar R)\,,$ we have
$${\rm Res}_{N_{\skew3\hat{\bar G}}(\bar R)}^{\skew3\hat{\bar G}}(\varphi_{\bar M}) = 
\varphi_{\,{}^{{}^2\omega}\bar V} + 2\.\bar\psi
\eqno £A.9.10.$$
On the other hand, denoting by $W$ a {\it multiplicity $\skew6\widehat{\bar N}^{^S}_G (R)\-$module\/} of $N\,,$
it follows from isomorphisms~£A.9.7 and~£A.9.8 that we have a $k_*\skew3\hat{\bar N}_G (R_E)\-$module isomorphism
$$V\cong W\otimes_k {\rm Res}_{\pi} (\bar V)
\eqno £A.9.11;$$
moreover, it follows from the commutativity of diagram~£A.3.1 that we still have a $k_*\bar N_{\hat G} (R)\-$module isomorphism
$${}^{{}^2\omega} V\cong {}^{{}^2\omega} W\otimes_k {\rm Res}_{\hat\pi_R} ({}^{{}^2\omega}\bar V)
\eqno £A.9.12\phantom{.}$$
where $\hat\pi_R$ denotes the restriction to $\bar N_{G^{\skew4\hat {}}}\, (R_E)$ of $\hat\pi\,;$ consequently, with evident notation, we get
$$\varphi_{\,{}^{{}^2\omega} V} = \varphi_{\,{}^{{}^2\omega} W}\.{\rm Res}_{\hat\pi_R} (\varphi_{\,{}^{{}^2\omega}\bar V})
\eqno £A.9.13.$$

\smallskip
But, according to Theorem~£A.10 below, we also have
$${\rm Res}_{N_{{}^{\skew6\hat {}} G}(R)}^{{}^{\skew6\hat {}} G}(\varphi_N) = 
\varphi_{\,{}^{{}^2\omega} W}  + 2\.\eta
 \eqno £A.9.14\phantom{.}$$
 for a suitable Brauer character $\eta$ over  $N_{{}^{\skew6\hat {}} G}(R)\,.$ In conclusion, from equalities~£A.9.5, £A.9.10
 and £A.9.14 we get
 $$\eqalign{{\rm Res}_{N_{\hat G}(R)}^{\hat G}(\varphi_M) 
 &= {\rm Res}_{N_{{}^{\skew6\hat {}} G}(R)}^{{}^{\skew6\hat {}} G}(\varphi_N)\. {\rm Res}_{\hat\pi_R}
 \big( {\rm Res}_{N_{\skew3\hat{\bar G}}(\bar R)}^{\skew3\hat{\bar G}} (\varphi_{\bar M})\big)\cr
&= (\varphi_{\,{}^{{}^2\omega} W}  + 2\.\eta)\.{\rm Res}_{\hat\pi_R}(\varphi_{\,{}^{{}^2\omega}\bar V} + 2\.\bar\psi)\cr
&= \varphi_{\,{}^{{}^2\omega} V} + 2\.\psi \cr}
\eqno £A.9.15\phantom{.}$$
where $\psi = \eta\.{\rm Res}_{\hat\pi_R}(\varphi_{\,{}^{{}^2\omega}\bar V} + 2\.\bar\psi) + 
\varphi_{\,{}^{{}^2\omega} W}\.{\rm Res}_{\hat\pi_R}(\bar\psi)\,.$ We are done.

\bigskip
\noindent
{\bf Theorem~£A.10.} {\it Let $M$ be a $k_*\hat G\-$module such that ${\rm End}_k (M)$ is genera-ted
by a $G\-$stable $k^*\-$subgroup $\hat K$ of $\,{\rm End}_k (M)^*$ which is the central pro-duct of $k^*$
 by an extra-special $\ell\-$subgroup for an odd prime number $\ell\not= p\,.$ For any local pointed group $R_E$
 on ${\rm End}_k (M)\,,$ denoting by $V$ a multiplicity $\skew3\hat{\bar N}_G (R_E)\-$module of $R_E$
 and by ${}^{{}^2\omega} V$ the restriction of $\,V$ via the isomorphism $\bar N_{\hat G}(R)\cong 
\skew3\hat{\bar N}_G (R_E)$ determined by the unique odd-polarization ${}^2\omega$ such that diagram~{\rm £4.3.1} is commutative, we have
$${\rm Res}_{N_{\hat G}(R)}^{\hat G}(\varphi_M) = \varphi_{\,{}^{{}^2\omega} V} + 2\.\psi
\eqno £A.10.1\phantom{.}$$
where  $\varphi_M$ and $\varphi_{\,{}^{{}^2\omega} V}$ denote the respective Brauer characters of $M$ and of 
$\,{}^{{}^2\omega} V$ considered as a $k_* N_{\hat G} (R)\-$module, and $\psi$ is a Brauer 
character over~$N_{\hat G}(R)\,.$\/}

\medskip
\noindent
{\bf Proof:} We actually may assume that $M$ is faithful and that $\hat G = \bar N_{\hat G}(R)\,;$ then, $G$ stabilizes the decomposition
[4,~Ch.~5, Theorem~2.3]
$$K = C_K (R)\times [K,R]
\eqno £A.10.2\phantom{.}$$
of the $k^*\-$quotient of $\hat K\,;$ thus, setting $S = {\rm End}_k (M)$ and denoting by $S'$ and $S''$ the $k\-$subalgebras
of $S$ generated by the respective converse images $\hat K'$ of~$C_K (R)$ and $\hat K''$ of $[K,R]\,,$ we have 
$S = S'\otimes_k S''\,,$  $\hat K'$ and $\hat K''$ are also  central products of $k^*$ by  extra-special $\ell\-$subgroups
(here  we also consider $\Bbb Z/\ell\Bbb Z$ as an {\it extra-special\/} $\ell\-$group) and
 $G$ still stabilizes 
 $$S' = {\rm End}_k (M')\qq S'' = {\rm End}_k (M'')
 \eqno £A.10.3.$$

 \smallskip
Consequently, as in the proof above, it follows  from the commutativity of diagram~£A.3.1 that it suffices to prove 
the theorem for $M'$ and for $M''\,.$ That is to say, we may assume that either $K = C_K (R)$ or $K = [K,R]\,;$ in the~first case, $R$ centralizes $\hat K$ [4,~Ch.~5, Theorem~1.4], so that it centralizes~$S$ which forces $R = \{1\}\,;$ then, we have $ M = V\,,$ $E = k$ and $\hat F_S (R) = k^*\,,$ and~by~the very definition of~$\skew3\hat{\bar N}_G (R_E)$
(cf.~£2.5) we get an isomorphism $\skew3\hat{\bar N}_G (R_E)\cong \bar N_{\hat G} (R) = \hat G$ compatible with the 
canonical $k^*\-$group homomorphism~£2.8.1, so that equality~£A.10.1 is trivially true with $\psi = 0\,.$

\smallskip
Following the notation in~£A.11 and according to isomorphism~£A.12.2 below, let us denote by $H$ the image of $Sp(K,\kappa)$ in $N_{S^*} (\hat K)\,;$
in particular, the nontrivial element in~$Z(H)$ is an involution $s\in S$ which stabilizes
$\hat K$ and induces $-{\rm id}_K$ over $K\,;$ note that, if $s'\in S$ is such an involution then $s's$ stabilizes $\hat K$
acting trivially on $K$ and therefore, according again to isomorphism~£A.12.2 below, $s'$ belongs to~$\hat K\.\langle s\rangle\,,$ so that we have $s'\in \{s^x,-s^x\}$ for a suitable $x\in \hat K\,.$

\smallskip
In the second case above, we have $C_K (R) = \{1\}$ and therefore $R$ fixes a unique pair of such involutions $\{s,-s\}$
which by oddness forces 
$$\hat G = N_{\hat G} (R) \i C_{S^*}(s)
\eqno £A.10.4;$$
consequently, $\hat G$ is contained in the intersection
$$N_{S^*} (\hat K)\cap C_{S^*}(s) =  k^*\times H
\eqno £A.10.5\,.$$
Moreover, since $K$ indexes an $R\-$stable basis of  $S = k_*\hat K$ (cf.~£A.11 below), we have $S(R)\cong k$ which
forces $V\cong k\,;$ hence, by the very definition of the $k^*\-$group $\skew3\hat{\bar N}_G (R_E)\,,$ in this case we get a $k^*\-$group isomorphism (cf.~£2.5)
$$\skew3\hat{\bar N}_G (R_E)\cong k^*\times \bar N_G (R)
\eqno £A.10.6\,.$$

\smallskip

At this point, it follows from Lemma~£A.14 below that the decomposition
$$\hat G = N_{\hat G}(R)\cong \hat{N}_G (R_E)\cong k^*\times G
\eqno £A.10.7\phantom{.}$$
determined by the $k^*\-$group homomorphism ${}^2 \omega_{(R,S)}\,\colon \hat F_S (R)\to k^*$ coincides\break
\eject
\noindent with  the decomposition
$$\hat G =  k^*\times (\hat G \cap H)
\eqno £A.10.8\phantom{.}$$
obtained from the inclusion $\hat G\i N_{S^*} (\hat K)\cap C_{S^*}(s)\,.$ In particular, the restriction of 
$\varphi_{\,{}^{{}^2\omega} V}$ to $\hat G \cap H$ is just the {\it trivial\/} character.

\smallskip
On the other hand, for any $y\in \hat G\cap H\,,$ acting over $K$ the product $sy$ only fix the trivial element $1\,;$ indeed, 
if $syx(sy)^{-1} = x$ for some $x\in K$ then $yx^{-1}y^{-1} = x$ and therefore $\{x,x^{-1}\}$ is an orbit of 
$\langle y\rangle$ which forces $x =x^{-1}\,,$ so that $x = 1\,;$ in particular, we get
$${\rm tr}_M (sy)\.{\rm tr}_{M^*} (sy) = {\rm tr}_S (sy) = 1
\eqno £A.10.9.$$
But, we clearly have
$$k\langle s \rangle = k\.{\rm id}_M + k\.s = k\.i +k\.i'
\eqno £A.10.10$$
for suitable mutually orthogonal idempotents $i$ and $i'$ of $S\,,$ and we choose the notation in such a way that
${\rm dim}\big(i(M)\big)\ge {\rm dim}\big(i'(M)\big)\,.$

\smallskip
Then, denoting by $\varphi_M\,,$ $\varphi_{i(M)}$ 
and $\varphi_{i(M')}$ the respective Brauer charac-ters of~$M\,,$~$i(M)$ and~$i'(M)\,,$ let us consider the  ordinary characters
$\chi_M\,,$ $\chi_{i(M)}$ and $\chi_{i'(M)}$ over~$\hat G\cap H$ which respectively lift the restrictions to $\hat G\cap H$  
of $\varphi_M\,,$  $\varphi_{i(M)}$ and $\varphi_{i'(M)}$ to the set of characters $\chi$ fulfilling $\chi (y) = \chi (y_{p'})$ for any  $y\in \hat G\cap H\,;$ consequently, we clearly have $\chi_M = \chi_{i(M)} + \chi_{i(M')}$ and moreover
$$ 1 =\big(\chi_{i(M)}(y) - \chi_{i'(M)}(y)\big)\big(\bar\chi_{i(M)}(y) - \bar\chi_{i'(M)}(y)\big)
\eqno £A.10.11\phantom{.}$$
for any $y\in \hat G\cap H$ (cf.~£A.10.9); in particular, the {\it norm\/} of $\chi_{i(M)} - \chi_{i'(M)}$ is equal to $1$ and, according to our choice  of notation, we still have 
$$1 = \chi_{i(M)}(1) - \chi_{i'(M)}(1)
\eqno £A.10.12\phantom{.}$$
hence, for a suitable linear character $\lambda$ of~$\hat G\cap H\,,$ 
 we get $\chi_{i(M)} = \lambda  + \chi_{i'(M)}$ or, equivalently,
$$\chi_M = \lambda + 2\.\chi_{i'(M)}
\eqno £A.10.13.$$

 \smallskip
 Now, it suffices to prove that $\lambda$ is the trivial character. Note that, denoting by $k'$ the subfield of $k$ generated
 by the $\ell\-$th roots of unity, we still can define a $k'^*\-$group $\hat K' = k'^*\times K$ as in~£A.11.1 below and, setting
 $S' = k'_*\hat K'\,,$ we have $S =  k\otimes_{k'} S'$ and $H$ is contained in $1\otimes S'\,,$ so that $i$ and~$i'$ 
 also belong to $1\otimes S'\,;$  hence, the values of
the  ordinary characters $\chi_M\,,$ $\chi_{i(M)}$ and~$\chi_{i(M')}$ are contained in the extension of $\Bbb Q$ by
 the $\ell\-$th roots of unity. Consequently, it suffices to prove that the restriction of $\lambda$ to a Sylow $\ell\-$subgroup
 $L$ of $\hat G\cap H$ is trivial.

 \smallskip
 But, it is well-known that for a maximal Abelian $k^*\-$subgroup $\hat A$ of $\hat K$ and  a $k^*\-$group homomorphism
$\zeta\,\colon \hat A\to k^*\,,$ denoting by $k_\zeta$ the corresponding $k_*\hat A\-$module, we have
 $$M \cong  {\rm Ind}_{\hat A}^{\hat K}(k_\zeta)
 \eqno £A.10.14;$$
 moreover, since $H$ acts over $\hat K\cong k^*\times K$ stabilizing $1\times K\,,$ it is easily checked that $L$ stabilizes a suitable choice of ${\rm Ker}(\zeta)\i 1\times K$ and therefore, choosing a complement $X$ of ${\rm Ker}(\zeta)$ in $K\,,$
it stabilizes  the basis $\{(1,x)\otimes 1\}_{x\in X}$ of $M\,;$ then, $L$ fixes $(1,1)\times 1$ and, for any
$L\-$orbit $O$ in $X-\{1\}\,,$ $\{(1,x)\otimes 1\}_{x\in O}$ and $\{(1,x^{-1})\otimes 1\}_{x\in O}$ are {\it different\/}
orbits of $L$ in this basis, since $\vert O\vert$ is odd; that is to say, the number of orbits of $L$ in this basis is odd.

\smallskip

In conclusion, since $L$ is an $\ell\-$group and $\ell\not= p\,,$ the multiplicity of $k$ in  $M$ considered as a $kL\-$module 
is an odd number; then, the restriction of equality~£A.10.13 to $L$ proves that the restriction of $\lambda$ to $L$ is trivial.
We are done.

\medskip
 £A.11.  Let $\ell$ be an odd prime number different from $p$ and $\hat K$  a  $k^*\-$group which is the central product of $k^*$ by an extra-special $\ell\-$group and, for our purposes, we also consider $\Bbb Z/\ell\Bbb Z$ as an {\it extra-special\/}
 $\ell\-$group.
 Denote by $\kappa$ the  {\it non-singular skew symmetric scalar  product\/} over the $k^*\-$quotient $K$ induced by the  commutator in $\hat K\,;$ thus, we have $\vert K\vert = \ell^{2n}$ and note that the case $n=0$ is not excluded.
 Then,  it is easily checked that $\hat K$ is isomorphic to $k^*\times K$ endowed with the product defined by
$$(\lambda, x)\.(\lambda', x') = (\lambda \lambda' \kappa (x,x')^{1\over 2}, xx')
\eqno £A.11.1,$$
for any $\lambda,\lambda'\in k^*$ and any $x,x'\in K\,,$ and with the group homomorphism $k^*\to k^*\times K$
mapping $\lambda\in k^*$ on $(\lambda,1)\,.$ It is quite clear that the corresponding symplectic group $Sp(K,\kappa)$ acts over this $k^*\-$group and we actually have
$${\rm Aut}_{k^*}(\hat K)\cong K\rtimes Sp(K,\kappa)
\eqno £A.11.2.$$

\medskip
£A.12. Moreover, it is well-known that $S= k_*\hat K$ is a simple $k\-$algebra and  $Sp(K,\kappa)$ clearly  acts over 
this $k\-$algebra stabilizing~$\hat K\,;$ thus,
since any central $k^*\-$extension of $Sp(K,\kappa)$ is trivial, this action can be lifted to a group homomorphism
$$Sp(K,\kappa) \too N_{S^*} (\hat K)
\eqno £A.12.1;$$
then, since $C_{S^*}(\hat K) = k^*\.{\rm id}_S\,,$ from isomorphism~£A.11.2 we easily get
$$N_{S^*} (\hat K) \cong \hat K\rtimes Sp(K,\kappa)
\eqno £A.12.2;$$
let us identify $Sp(K,\kappa)$ with its image in $N_{S^*} (\hat K)$ (for the choice of a $k^*\-$group isomorphism 
$\hat K\cong k^*\times K\,!$).

\medskip
£A.13. Identifying $\hat K$ with $k^*\times K\,,$ it is clear that a $p\-$subgroup $R$ of $Sp(K,\kappa)$ stabilizes the basis $\{(1,x)\}_{x\in K}$ of $S$
and therefore $S$ becomes a {\it Dade $R\-$algebra\/}; moreover, it is easily checked that the restriction  $\kappa_R$ 
of $\kappa$ to $C_K (R)$ remains a  {\it non-singular skew symmetric scalar  product\/} and therefore
$C_{\hat K}(R) \cong k^*\times C_K (R)$ is also the central product of $k^*$ by an extra-special $\ell\-$group. 
Then, it is clear that the Brauer homomorphism induces a $k\-$algebra isomorphism [11,~statement~2.8.4]
$$k_*C_{\hat K}(R)\cong S(R)
\eqno £A.13.1\phantom{.}$$
and it is easily checked that the action of $N_{Sp(K,\kappa)}(R)$ 
over $C_{\hat K}(R)$ is contained in the corresponding symplectic group $Sp\big(C_K (R),\kappa_R\big)\,;$
that is to say, the Brauer homomorphism can be extended to a group homomorphism
$${\rm Br}^*_R : N_{Sp(K,\kappa)}(R)\too Sp\big(C_K (R),\kappa_R\big)\i S(R)^*
\eqno £A.13.2\phantom{.}$$ 
such that the action of $x\in N_{Sp(K,\kappa)}(R)$ coincides with the conjugation by~${\rm Br}^*_R(x) $ on $S(R)\,;$
thus, choosing an element $a\in S^R$ lifting ${\rm Br}^*_R(x) $ and  an idempotent $i$ in the unique local point 
of $R$ on~$S\,,$ and denoting by $\skew4\widetilde{x}^R$ the image of $x$ in $F_S (R)$ and by 
$\overline{i x  a^{-1} i}^S$ the image of 
the product $i x  a^{-1} i$
in the quotient
$$N_{(iSi)^*}(Ri)\Big/\Big(i + J\big((iSi)^R\big)\Big)
\eqno £A.13.3,$$
the pair $(\skew4\widetilde{x}^R,\overline{i x a^{-1} i}^S\,)$ is an element of $\hat F_S (R)$ [11,~Proposition~6.10].

\bigskip
\noindent
{\bf Lemma~£A.14.} {\it With the notation above, denote by  $\,{}^2\omega$ the unique odd-polarization such that diagram~{\rm £4.3.1} is commutative and let $R$ be a $p\-$subgroup  of~$Sp (K,\kappa)\,.$ For any 
$x\in N_{Sp(K,\kappa)}(R)$ of odd order, choosing an element $a\in S^R$ lifting ${\rm Br}^*_R(x) $ and  an idempotent $i$ in 
the unique local point of $R$ on~$S\,,$ and denoting by $\skew4\widetilde{x}^R$ the image of $x$ in $F_S (R)$ and by 
$\overline{i x  a^{-1} i}^S$ the image of  the product $i x  a^{-1} i$ in the quotient
$$N_{(iSi)^*}(Ri)\Big/\Big(i + J\big((iSi)^R\big)\Big)
\eqno £A.14.1,$$
 we have 
$${}^2 \omega_{(R,S)}(\skew4\widetilde{x}^R,\overline{i x  a^{-1} i}^S) = 1
\eqno £A.14.2.$$\/}

\par
\noindent
{\bf Proof:} Arguing by induction on $\vert R\vert\,,$ set $Z = \Omega_1\big(Z(R)\big)\,,$ $\bar R = R/Z$ and
 $T = {\rm Res}_Z^R(S)\,;$ it follows   from Theorem~£A.3 above that we have the commutative~diagram 
$$\matrix{\Bbb O^2\big(\hat F_S (R)_Z\big)Ê\hskip-5pt&\buildrel \Delta_{R,S,Z}
\over{\hbox to 35pt{\rightarrowfill}} &\hskip-5pt \Bbb O^2\big(\hat F_T (Z)\big)
\,\hat\times\, \Bbb O^2\big(\hat F_{S(Z)}(\bar R)\big)\cr 
{\scriptstyle {}^2\omega_{(R,S)}} \searrow \hskip-20pt&\phantom{\bigg\uparrow} &
\hskip-20pt\swarrow\hskip5pt \scriptstyle
{{}^2\omega_{(Z,T)}\,\hat\times\,{}^2\omega_{(\bar R,S(Z))}}\cr  &\ k^* &\cr}
\eqno £A.14.3.$$
But,  choosing an element $c\in S^Z$ lifting $\bar x = {\rm Br}^*_Z(x)$ and  an idempotent $j$ in the unique local point 
of $Z$ on~$S$ fulfilling $ji = j = ij\,,$ and setting $\bar a ={\rm Br}_Z(a)$ and $\bar\imath ={\rm Br}_Z(i)\,,$ 
it is easily checked  from [15,~Proposition~9.11] that we have 
$$\Delta_{R,S,Z}(\skew4\widetilde{x}^R,\overline{i x  a^{-1} i}^S) = (\skew4\widetilde{x}^Z,\overline{jxc^{-1} j}^T)\.
(\skew4\widetilde{\bar x}^{\bar R},\overline{\bar\imath\bar x \bar a^{-1}\bar\imath}^{S(Z)})
\eqno £A.14.4;$$
hence, since we  clearly have
$${\rm Br}^*_{\bar R}\big({\rm Br}^*_Z(x)\big) = {\rm Br}^*_R (x)
\eqno £A.14.5,$$
if $Z\not= R$ then from the induction hypothesis we get
$${}^2 \omega_{(Z,T)}(\skew4\widetilde{x}^Z,\overline{jxc^{-1} j}^T) = 1=
{}^2 \omega_{(\bar R,S(Z))}(\skew4\widetilde{\bar x}^{\bar R},\overline{\bar\imath\bar x \bar a^{-1}\bar\imath}^{S(Z)}) 
\eqno £A.14.6.$$
Now, equality~£A.14.2 follows from the commutativity of~diagram~£A.14.3.

\smallskip
From now on, we assume that $R$ is $p\-$elementary Abelian. Arguing by induction on $\vert K\vert\,,$ if $K$ decomposes on
a {\it direct  orthogonal\/} sum of two $R\.\langle x\rangle\-$stable nontrivial subspaces
$$K = K'\,\bot\, K''
\eqno £A.14.7\phantom{.}$$
then $\hat K$ is the central product of
 the converse images $\hat K'$ of $K'$ and $\hat K''$ of~$K''\,,$ and, setting 
 $$S' =k_*\hat K'\qq S'' =k_*\hat K''
 \eqno £A.14.8,$$
$S'$ and $S''$ are also {\it Dade $R\-$algebras\/} and we have $S \cong S'\otimes_k S''\,;$ in particular, it follows   from Theorem~£A.3 above that we have the commutative~diagram 
$$\matrix{\Bbb O^2\big(\hat F_{S'} (R)\big)\,\,\hat\cap\,\,\Bbb O^2\big(\hat F_{S''} (R)\big)
&\buildrel \hat\nu_{R,S',S''}\over{\hbox to 30pt{\rightarrowfill}} &\Bbb O^2\big(\hat F_{S} (R)\big)\cr
{\scriptstyle {{}^2\omega_{(R,S')}\,\hat\times\,{}^2\omega_{(R,S'')}}}\; \searrow\hskip-15pt
&\phantom{\bigg\uparrow} & \swarrow\hskip5pt
\scriptstyle {{}^2\omega_{(R,S)}}\cr &\ k^* &\cr}
\eqno £A.14.9.$$

\smallskip
But, denoting by $\kappa'$ and $\kappa''$ the respective restrictons of $\kappa$ to $K'$ and~$K''\,,$ it is clear that 
$x$ is the image of $x'\otimes x''$ for suitable elements $x'\in Sp(K',\kappa')$ and $x''\in Sp (K'',\kappa'')$
normalizing the respective images $R'\i Sp(K',\kappa')$ and~$R''\i Sp (K'',\kappa'')$ of $R\,.$ Moreover,
 choosing  elements $a'\in S'^{R'}$ and  $a''\in S''^{R''}$ respectively lifting ${\rm Br}^*_{R'} (x') $ and 
 ${\rm Br}^*_{R''} (x'') \,,$ and  idempotents $i'$ and $i''$ in the respective unique local points of $R$ on~$S'$
 and on $S''\,,$ we clearly may choose the element $a$ equal to the image of $a'\otimes a''$ and the idempotent $i$ in an orthogonal decomposition of the image of $i'\otimes i''\,,$ so that ${\rm Br}_{R}(x)$ is equal to the corresponding image of
${\rm Br}_{R'}(x')\otimes {\rm Br}_{R''}(x'')$ {\it via\/} the isomorphism
$$S'(R')\otimes_k S''(R'')\cong S(R)
\eqno £A.14.10.$$
Then, it easily follows from [15,~9.15] that we have
$$\eqalign{\hat\nu_{R,S',S''}\big((\skew4\widetilde{x'}^{R'},\overline{i 'x'  a'^{-1} i'}^{S'})\.
(\skew4\widetilde{x''}^{R''},\overline{i'' x''  a''^{-1} i''}^{S''})&\big)\cr
 = (&\skew4\widetilde{x}^R,\overline{i x  a^{-1} i}^S)\cr}
\eqno £A.14.11.$$
Now, since the induction hypothesis implies that
$${}^2 \omega_{(R',S')}(\skew4\widetilde{x'}^{R'},\overline{i 'x'  a'^{-1} i'}^{S'}) = 1 =
{}^2 \omega_{(R'',S'')}(\skew4\widetilde{x''}^{R''},\overline{i'' x''  a''^{-1} i''}^{S''})
\eqno £A.14.12,$$
 equality~£A.14.2 follows from the commutativity of~diagram~£A.14.9.

 \smallskip
Thus, we may assume that $K$ does not admit a decomposition on a {\it direct  orthogonal\/} sum of two 
$R\.\langle x\rangle\-$stable nontrivial subspaces. Denoting by $\Bbb F$ the field of cardinal $\ell\,,$
if $L$ is a simple $\Bbb F R\-$submodule of $K$ then the restriction of $\kappa$ to $L$ is either {\it non-singular\/}
or, denoting by $L^{\bot}$ the {\it orthogonal\/} space of~$L\,,$ $L^\bot$ contains $L$ and we have a canonical isomorphism $K/L^\bot\cong L^*\,,$ so that, if~$L'$ is an $\Bbb F R\-$complement of $L^{\bot}$ in $K\,,$ the restriction of $\kappa$ to 
$L\oplus L'$ is {\it non-singular\/}. Consequently, the dimensions of all the simple $\Bbb F R\-$submodules of~$K$ 
have the same {\it parity.\/}

\smallskip
Firstly assume that the dimensions of all the simple $\Bbb F R\-$submodules $L$ of~$K$ are~{\it odd\/}; in this case, $L^*$ is also a $\Bbb F R\-$submodule of $K$ {\it not\/} isomorphic to~$L\,.$ Thus, since $\vert\langle x\rangle\vert$ is odd, 
the group $\langle x\rangle$  has exactly two orbits in the set of {\it isotypic components\/} of the $\Bbb F R\-$module~$K$
and then, denoting by $A$ and~$B$ the sums  of {\it isotypic components\/} in each $\langle x\rangle\-$orbit, 
$A$ and $B$ are maximal {\it totally singular\/} subspaces fulfilling $K = A\oplus B\,.$ Hence, the converse images
$\hat A$ of $A$ and $\hat B$ of $B$ in $\hat K$ are maximal Abelian subgroups and  it is well-known that, for  a $k^*\-$group homomorphism $\zeta\,\colon \hat A\to k^*$ that we may choose $R\.\langle x\rangle\-$stable (cf.~£A.11.1), we have
 $$M \cong  {\rm Ind}_{\hat A}^{\hat K}(k_\zeta)
 \eqno £A.14.13\phantom{.} $$
where $k_\zeta$ denotes the corresponding $k_*\hat A\-$module.

\smallskip
In this situation, the group $R\.\langle x\rangle$  stabilizes  the basis $\{(1,y)\otimes 1\}_{y\in B}$ of~$M\,,$
so that the {\it Dade $R\-$algebra\/} $S$ is {\it similar\/} to $k\,;$ in particular, identifying $S$ with the {\it induced 
$\hat K\-$interior algebra\/} ${\rm Ind}_{\hat A}^{\hat K}(k_\zeta)$ (cf.~£2.12) where~$k_\zeta$ still denotes the corresponding $\hat A\-$interior algebra, the primitive idempotent 
$$ i = (1,1)\otimes 1\otimes (1,1)
\eqno £A.14.14\phantom{.}$$ 
actually belongs to the unique local point of $R$ on $S\,;$ now,  $x$ and ${\rm Br}^*_R(x)$ respectively centralize~$i$ and
${\rm Br}_R(i)\,,$ and, with the notation above, it is easily checked that $\overline{i x  a^{-1} i}^S = \overline{i}^S $
which proves equality~£A.14.2 in this case.

\smallskip
Finally assume that the dimensions of all the simple $\Bbb F R\-$submodules $L$ of~$K$ are~{\it even\/};
in this case, the image of $\Bbb F R$ in ${\rm End}_{\Bbb F}(L)$ is an extension $\Bbb F_L$ of $\Bbb F$ of even degree
and therefore it contains a {\it primitive fourth\/} root $\tau_L$ of unity;\break
\eject
\noindent 
moreover, since $\vert\langle x\rangle\vert$ is odd, the stabilizer in $\langle x\rangle$ of  the {\it isotypic component\/} 
containing $L$ acts on $\Bbb F_L$ fixing $\tau_L\,.$ Consequently, considering all the orbits of~$\langle x\rangle\,,$ 
we get a {\it self-adjoin\/}  endomorphism $\tau$ of $K$
which centralizes $R\.\langle x\rangle$ and fulfills $\tau^2 = -{\rm id}_K\,.$

\smallskip
At this point, we consider the central product $\hat K\,\hat\times\, \hat K\,,$ and in the $k^*\-$quo-tient $K\times K$ we set
$$A = \{(y,\tau (y))\}_{y\in K}\qq B = \{(-y,\tau (y))\}_{y\in K}
\eqno £A.14.15;$$
as above, we have $K = A\oplus B\,,$ $A$ is {\it totally singular\/}  since
$$\eqalign{(\kappa\times \kappa)\big((y,\tau (y)),(y',\tau (y'))\big) &= \kappa(y,y') \kappa\big(\tau(y),\tau(y')\big)\cr
&= \kappa(y,y') \kappa\big(\tau^2(y),y'\big)\cr
& = \kappa(y,y')\kappa(y,y')^{-1} = 1\cr}
\eqno £A.14.16\phantom{.}$$
for any $y\in K$ and, similarly,  $B$ is {\it totally singular\/} too. Once again, the converse images
$\hat A$ of $A$ and $\hat B$ of $B$ in $\hat K\,\hat\times\, \hat K$ are maximal Abelian subgroups;
hence, the argument above applied to the $p\-$subgroup $\Delta(R) = \{u\otimes u\}_{u\in R}$ and to the element
$x\otimes x$ of $Sp(K\times K,\kappa\times \kappa)$ proves that
$${}^2\omega_{(\Delta(R),S\otimes_k S)}\big(\widetilde{x\otimes x}^{\Delta(R)},
\overline{j(x\otimes x)(a\otimes a)^{-1}j}^{S\otimes_k S}\big) = 1
\eqno £A.14.17$$
for the choice of an idempotent $j$  in the unique local point of $\Delta(R)$ on~$S\otimes_k S\,.$

\smallskip
Consequently, since it follows again  from Theorem~£A.3 that we have the commutative~diagram 
$$\matrix{\Bbb O^2\big(\hat F_{S} (R)\big) \star \Bbb O^2\big(\hat F_{S} (R)\big)
&\buildrel \hat\nu_{R,S,S}\over{\hbox to 30pt{\rightarrowfill}} &\Bbb O^2\big(\hat F_{S\otimes_k S} (R)\big)\cr
{\scriptstyle {{}^2\omega_{(R,S)}\,\hat\times\,{}^2\omega_{(R,S)}}}\; \searrow\hskip-15pt
&\phantom{\bigg\uparrow} & \swarrow\hskip5pt
\scriptstyle {{}^2\omega_{(R,S\otimes_k S)}}\cr &\ k^* &\cr}
\eqno £A.14.18\phantom{.}$$ 
and since we clearly have
$$\eqalign{\hat\nu_{R,S,S}\big((\skew4\widetilde{x}^R,\overline{i x  a^{-1} i}^S)&\.
(\skew4\widetilde{x}^R,\overline{i x  a^{-1} i}^S)\big)\cr 
&= \big(\widetilde{x\otimes x}^{\Delta(R)},
\overline{j(x\otimes x)(a\otimes a)^{-1}j}^{S\otimes_k S}\big)\cr}
\eqno £A.14.19,$$
from equality~£A.14.17 we actually get 
$$\big({}^2 \omega_{(R,S)}(\skew4\widetilde{x}^R,\overline{i x  a^{-1} i}^S)\big)^2 = 1
\eqno £A.14.20$$
which forces ${}^2 \omega_{(R,S)}(\skew4\widetilde{x}^R,\overline{i x  a^{-1} i}^S)=1$ since $x$ has odd order.
We are done.

\bigskip\bigskip
\centerline{\large References}

\bigskip\noindent
[1]\phantom{.} Jon Alperin,  {\it Weights for finite groups\/}, in Proc. Symp.
Pure Math. 47 (1987) 369-379, Amer. Math. Soc., Providence.
\smallskip\noindent
[2]\phantom{.} Michel Brou\'e and Llu\'\i s Puig, {\it Characters and Local
Structure in $G\-$al-gebras,\/} Journal of Algebra, 63(1980), 306-317.
\smallskip\noindent
[3]\phantom{.} Walter Feit and John Thompson,  {\it Solvability of groups of
odd order\/}, Pacific Journal of Math., 13(1963), 775-1029. 
\smallskip\noindent
[4]\phantom{.} Daniel Gorenstein, {\it ``Finite groups''\/} Harper's Series,
1968, Harper and Row.
\smallskip\noindent
[5]\phantom{.} Marty Isaacs,  {\it Characters of Solvable and Symplectic groups\/}, 
Amer. Journal of Math., 95(1972), 594-635. 
\smallskip\noindent
[6]\phantom{.} Burkhard K\"ulshammer and Llu\'\i s Puig, {\it Extensions of
nilpotent blocks}, Inventiones math., 102(1990), 17-71.
\smallskip\noindent
[7]\phantom{.}  Gabriel Navarro, {\it Weights, vertices and a correspondence of characters in groups of odd
order\/}, Math. Z., 212(1993), 535-544
\smallskip\noindent
[8]\phantom{.} Tetsuro Okuyama, {\it Vertices of irreducible modules of $p\-$solvable groups\/}, preprint, 1981
\smallskip\noindent
[9]\phantom{.} Llu\'\i s Puig, {\it Pointed groups and  construction of
characters}, Math. Z. 176 (1981), 265-292. 
\smallskip\noindent
[10]\phantom{.} Llu\'\i s Puig, {\it Local fusions in block source algebras\/},
Journal of Algebra, 104(1986), 358-369. 
\smallskip\noindent
[11]\phantom{.} Llu\'\i s Puig, {\it Pointed groups and  construction of
modules}, Journal of Algebra, 116(1988), 7-129.
\smallskip\noindent
[12]\phantom{.} Llu\'\i s Puig, {\it Nilpotent blocks and their source
algebras}, Inventiones math., 93(1988), 77-116.
\smallskip\noindent
[13]\phantom{.} Llu\'\i s Puig, {\it Affirmative answer to a question of Feit},
Journal of Algebra, 131(1990), 513-526.
\smallskip\noindent
[14]\phantom{.} Llu\'\i s Puig, {\it Source algebras of $p\-$central group 
extensions\/}, Journal of Algebra, 235(2001), 359-398.
\smallskip\noindent
[15]\phantom{.} Llu\'\i s Puig, ``{\it Frobenius categories versus Brauer blocks\/}'',
Progress in Mahtematics, 274(2009), Bikh\"auser, Basel
\smallskip\noindent
[16]\phantom{.} Llu\'\i s Puig, {\it Block Source Algebras  in p-Solvable
Groups}, Michigan Math. Journal, 58(2009), 323-338
\smallskip\noindent
[17]\phantom{.} Jacques T\'evenaz, ``$G$-{\it Algebras and Modular
Representation Theory\/}'', Oxford Math. Mon., 1995, Clarendon Press, Oxford

\end